\pdfoutput=1
\date{} 
\title{The scaling limits of the Minimal Spanning Tree and Invasion Percolation in the plane}

\author{Christophe Garban \and G\'abor Pete \and Oded Schramm}

\documentclass[12pt,preprint, letterpaper]{article}

\usepackage[vmargin={2.5cm, 3cm},hmargin=2.6cm]{geometry}

\usepackage[colorlinks=false,pdfborderstyle={/W 1},pdfpagelabels]{hyperref}

\usepackage{latexsym} 
\usepackage{amsmath}
\usepackage{amssymb}
\usepackage{amsthm}
\usepackage{mathrsfs}
\usepackage{amsfonts}
\usepackage{graphicx}
\usepackage{url}
\input labelfig.tex  

\usepackage[usenames,dvipsnames]{color}
\usepackage{bbm,verbatim}

\definecolor{MyMagenta}{rgb}{0.9,0,0.9}

\usepackage{tikz}

\usepackage{multirow} 
\usepackage{array}
\newcolumntype{M}[1]{>{\centering}m{#1}}

\def \proof {{ \medbreak \noindent {\bf Proof. }}}
\def\proofof#1{{ \medbreak \noindent {\bf Proof of #1.}}}

\numberwithin{equation}{section}
\numberwithin{figure}{section}
\newtheorem{theorem}{Theorem}
\numberwithin{theorem}{section}
\newtheorem{corollary}[theorem]{Corollary}
\newtheorem{lemma}[theorem]{Lemma}
\newtheorem{proposition}[theorem]{Proposition}
\newtheorem{conjecture}[theorem]{Conjecture}

\newtheorem{question}[theorem]{Question}

\newtheorem{definition}[theorem]{Definition}
\theoremstyle{remark}\newtheorem{remark}[theorem]{Remark}
\def\eqref#1{(\ref{#1})}
\let\qqed=\qed

\def\bl{\begin{lemma}}
\def\el{\end{lemma}}
\def\bth{\begin{theorem}}
\def\eth{\end{theorem}}
\def\bc{\begin{corollary}}
\def\ec{\end{corollary}}
\def\bcj{\begin{conjecture}}
\def\ecj{\end{conjecture}}
\def\bpr{\begin{proposition}}
\def\epr{\end{proposition}}
\def\bde{\begin{definition}}
\def\ede{\end{definition}}
\newcommand{\be}{\begin{eqnarray}}
\newcommand{\ee}{\end{eqnarray}}
\newcommand{\bes}{\begin{eqnarray*}}
\newcommand{\ees}{\end{eqnarray*}}
\def\bq{\begin{question}}
\def\eq{\end{question}}

\def\bi{\begin{itemize}}
\def\ei{\end{itemize}}
\def\bnum{\begin{enumerate}}
\def\enum{\end{enumerate}}

\def\QED{\qqed\medskip}
\let\qed=\QED

\newcommand{\R}{\mathbb{R}}

\newcommand{\C}{\mathbb{C}}
\newcommand{\Z}{\mathbb{Z}}
\newcommand{\N}{\mathbb{N}}

\def\RR{\mathcal{R}}

\def\diam{\mathrm{diam}}

\def\dist{\mathrm{dist}}

\def \eps {\epsilon}
\def \P {{\bf P}}
\def \E {{\bf E}}
\def\md{\mid}
\def\Bb#1#2{{\def\md{\bigm| }#1\bigl[#2\bigr]}}

\def\Pb{\Bb\P}
\def\Eb{\Bb\E}

\def \p {{\partial}}

\def\calC{\mathcal{C}}

\def\1{\mathbf{1}}

\def\lora{\longrightarrow}

\def\Set{\mathcal{S}}
\def\For{\mathcal{F}}
\def\bigtimes{\mathop{\mathsf{X}}}

\def\MST{\mathsf{MST}}

\def\MSF{\mathsf{MSF}}
\def\Trunk{\mathsf{Trunk}}
\def\G{\mathcal{G}}

\def\vl{{\bar{\lambda}}}
\def\IP{\mathsf{InvPerc}}

\def\A{\mathcal{A}}

\def\Tg{\mathsf{T}} 
\def\Td{\mathsf{T}^*}
\def\Tor{\mathbb{T}^2}
\def\boxup{\boxdot}
\def\Quad{Q} 
\def\QUAD{\mathcal{Q}} 
\def\HH{\mathscr{H}}  
\def\T{\mathcal{T}}  

\def\PPP{\mathsf{PPP}}

\def\Net{\mathsf{N}}
\def\EnNet{\mathsf{EN}}
\def\meso{\textrm{-}\mathsf{meso}}

\def\Piv{\mathcal{P}}
\def\PivO{\Piv_{\mathrm{open}}}
\def\PivC{\Piv_{\mathrm{closed}}}

\def\Pimp{\mathscr{P}}
\def\PimpO{\Pimp_{\mathrm{open}}}
\def\PimpC{\Pimp_{\mathrm{closed}}}

\newcommand{\red}[1]{\textcolor{red}{#1}}
\newcommand{\blue}[1]{\textcolor{blue}{#1}}

\begin{document}
\maketitle

\begin{abstract}
We prove that the Minimal Spanning Tree and the Invasion Percolation Tree on a version of the triangular lattice in the complex plane have unique scaling limits, which are invariant under rotations, scalings, and, in the case of the $\MST$, also under translations. However, they are not expected to be conformally invariant. We also prove some geometric properties of the limiting $\MST$. The topology of convergence is the space of spanning trees introduced by Aizenman, Burchard, Newman \& Wilson (1999), and the proof relies on the existence and conformal covariance of the scaling limit of the near-critical percolation ensemble, established in our earlier works.
\end{abstract}

\vskip 0.8 cm

\begin{figure}[htbp]
\centerline{
\AffixLabels{
\includegraphics[width=0.45 \textwidth]{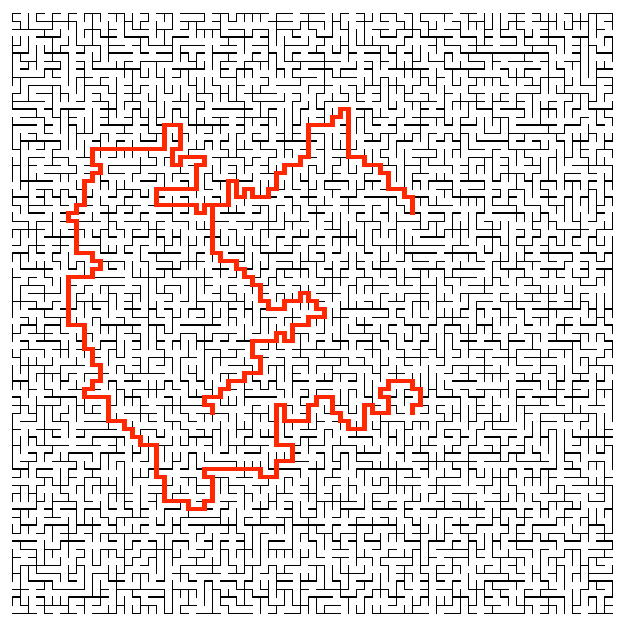}
\hskip 1 cm
\includegraphics[width=0.45 \textwidth]{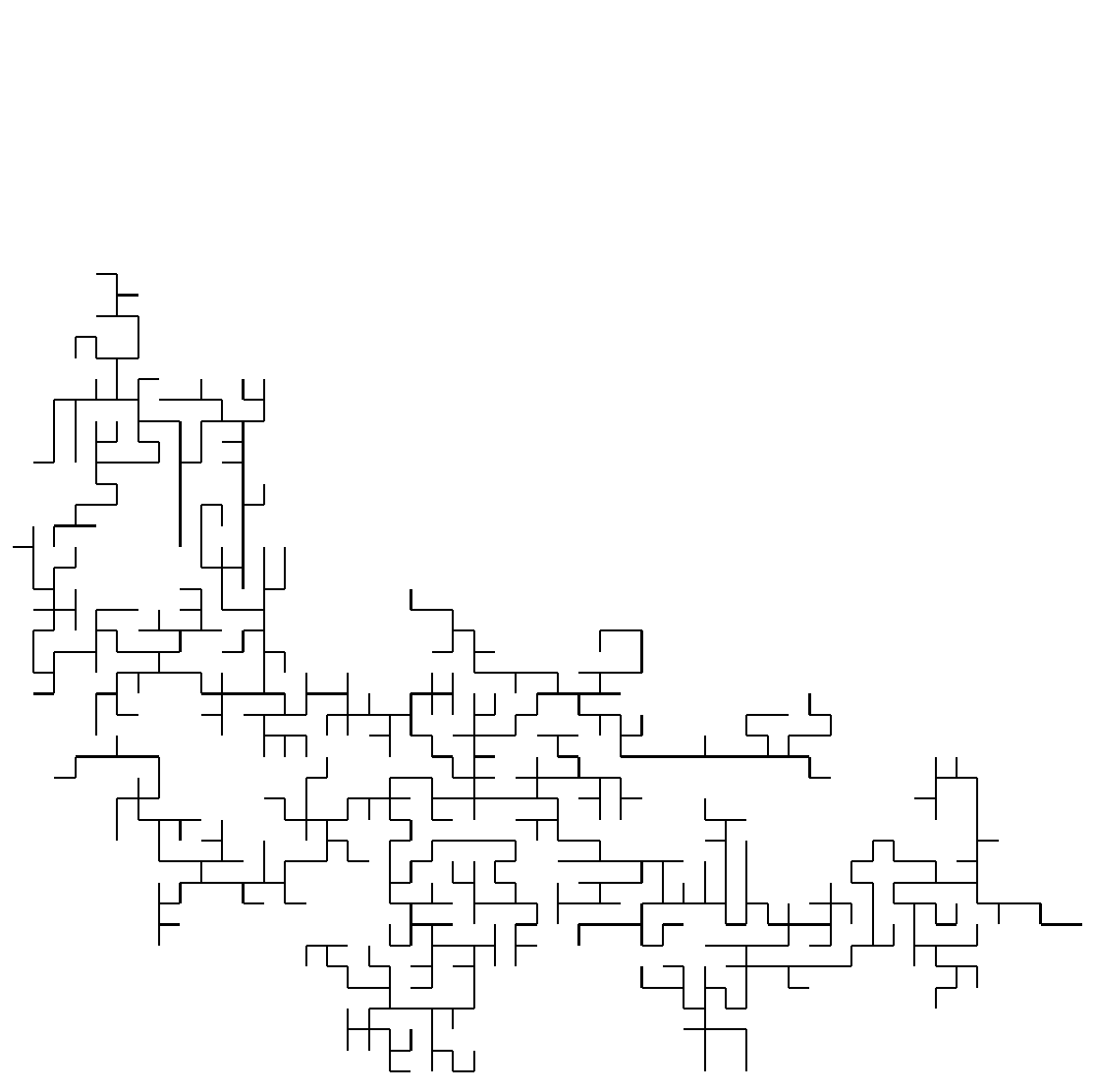}
}}
\vskip 0.3 cm
{The $\MST$ in a box, and $\IP$ started from the midpoint of the left boundary of the box until reaching the right boundary, on $\Z^2$.}
\label{f.first}
\end{figure}

\eject

\tableofcontents

\section{Introduction} 

The Minimal Spanning Tree of weighted graphs is a classical combinatorial object \cite{Bor,Kruskal,GraHell}, and is also very interesting from the viewpoint of probability theory and statistical physics:  when the weights on the edges of a graph are chosen at random, using i.i.d.~variables, then the resulting random tree turns out to be closely related to the near-critical regime of Bernoulli bond percolation on that graph.

In Bernoulli bond percolation at density $p\in[0,1]$, each edge of the graph is kept open with probability $p$ or becomes closed with probability $1-p$, independently, and then one looks at the connected open components, called clusters. In site percolation, the vertices are chosen to be open or closed instead of the edges. These are among the most important spatial stochastic processes, due to their simultaneous simplicity and richness \cite{BroHa,KestenBook,Grimmett}. The main interest is in the phase transition near the critical density $p_c$, below which all clusters are small, above which a cluster (sometimes clusters) of positive density emerge. The theory of critical percolation in the plane has seen a lot of progress lately, starting with Smirnov's proof of conformal invariance of crossing probabilities for site percolation on the triangular lattice \cite{Smirnov}, and with the introduction of the Stochastic Loewner Evolution \cite{Schramm} that describes the conformally invariant curves that are the scaling limits of interfaces between open and closed clusters. These SLE curves can be used to understand critical percolation in depth \cite{WWperc}, including the computation of critical exponents that had been predicted by physicists using non-rigorous conformal field theory techniques. 

Beyond the static critical system, it is natural to consider dynamical versions: first, to slowly change $p$ near $p_c$ and observe how the phase transition exactly takes place --- called near-critical percolation; second, to apply a stationary dynamics and observe how the critical system is changing in time --- called dynamical percolation. Indeed, by ``perturbing''  critical percolation, the static results of the previous paragraph have also given way to an exhaustive study of dynamical and near-critical percolation \cite{SchSt, GPS1,HammPS, GPS2a,DPSL}; see also the surveys \cite{Jeff,Buzios}. In particular, in  \cite{GPS2a,DPSL} we have proved the existence and conformal covariance of the scaling limit of the near-critical percolation ensemble, w.r.t.~the quad-crossing topology introduced in~\cite{SSblacknoise}. Very roughly, this near-critical scaling limit is constructed from the critical scaling limit, plus independent randomness that governs how macroscopic clusters merge as we raise $p$. 

It turns out that the macroscopic structure of the Minimal Spanning Tree ($\MST$) and the Invasion Percolation Tree ($\IP$) can also be described based on this merging process. Thus, building on \cite{GPS2a,DPSL}, in the present paper we prove the existence and some conformal properties of the scaling limits of $\MST$ and $\IP$ on the triangular lattice, in the space of essential spanning forests introduced in \cite{ABNW}. In that paper, tightness results were proved, implying that subsequential scaling limits of the Minimal and Uniform Spanning Trees in the plane exist. Our proof of the uniqueness of the scaling limit has the important implication that the conjectural universality of critical percolation implies universality for many processes related to the near-critical ensemble, including $\MST$ and $\IP$. That this program of describing near-critical objects from the critical scaling limit may have a chance to work was suggested in \cite{CFNb}. Another motivation for our work is that it leads to interesting new objects: these two scaling limits are invariant under rotations and scalings, but, conjecturally, not under general conformal maps. Furthermore, the methods developed to establish these scaling limits also give information about the large-scale geometry of the discrete trees.

\subsection{The Minimal Spanning Tree $\MST$}\label{ss.introMST}

For each edge of a finite graph, $e\in E(G)$, let $U(e)$ be an independent Unif$[0,1]$ label. The {\bf Minimal Spanning Tree}, denoted by $\MST$, is the spanning tree $T$ for which $\sum_{e\in T} U(e)$ is minimal. This is well-known to be the same as the union of lowest level paths between all pairs of vertices (i.e., the path between the two points for which the maximum label on the path is minimal). One can also use the so-called {\bf reversed Kruskal algorithm} to construct $\MST$: delete from each cycle the edge with the highest label $U$. This algorithm also shows that $\MST$ depends only on the ordering of the labels, not on the values themselves. Moreover, this algorithm also makes sense on any infinite graph, and produces what in general is called the Free Minimal Spanning Forest ($\mathsf{FMSF}$) of the infinite graph. The Wired Minimal Spanning Forest ($\mathsf{WMSF}$) is the one when we also remove the edge with the highest label (if such edge exists) from each cycle that ``goes through infinity'', i.e., which is the union of two disjoint infinite simple paths starting from a vertex. For the case of Euclidean planar lattices, these two measures on spanning forests are known to be the same, again denoted by $\MST$, and it almost surely consists of a single tree \cite{AleMol}. This measure can also be obtained as a {\bf thermodynamical limit}: take any exhaustion by finite subgraphs $G_n(V_n,E_n)$, introduce a boundary condition by identifying some of the vertices on the boundary of $G_n$ (i.e., elements of $V_n$ that have neighbors in $G$ outside of $V_n$), and then take the weak limit. On a general infinite graph, when no identifications are made in the boundary, one gets the $\mathsf{FMSF}$, and when all vertices are glued into a single vertex, one gets the $\mathsf{WMSF}$. Studying these measures has a rich history on $\Z^d$, on point processes in $\R^d$, and on general transitive graphs; see \cite{Alexander, Penrose, Yukich, objective, LPS, Timar, ChatSen, TimarIndist, NTW, LPbook} and the references therein.

One can use the same Unif$[0,1]$ labels that defined the $\MST$ to obtain a coupling of percolation for all densities $p\in[0,1]$: an edge is ``open at level $p$'' if $U(e)\leq p$. This way we get a {\bf coupling} between the $\MST$ and the {\bf percolation ensemble}. Moreover, as we explain in the next paragraph, 
the macroscopic structure of the $\MST$ is basically determined by the labels in the near-critical regime of percolation, and hence one may hope that the scaling limit of the $\MST$ is determined by the scaling limit of the near-critical ensemble.

\begin{figure}[htbp]
\centerline{
\AffixLabels{
\includegraphics[width=0.45 \textwidth]{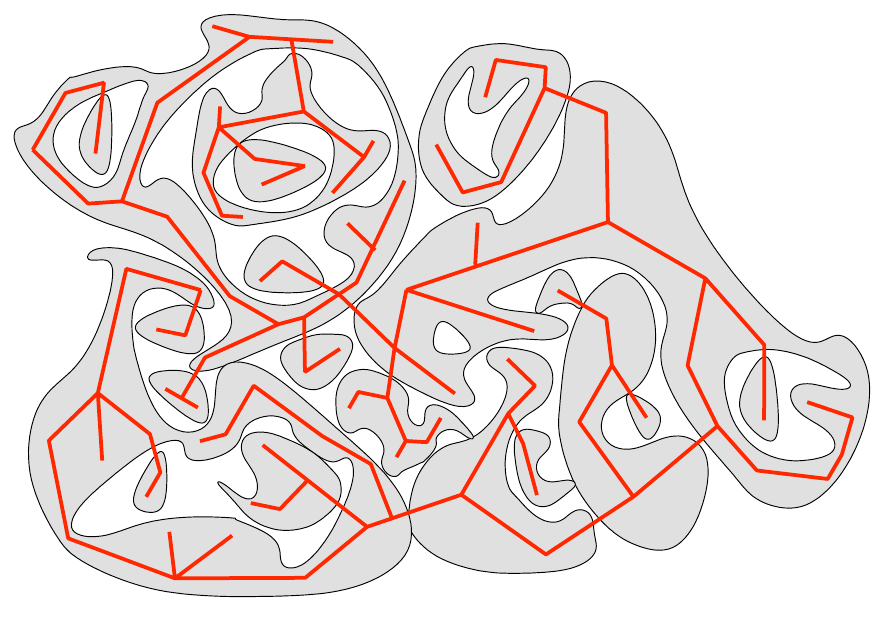}
}}
\caption{The $\MST$ connects the percolation $p$-clusters without creating cycles, yielding the cluster-tree $\MST^p$.}
\label{f.clusterTree}
\end{figure}

Consider the {\bf $p$-clusters} (i.e., open components at level $p$) in the percolation ensemble on some large finite graph. Contract each component into a single vertex, keeping the edges (together with their labels) between the clusters, resulting in the ``cluster graph''. It is easy to verify that making these contractions on the $\MST$ we get exactly the $\MST$ on the cluster graph. We denote this {\bf cluster tree} by $\MST^p$. See Figure~\ref{f.clusterTree}. Now assume that $p_1$ is small enough so that even the largest $p_1$-clusters are of small macroscopic size --- then the tree $\MST^{p_1}$ will tell us the macroscopic structure of $\MST$. On the other hand, if $p_2>p_1$ is large enough, then most sites are in just one giant $p_2$-cluster. Note that, for any $p>p_1$, we get the tree $\MST^{p}$ from $\MST^{p_1}$ by contracting the edges with labels in $(p_1,p]$. Thus, if we have the collection of {\it all} the $p$-clusters for all $p\in(p_1,p_2)$, then by following how they merge as we are raising $p$, we can reconstruct the tree $\MST^{p_1}$. Now, one may hope that in order to tell the macroscopic structure of $\MST^{p_1}$, it is enough to know only the {\it macroscopic} $p$-clusters for all $p\in(p_1,p_2)$ and follow how {\it those} merge. The {\bf near-critical window} of percolation is exactly the window $(p_1,p_2)$ in which the above phase transition of the cluster sizes takes place, and the scaling limit of the near-critical ensemble is exactly the object that describes the macroscopic $p$-clusters in this window. Therefore, the above hope has the interpretation that the scaling limit of the near-critical ensemble should describe the scaling limit of the $\MST$. This, of course, raises several questions: May the dust of microscopic $p$-clusters condensate into a new macroscopic $p'$-cluster at some $p'>p$, ruining the strategy of ``following how macroscopic clusters merge''? Could $\MST^{p_1}$ go through microscopic $p_1$-clusters in a way that significantly influences its macroscopic structure? 

Our work addresses these questions in the case of {\bf planar lattices}. The near-critical window for Bernoulli$(p)$ percolation on the triangular lattice $\eta\Tg$ or the square lattice $\eta\Z^2$ with mesh $\eta>0$ is given by
\be\label{e.nc}
p=1/2+\lambda r(\eta)\text{ with }\lambda\in(-\infty,\infty)\text{ fixed and }\eta\to 0\,,
\ee
where $r(\eta)=\eta^{2}/\alpha_4(\eta,1)$, with $\alpha_4(\eta,1)$ being the alternating 4-arm probability of critical percolation \cite{WWperc}. It was proved on $\eta\Tg$ using SLE$_6$ computations \cite{SW} that $r(\eta)=\eta^{3/4+o(1)}$. As shown in \cite{Kesten}, for $\lambda\ll -1$ we are at the subcritical end of the near-critical window, for $\lambda\gg 1$ we are at the supercritical end, and for any fixed $\lambda\in\R$, box-crossing probabilities are comparable to the critical case (just they are close to 0 for $\lambda\ll -1$, and close to 1 for  $\lambda\gg 1$). That is,~(\ref{e.nc}) is indeed the near-critical window. Then it was proved in \cite{GPS2a,DPSL} that for any $\lambda\in\R$ there is a unique scaling limit as $\eta\to 0$; moreover, the entire coupled percolation ensemble, viewed near the critical point via the parametrization~(\ref{e.nc}), where all the macroscopic changes happen, has a scaling limit as a Markov process in $\lambda\in\R$. It is important to keep in mind that even for any given $\lambda\not=0$, this scaling limit is an interesting new object, known to be different from the critical scaling limit: the interfaces are singular w.r.t.~SLE$_6$ \cite{NW}. (See also \cite{Au} and \cite[Theorem 13.4]{DPSL} for the much simpler result that the full scaling limits are singular.)

\begin{figure}[htbp]
\SetLabels
\endSetLabels
\centerline{
\AffixLabels{
\includegraphics[width=0.3\textwidth]{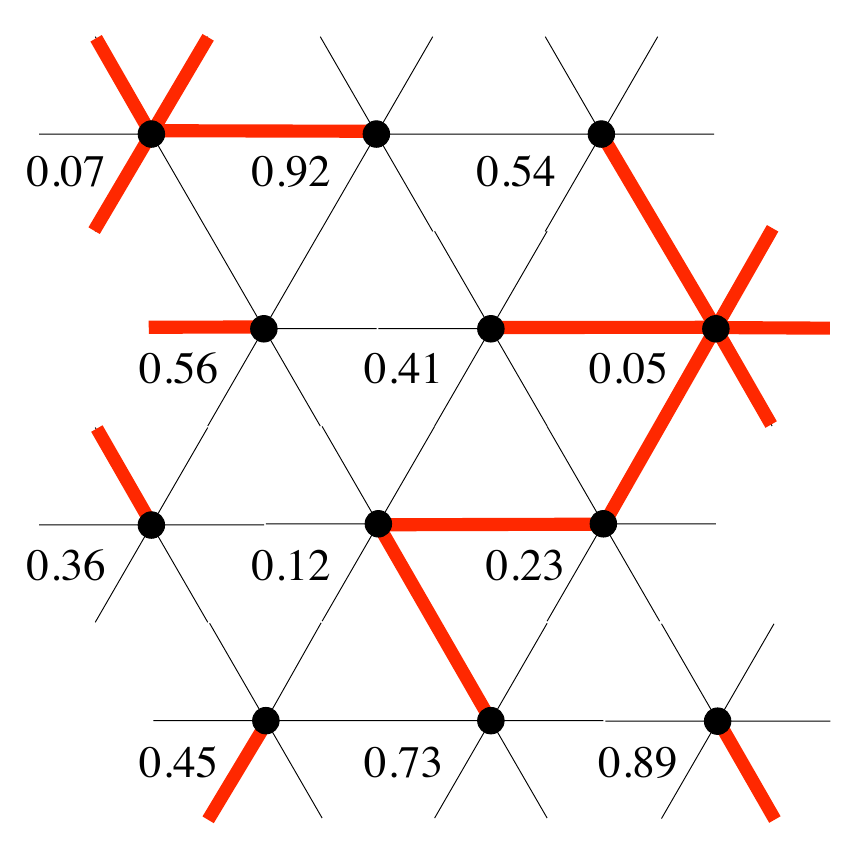}
}}
\caption{The minimal spanning tree associated to vertex labels of the triangular lattice $\Tg$, with a periodic boundary condition.} 
\label{f.vertexTree}
\end{figure}

Since we have a proof of the existence and properties of the scaling limit of the near-critical ensemble only for site percolation on the triangular lattice $\Tg$, if we want to use that to build the $\MST$ scaling limit, we will need a version of the $\MST$ that uses Unif$[0,1]$ vertex labels $\{V(x)\}$ on $\Tg$. So, assign to each edge $e=(x,y)$ the vector label 
\be\label{e.label}
U(e):=\big(V(x)\vee V(y), V(x)\wedge V(y)\big)\,,
\ee 
and consider the lexicographic ordering on these vectors to determine the $\MST$. See Figure~\ref{f.vertexTree}. With a slight abuse of terminology, this is what we will call the $\MST$ on the lattice $\Tg$. Our strongest results will apply to this model, but some of them will also hold for subsequential limits of the usual $\MST$ on $\Z^2$, known to exist by \cite{ABNW}.
 
Let us make an important remark here. The use of the lexicographic ordering for the vector labels~(\ref{e.label}) is somewhat arbitrary, and starting from the same vertex labels, using a different way to get edge labels or using a different natural ordering, one could a priori get an $\MST$ with a very different global structure. In fact, this does happen if the vertex labels are assigned maliciously. Nevertheless, with the Unif$[0,1]$ labels, for any rule to construct the $\MST$ on $\Tg$ that ensures that {\it any two $p$-clusters are connected by a unique path of this $\MST$} (which is exactly how our definition works), our approximation of the macroscopic structure of the $\MST$ using the near-critical ensemble will work with large probability, and hence the scaling limit will be the same.

We can now state our main theorem:

\bth[Limit of $\MST_\eta$ in $\C$]\label{t.main}
As $\eta\to 0$, the spanning tree $\MST_\eta$ on $\eta\Tg$ converges in distribution, in the metric $d_\Omega$ of Definition~\ref{d.ESF} below, to a unique scaling limit $\MST_\infty$ that is invariant under translations, scalings, and rotations.
\eth

The strategy of the proof  will be described in Subsection~\ref{ss.strategy}. As a key step, we also prove convergence in any fixed torus $\Tor_M$; see Theorem~\ref{t.torus}. We work in tori to avoid the technicalities related to boundary issues, but with not too much additional work the extension to finite domains with free or wired boundary conditions would be certainly doable.


In Section~\ref{s.MSTgeom}, strengthening the results of \cite{ABNW}, we study the geometry of the limiting tree $\MST_\infty$. The degree of a vertex in a tree graph has the usual meaning, but the degree of a point in a spanning forest of the plane needs to be defined carefully, which we will do in Subsection~\ref{ss.degtype}. To give an example, a pinching point on an $\MST_\infty$ path should not be called a branching point, but it still gives rise to a degree 4 point. Consequently, stating the results on the geometry of the limiting tree also needs some care, to be done precisely only in Theorem~\ref{t.types}. Nevertheless, here are some of the earlier results and our new ones in rough terms. It was proved in \cite{ABNW} that there is an unspecified absolute bound $k_0$ such that almost surely all degrees in any subsequential limit of $\MST_\eta$ are at most  $k_0$. Furthermore, the set of branching points was shown to be almost surely countable. Here, we will prove that there are almost surely no pinching points, all degrees are bounded by 4, and the set of points with degree 4 is at most countable. We will also prove, in Subsection~\ref{ss.trunk}, that the Hausdorff dimension of the trunk is strictly below $7/4$. However, we do not have a guess for the exact dimension; the situation is similar to the somewhat related problem of finding the percolation chemical distance exponent \cite{Damron}.
\medskip

To conclude this subsection, let us note that the recent works \cite{ABG,ABGM} follow a strategy similar to ours, but in a very different setting: namely, in the {\bf mean-field} case. It is well-known that there is a phase transition at $p=1/n$ for the {\bf Erd{\H o}s-R\'enyi random graphs} $G(n,p)$. Similarly to the above case of planar percolation, it is a natural problem to study the geometry of these random graphs near the transition $p_c=1/n$. It turns out in this case that the non-trivial rescaling is to work with $p=1/n + \lambda/ n^{4/3},\, \lambda\in \R$. If $R_n(\lambda)=(C_n^1(\lambda), C_n^2(\lambda),\ldots)$ denotes the sequence of clusters at $p=1/n + \lambda/ n^{4/3}$, ordered in decreasing order of size, say, then it is proved in \cite{ABG} that as $n\to \infty$, the normalized sequence $n^{-1/3}\, R_n(\lambda)$ converges in law to a limiting object $R_\infty(\lambda)$ for a certain topology on sequences of compact spaces which relies on the Gromov-Hausdorff distance.  This near-critical ensemble $\{ R_\infty(\lambda)\}_{\lambda\in \R}$ has then been used in \cite{ABGM} to obtain a scaling limit as $n\to \infty$ (in the Gromov-Hausdorff sense) of the $\MST$ on the complete graph with $n$ vertices. One could say that \cite{DPSL} is the Euclidean ($d=2$) analogue of the mean-field case \cite{ABG}, and our present paper is the analogue of \cite{ABGM}. However, an important difference is that in the mean-field case one is interested in the intrinsic metric properties (and hence works with the Gromov-Hausdorff distance between metric spaces), while in the Euclidean case one is first of all interested in how the graph is embedded in the plane.

\subsection{The Invasion Percolation Tree $\IP$}\label{ss.introIP}

The connection between $\mathsf{WMSF}$ and critical percolation on infinite graphs can also be seen through {invasion percolation}. For a vertex $x$ in an infinite graph $G(V,E)$, and the labels $\{U(e)\}$, let $T_0=\{x\}$, then, inductively, given $T_n$, let $T_{n+1}=T_n \cup \{e_{n+1}\}$, where $e_{n+1}$ is the edge between $T_n$ and $V\setminus T_n$ that has the smallest label $U$. The {\bf Invasion Percolation Tree} of $x$ is then $\IP(x):=\bigcup_{n\geq 0} T_n$. It is easy to see that, even deterministically, if $U: E(G)\longrightarrow \R$ is an injective labelling of a locally finite graph, then $\mathsf{WMSF}=\bigcup_{x\in V(G)} \IP(x)$.

Once the invasion tree enters an infinite $p$-cluster $\calC$, it will not use edges outside it. Furthermore, it is not surprising (though non-trivial to prove, see \cite{relentless}) that for any transitive graph $G$ and any $p>p_c(G)$, the invasion tree eventually enters an infinite $p$-cluster. Therefore, $\limsup \{ U(e) : e\in \IP(x)\}=p_c(G)$ for any $x\in V(G)$. This way, invasion percolation can be considered as a ``self-organized criticality'' version of critical percolation; finer results for the planar case are given in \cite{ChChN1, Damron1, Damron2}. Moreover, $\IP$ can be used to study Bernoulli percolation itself: e.g., for the well-behavedness of the supercritical phase on $\Z^d$, $d>2$ \cite{ChChN2}, and for uniqueness monotonicity on non-amenable graphs \cite{relentless}. Invasion percolation can be analyzed very well on regular trees \cite{invtree}, with a scaling limit that can be described using diffusion processes \cite{invtreeSL}.

For planar lattices, since $\IP_\eta$ is so intimately related to $\MST_\eta$, it will be quite easy to modify the proof of Theorem~\ref{t.main} for the case of $\IP$; see Section~\ref{s.IP}.

\subsection{The scaling limit of the near-critical ensemble}\label{ss.NCE}

We need to recall how the scaling limit of the near-critical ensemble is constructed in \cite{GPS2a,DPSL}, because the present paper is heavily built on this. To start with, we slightly change the near-critical parametrization given in~(\ref{e.nc}):

\begin{definition}\label{d.NCE}
The {\bf near-critical ensemble} $(\omega_\eta^\lambda)_{\lambda\in \R}$ will denote the following process: 
\bi
\item[(i)] Sample $\omega_\eta^{\lambda=0}$ according to $\P_\eta$, the  law of critical percolation on $\eta \Tg$. We will sometimes represent this as a black-and-white coloring of the faces of the dual hexagonal lattice, with white hexagons standing for closed (empty) sites. 
\item[(ii)] As $\lambda$ increases, closed sites (white hexagons) switch to open (black) at an exponential rate $r(\eta)$, as given after~(\ref{e.nc}).
\item[(iii)] As $\lambda$ decreases, black hexagons switch to white at rate $r(\eta)$.
\ei
Note that, for any $\lambda\in \R$, the near-critical percolation $\omega_\eta^\lambda$ corresponds exactly to a percolation configuration on $\eta \Tg$ with parameter
\[
\begin{cases}
p=p_c + (1-p_c)\, (1-e^{-\lambda\, r(\eta)}) \quad \text{if } \lambda\geq 0 \\
p=p_c \, e^{-|\lambda|\, r(\eta)} \quad \text{if } \lambda< 0\,.
\end{cases}
\]
For any site $x$, the value $\lambda(x)\in\R$ where $x$ switches from closed to open will be called the near-critical percolation label of $x$.

The same definitions can be made on $\eta\Z^2$.
\end{definition}

It is easy to understand intuitively why $r(\eta)$ is the right time rescaling to obtain the near-critical window. Assume that in the unit square there is no left-right crossing in $\omega^{\lambda=0}_\eta$. Then the expected number of those sites that are closed at $\lambda=0$ but are {\bf  pivotal} for the left-right crossing (i.e., opening any of them would establish the crossing) and which actually become open in $\omega^\lambda_\eta$ is known to be of order $\lambda$. Therefore, for $\lambda>0$ small, it is unlikely that a left-right crossing has been established if it was not already there, hence the system must have stayed very close to critical; on the other hand, one may expect that for $\lambda\gg 1$ a crossing is already quite likely, hence the system should already be quite supercritical. This was rigorously proved in \cite{Kesten}. Then, if one wants to describe the {\bf scaling limit} of $\omega^\lambda_\eta$ as $\eta\to 0$, a natural idea that was detailed in \cite{CFNb} is that this should be possible by following which of those points get opened (for $\lambda>0$) or get closed (for $\lambda<0$) that were pivotal at $\lambda=0$ for at least some small macroscopic distance $\eps>0$. To this end, one should look at the counting measure on $\eps$-pivotal points at criticality, normalized such that the measure stays non-trivial as $\eta\to 0$, and hope that these $\eps$-pivotal measures have limits that are measurable w.r.t.~the scaling limit of critical percolation itself. This is the main result of \cite{GPS2a} (with a slight change of what $\eps$-pivotal means). Then, the scaling limit of the near-critical ensemble may be described by taking Poisson point processes of switch times, with intensity measures being these $\eps$-pivotal measures, and by updating the crossings of all the quads (certain generalized rectangles) according to these pivotal switches. This is done in \cite{DPSL}. Here there are roughly two main issues: firstly, it is not immediately clear how one can update the crossings of {\it all} the quads by pivotal switches that are happening at {\it all} spatial and time scales. For this, one should code the percolation configuration in a suitable manner that is minimal enough so that the updates can be done, but rich enough so that it contains all the relevant information. This coding and updating takes up a large part of \cite{DPSL}, done through the so-called $\eps$-networks that we will actually recall in Section~\ref{s.enhanced}. The second main issue is that one needs to prove that despite all the switches that take place as $\lambda$ increases, following the switches of all the initially $\eps$-pivotal sites gives a good idea about the $\eps$-pivotal switches at later times. For this, the key discrete result from \cite{Damron1,DPSL} is the following proposition, which we will often use also in the present paper:

\bpr[Near-critical stability]\label{p.stab}
For any fixed $-\infty<\lambda<\lambda'<\infty$, in the near-critical ensemble on $\eta\Tg$, let $\A^{\lambda,\lambda'}_k (r,R)$ denote the following {\bf near-critical polychromatic $k$-arm event}: there exist $k\geq 2$ disjoint paths in the lattice that connect the boundary pieces of the annulus $B_R(0)\setminus B_r(0)$, each called either ``primal'' or ``dual'', and all the percolation ensemble labels along all the primal arms are at most $\lambda'$, while all the labels along the dual arms are at least $\lambda$. Note that $\lambda=\lambda'$ gives back the usual notion of primal and dual arms in the percolation configuration $\omega^\lambda_\eta$. Then, 
$$
\Pb{\A^{\lambda,\lambda'}_k (r,R)} \leq C_{\lambda,\lambda'} \, \alpha_k(r,R)\,,
$$
where $\alpha_k(r,R)=\alpha^\eta_k(r,R)$ is the polychromatic $k$-arm probability in critical percolation on the same lattice. Similarly, for the {\bf monochromatic $k$-arm events}, $k\ge 1$, where all arms are primal,
$$
\Pb{{\A}^{\lambda'}_k (r,R)} \leq C'_{\lambda,\lambda'} \, \alpha'_k(r,R)\,,
$$
where $\alpha'_k(r,R)={\alpha'}^\eta_k(r,R)$ is the monochromatic $k$-arm probability at criticality. For $\alpha'_1(\cdot,\cdot)$, we will just use $\alpha_1(\cdot,\cdot)$.

The same statements hold for bond percolation on $\eta\Z^2$, just with dual arms being paths in the dual lattice, in the usual manner.
\epr

\begin{remark}\label{r.arms}
For fixed radii $0<r<R<\infty$, the discrete multi-arm probabilities $\alpha^\eta_k(r,R)$ and ${\alpha'}^\eta_k(r,R)$ converge, as $\eta\to 0$, to their SLE$_6$ counterparts (see \cite{WWperc}). In the present paper, we will be interested in these quantities only up to constant factors, not in the details of their convergence, hence their dependence on $\eta$ is not important and will be omitted from the notation. Formulas like $\alpha_k(\eta,1)$ will also be understood on the discrete lattice, always with mesh $\eta$. We will also use the quasi-multiplicativity of multi-arm probabilities (both for the discrete and continuum versions): for any $k\geq 1$, there exists $c_k>0$ such that
$$
c_k\, \alpha_k(r_1, r_2) \, \alpha_k(r_2, r_3) \leq \alpha_k(r_1, r_3) \leq \alpha_k(r_1, r_2) \,\alpha_k(r_2, r_3)\,,
$$
for all $0<r_1<r_2<r_3<\infty$. Similarly for $\alpha'_k$. See, e.g., \cite[Subsection 4.5]{Nolin}.
\end{remark}

The proof of Proposition~\ref{p.stab} for the alternating 4-arm event is given in \cite[Lemma 6.3]{Damron1}, or follows directly from \cite[Lemma 8.4]{DPSL}, which is more general in that it does not assume that the dynamics is monotone in $\lambda$. 
For general $k$, the case of $\lambda=\lambda'$ is known as Kesten's near-critical stability \cite{Kesten}. And just as in Kesten's approach, the proof for general $k$ and general $\lambda<\lambda'$ is a simple modification of the proof for the alternating 4-arm event: the key point is that the pivotality of a site for a general $k$-event still depends on an alternating 4-arm event around that site, and hence the near-critical stability of the alternating 4-arm probability, proved using a recursion in \cite{DPSL}, easily implies the stability of the general $k$-arm event, as well. We omit the details.

The above sketch of the contents of \cite{GPS2a,DPSL} should make it clear that the scaling limit of the near-critical ensemble is constructed entirely from the critical scaling limit, plus independent randomness of the pivotal switch times. Moreover, all the proofs in \cite{GPS2a,DPSL} are universal in the sense that they use lattice-independent discrete percolation technology that have been available since \cite{Kesten}. Altogether, once one proves Cardy's formula for critical percolation on $\eta\Z^2$, which would imply the same scaling limit as on $\eta\Tg$, we would also immediately get that the scaling limit for the entire near-critical ensemble is the same. This universal aspect remains true for the present paper.

\subsection{Strategy of the proof and organization of the paper}\label{ss.strategy}

First of all, in Subsection~\ref{ss.ESF}, we describe the topological space in which the convergence of our random trees will take place:  the space of essential spanning forests in $\C$, introduced in \cite{ABNW}. There are possible alternatives to using this topology, such as the quad-crossing topology of \cite{SSblacknoise} (suggested to us for this purpose by Nicolas Broutin) or the topology introduced in \cite{Schramm} for the scaling limit of the Uniform Spanning Tree. Especially the quad-crossing topology (recalled in Subsection~\ref{ss.quads}) would seem natural, since the scaling limit of near-critical percolation is taken in this space. Nevertheless, we chose the topology of \cite{ABNW} for several reasons: that was the first paper dealing with subsequential scaling limits of $\MST_\eta$, proving results that we are sharpening here; using this topology to describe paths in the spanning trees is not harder than using quad-crossings, while it also gives a natural way to glue the paths into more complicated trees; there is a simple explicit metric generating this topology. However, we will unfortunately need more topological preparations than just recalling these definitions, because the minimalist structure, based on just the pivotal measures of \cite{GPS2a}, which was enough to describe the scaling limit of the near-critical ensemble in \cite{DPSL}, will not be enough for the tree structures of the present paper. In particular, in Proposition~\ref{p.PivSetCol}, we will prove that that {\bf set of colored pivotals} also has a limit as $\eta\to 0$. 

In Section~\ref{s.enhanced}, we first recall the definition of the {\bf networks} $\Net^{\vl,\eps}_\eta$ and $\Net^{\vl,\eps}_\infty$ introduced in \cite{DPSL}, where $\vl=(\lambda,\lambda')$ is a pair of near-critical parameters with $\lambda<\lambda'$. These are graphs with vertex sets $X$ given by those $\eps$-pivotals in the configuration $\omega^\lambda$ on a torus $\Tor_M$ that experience a switch between level $\lambda$ and $\lambda'$, and edges given roughly by the primal and dual connections in $\omega^\lambda\setminus X$. Then we need to add a bit more structure to these networks, creating the so-called {\bf enhanced networks}: roughly, we will need to know which of these pivotals are connected together by an open cluster of $\omega^\lambda\setminus X$, and will need to know the colors of these pivotals in $\omega^\lambda$. For this, we will use Proposition~\ref{p.PivSetCol} mentioned in the previous paragraph and Proposition~\ref{p.clusterdim} saying that clusters of large diameter also have large volume (which excludes certain pathological geometric behaviour that would ruin the construction). From these {enhanced networks}, we will obtain finite labelled graphs whose vertices will basically be open $\lambda$-clusters that have $\eps$-pivotals switching in the time interval $(\lambda,\lambda')$, with edges labelled by the times of the pivotal switches, showing how the $\lambda$-clusters merge. We will define the $\MST$ on this finite labelled graph, denoted by $\MST^{\vl,\eps}_\eta$ in the discrete and $\MST^{\vl,\eps}_\infty$ in the continuum case --- these are basically the macroscopic approximations to the cluster trees that we discussed in Subsection~\ref{ss.introMST}. To be more precise, in Section~\ref{s.enhanced} we define only some Minimal Spanning Forests, and we need a bit more work until in Lemma~\ref{l.giant} we can actually define the trees. The fact that these approximating {\bf cut-off trees} $\MST^{\vl,\eps}_\eta$ and $\MST^{\vl,\eps}_\infty$ are close to each other if the underlying near-critical ensembles $\omega^{[\lambda,\lambda']}_\eta$ and $\omega^{[\lambda,\lambda']}_\infty$ are close follows easily from \cite{DPSL}. 

In Section~\ref{s.epsapprox} we prove that the cut-off trees $\MST^{\vl,\eps}_\eta$ are close to the true $\MST_\eta$ if $\lambda\ll -1$, $\lambda'\gg 1$, and $\eps>0$ is small. Here the key technique is near-critical stability, Proposition~\ref{p.stab}. 

Summarizing, we get that $\MST_\eta$ is close to $\MST^{\vl,\eps}_\infty$. Since the latter does not depend on $\eta$, while the former does not depend on $\vl$ and $\eps$, they both need to be close to an object that does not depend on any of these parameters: this will be the scaling limit $\MST_\infty$. To give a succinct pictorial summary of this strategy:
\begin{center}
\begin{tikzpicture}
  \node (A) {$\MST_\eta$};
  \node (B) [node distance=2.5cm, right of=A] {$\MST^{\vl,\eps}_\eta$};
  \node (C) [node distance=2.5cm, right of=B] {$\MST^{\vl,\eps}_\infty$};
  \node (D) [node distance=1.5cm, below of=B] {$\MST_\infty$};
  \draw[<->] (A) to (B);
  \draw[<->] (B) to (C);
  \draw[->] (A) to (D);
  \draw[->] (C) to (D);
\end{tikzpicture}
\end{center}
This conclusion will be materialized in Section~\ref{s.main}, together with the extension from the case of the tori $\Tor_M$ to the full plane, and with the proof of the claimed invariance properties.

As already advertised in Subsections~\ref{ss.introMST} and~\ref{ss.introIP}, the results on the geometry of $\MST_\infty$ are discussed in Section~\ref{s.MSTgeom}, while Section~\ref{s.IP} establishes the existence and invariance properties of $\IP_\infty$. We conclude the paper with some open problems in Section~\ref{s.conj}.

\bigskip
\noindent{\bf Acknowledgments.} We thank  Louigi Addario-Berry, Nicolas Broutin, Laure Dumaz, Gr\'egory Miermont and David Wilson for stimulating discussions, Rob van den Berg for pointing out the connection between Proposition~\ref{p.clusterdim} and \cite{Jarai}, and Alan Hammond and two amazing anonymous referees for very good comments on the manuscript.  

Part of this work was done while all authors were at Microsoft Research, Redmond, WA, or GP was visiting CG at ENS Lyon. CG was partially supported by the ANR grant MAC2 10-BLAN-0123. GP was supported by an NSERC Discovery Grant at the University of Toronto, an EU Marie Curie International Incoming Fellowship at the Technical University of Budapest, and partially supported by the Hungarian National Research, Development and Innovation Office, NKFIH grant K109684, and by the MTA R\'enyi Institute ``Lend\"ulet'' Limits of  Structures Research Group.

 
\section{Topological and measurability preliminaries}\label{s.topology}

\subsection{The space of essential spanning forests }\label{ss.ESF}
 
The following topological setup for discrete and continuum spanning trees was introduced in \cite{ABNW}. We are summarizing here the definitions and the notation, with small modifications; the main difference is roughly that $\Omega$ will also contain spanning trees of subsets of the complex plane, to accommodate the invasion percolation tree $\IP$ and our approximating trees $\MST^{\vl,\eps}$.
 
We will work in a one-point compactification of $\C=\R^2$, denoted by $\hat\C=\C\cup\{\infty\}$, with the Riemannian metric 
\be\label{e.metric}
\frac{4}{(1+x^2+y^2)^2} \big(dx^2+dy^2\big)\,;
\ee
by stereographic projection, $\hat\C$ is isometric with the unit sphere. Note that this metric is equivalent to the Euclidean metric in bounded domains, while the distance between any two points outside the square of radius $M$ around the origin in $\C$ is at most $O(1/M)$. This will imply that convergence of spanning trees in $\hat\C$ is the same as convergence within bounded subsets of $\C$. This is necessary, since convergence of random spanning trees cannot be uniform in $\C$: on $\eta\Z^2$, inside the infinitely many pieces $[i,i+1)\times[j,j+1)$, $i,j\in\Z$, one can find arbitrary topological behavior (e.g., macroscopically vanishing areas with arbitrarily large numbers of macroscopic branches emanating from them) that will be very far from the almost sure behavior of the continuum tree.

Spanning trees on infinite graphs are usually defined and studied as weak limits of spanning trees in finite subgraphs exhausting the infinite graph. For these finite graphs, one may consider different boundary conditions: most importantly, free or wired. As mentioned in the Introduction, for the $\MST$ on Euclidean planar lattices, all such boundary conditions give the same limit measure, and we will work in the tori $\Tor_M$ of side-length $2M$, which can be realized as the subdomains $[-M,M)^2$ of $\C$, or even as subgraphs of $\eta\Tg$ for suitable values of $M$, with a periodic boundary condition (which is sandwiched between the free and the wired conditions). See Figure~\ref{f.vertexTree} in the Introduction.

\begin{definition}\label{d.immertrees}
A {\bf  reference tree} $\tau$ is a tree with a finite set  of leaves  (or external vertices), denoted by $\xi(\tau)$, with each edge considered to be a unit interval. A {\bf reparametrization} is a continuous map $\phi:\tau\lora \tau$ that fixes all the vertices and is monotone on the edges. An {\bf immersed tree indexed by $\tau$} is an equivalence class of continuous maps $f:\tau\lora\hat\C$, where $f_1$ and $f_2$ are considered equivalent if there exist reparametrizations $\phi_1,\phi_2$ with $f_1\circ\phi_1 = f_2\circ\phi_2$. 
The collection of immersed trees indexed by $\tau$ is denoted by $\Set_\tau$, and we set 
\be\label{e.SetUnion}
\Set^{(\ell)}:=\bigcup_{\tau:\, |\xi(\tau)|=\ell} \Set_\tau\,.
\ee
Immersed trees with leaves $x_1,\dots,x_\ell\in\hat\C$ will often be denoted by $T(x_1,\dots,x_\ell)\in \Set^{(\ell)}$. 

We will also consider trees immersed into the torus $\Tor_M$ with the flat Euclidean metric; the corresponding collection of immersed trees with $\ell$ leaves is denoted by $\Set^{(\ell)}_M$.

One may consider trees immersed not just into $\hat\C$ or $\Tor_M$, but into a graph $G(V,E)$ that is embedded into $\hat\C$ or $\Tor_M$, and then the image of $\tau$ is required to be a subtree of $G(V,E)$, with its vertices mapped into $V$ and any of its edges mapped to a union of edges from $E$.
\end{definition}

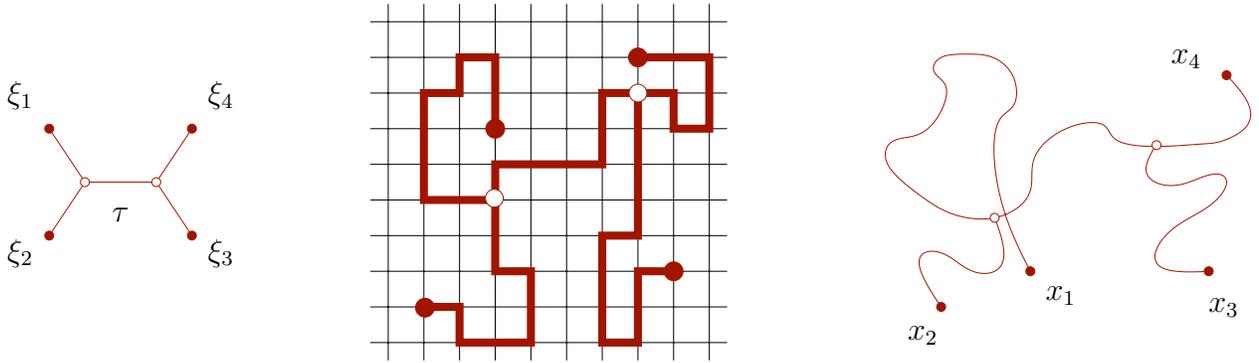
\begin{figure}[htbp]
\SetLabels
(0.08*0.4)$\tau$\\
(0*0.7)$\xi_1$\\
(0*0.3)$\xi_2$\\
(0.16*0.3)$\xi_3$\\
(0.16*0.7)$\xi_4$\\
(0.83*0.2)$x_1$\\
(0.72*0.1)$x_2$\\
(0.96*0.17)$x_3$\\
(0.93*0.8)$x_4$\\
\endSetLabels
\centerline{
\AffixLabels{
\includegraphics[width=\textwidth]{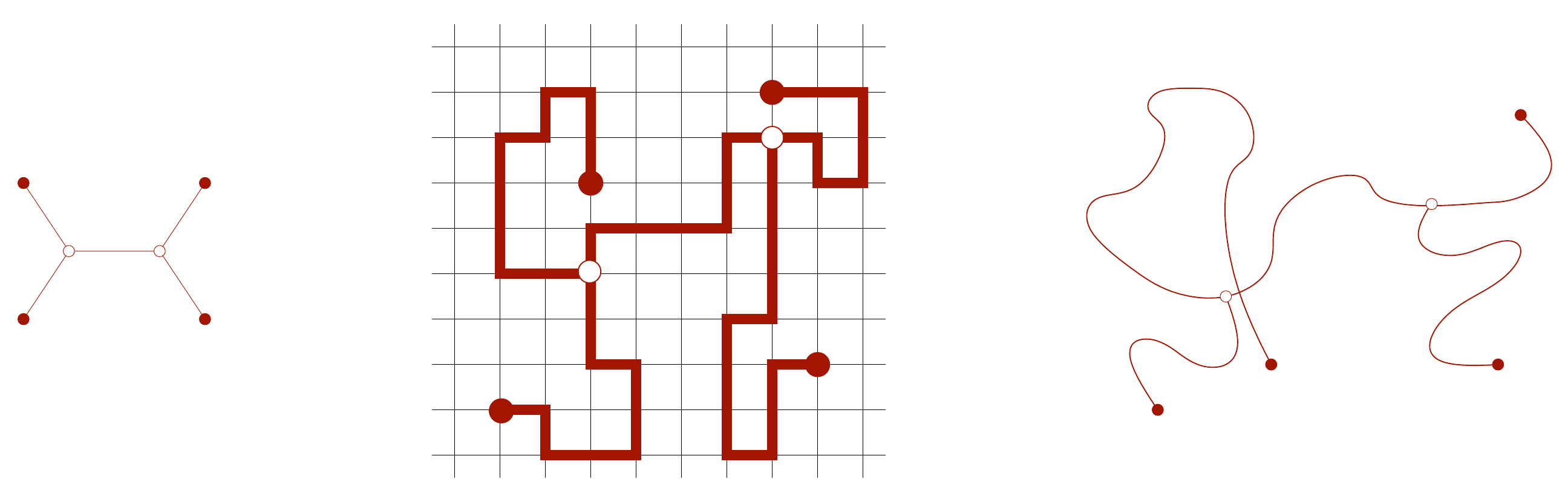}
}}
\caption{A reference tree $\tau$ with four leaves, with one immersion into $\Z^2$ and another into $\C$. The {\it image} in $\C$ is not a tree, but this is allowed. In the scaling limit of any discrete random tree in $\hat\C$ one cannot see such self-intersections, but could see touch-points, and self-intersections might happen in scaling limits in higher dimensions.}
\label{f.immersed}
\end{figure}

Note that if a reference tree $\tau'$ is given by contracting some edges of some $\tau$, denoted by $\tau'\prec\tau$, then $\Set_{\tau'}$ is naturally a subset of $\Set_\tau$, represented by maps $f:\tau\lora\hat\C$ that are constants on the contracted edges. By contractions in two non-isomorphic trees, $\tau_1$ and $\tau_2$, we may reach the same tree $\tau'\prec\tau_i$, hence $\Set^{(\ell)}$ may be viewed as covered by patches $\Set_\tau$ that are sewn together along ``smaller dimensional'' patches $\Set_{\tau'}$, similarly to a simplicial complex. (In particular, after these identifications, (\ref{e.SetUnion}) stops being a disjoint union.)

We now equip each $\Set_\tau$ with a very natural metric, extending the notion of uniform closeness up to reparametrization of curves: for two immersed trees $f_1,f_2:\tau\lora\hat\C$,
\be\label{e.TreeMetric}
\dist_\tau(f_1,f_2)=\inf_{\phi_1,\phi_2} \sup_{t\in\tau} \dist_{\hat\C}\big(f_1\circ\phi_1(t), f_2\circ\phi_2(t)\big)\,,
\ee
where the $\phi_i$'s run over all reparametrizations of $\tau$. This can be easily extended to immersed trees indexed by different reference trees: by the above remark about patches, for any pair of reference trees $\tau,\tau'$ there exist sequences $\tau=\tau_0,\tau_1,\dots,\tau_m=\tau'$ such that $\tau_i\prec \tau_{i+1}$ or $\tau_i\succ \tau_{i+1}$  for all $i=0,\dots,m-1$, and then for any $f:\tau\lora\hat\C$ and $f':\tau'\lora\hat\C$ we can take 
$$\dist(f,f')=\inf\left\{ \sum_{i=0}^{m-1} \dist_{\tau_i\curlyvee\tau_{i+1}}(f_i,f_{i+1}) : f_0=f,\  f_m=f',\   \tau_i\mathop{\lora}\limits^{f_i}\hat\C\text { for }i=1,\dots,m-1\right\},$$ where, with a rather obvious notation, $\tau_i\curlyvee\tau_j=\tau_i$ if $\tau_i\succ\tau_j$. For instance, for any $\tau,\tau'$ there exists $\tau''$ with $\tau,\tau'\prec\tau''$, hence $\dist(f,f')\leq \dist_{\tau''}(f,f')$.

With this metric, $\Set^{(\ell)}$ is clearly a complete separable metric space, called the {\bf space of $\ell$-trees}. Of course, a Cauchy sequence of trees contained fully in $\C$ might have a limit that has an edge going through $\infty$. Similarly, $\Set^{(\ell)}_M$ is complete and separable with the analogous metric, just using the Euclidean metric on $\Tor_M$ in~(\ref{e.TreeMetric}).

Now that we have a definition for the space of finite trees immersed in $\hat\C$ or $\Tor_M$, we can start defining what a spanning tree of $\hat\C$ or $\Tor_M$ should be: a set of finite trees that satisfy certain compatibility conditions.

The set of non-empty closed subsets of $\Set^{(\ell)}$ in the above metric, equipped with the Hausdorff metric, is denoted by $\Omega^{(\ell)}$. We will consider graded sets 
$$
\For=\big( \For^{(\ell)} \big)_{\ell\geq 1} ~\in~ \Omega^\times:=\bigtimes_{\ell\geq 1} \Omega^{(\ell)}\,,
$$ 
with the product topology. Clearly, $\Omega^\times$ is again complete, separable and metrizable; in one word, it is a Polish metric space.

Extending the map $\tau\mapsto \xi(\tau)$ giving the external vertices of an index tree,  for any $\For\in\Omega^\times$ we can define
$$
\xi(\For):=\bigcup \Big\{ f(\xi(\tau)) : \tau\mathop{\lora}\limits^{f} \hat\C \in\For^{(\ell)},\ \ell\geq 1 \Big\} \subset \hat\C\,,
$$ 
which gives the set of external vertices occurring in $\For$. It is clearly a Borel measurable function, since for any open $U\subset\hat\C$, the preimage $\xi^{-1}(U)$ is a countable intersection (over $\ell\ge 1$) of open sets.

Let $\Set_{B_1,\dots,B_\ell}$ be the set of immersed trees with endpoints $x_i \in B_i$, where each $B_i$ is a closed subset of $\hat\C$. Note that this is a closed subset of $\Set^{(\ell)}$. The non-empty closed subsets of $\Set^{(\ell)}$ that do not intersect $\Set_{B_1,\dots,B_\ell}$ form an open set in the Hausdorff metric, hence the map 
$$\Omega^\times\lora\Omega_{B_1,\dots,B_\ell}\subseteq \Omega^{(\ell)}\,,\qquad \For \mapsto\For^{(\ell)}\cap \Set_{B_1,\dots,B_\ell}$$ 
is measurable. In words, extracting  the subtrees of $\For$ with leaves in prescribed closed sets (e.g., the branches of $\For$ connecting two given points) is a measurable map.

\bde\label{d.ESF}
A graded set $\For=\big( \For^{(\ell)} \big)_{\ell\geq 1} \in \Omega^\times$ is called an {\bf essential spanning forest} on its external vertices $\xi(\For)$ if it satisfies the following properties:
\begin{itemize}
\item[{\bf (i)}] for each $\ell\in\N^+$ and any $\ell$-tuple $\{x_1,\dots,x_\ell\}$ of vertices in $\xi(\For)$, there exists at least one immersed tree $T(x_1,\dots,x_\ell)\in\For^{(\ell)}$ with those leaves; 
\item[{\bf (ii)}] for any immersed tree $T\in\For^{(\ell)}$, any subtree $T'\subset T$ (given by restricting the immersion to a combinatorial subtree of the index tree $\tau$) is again in some $\For^{(\ell')}$;
\item[{\bf (iii)}] for any two trees $T_i\in\For^{(\ell_i)}$, $i=1,2$, there is a tree in some $\For^{(\ell)}$ that contains both $T_i$'s as subtrees and has no leaves beyond those of the $T_i$'s.
\end{itemize}

Note that {\em (ii)} implies that $\xi(\For)$ contains all the vertices of all the embedded trees.

An essential spanning forest $\For$ is called a {\bf spanning tree} if $\xi(\For)\subset\C$ and every path $T(x,y)\in\For^{(2)}$ stays within a bounded region of $\C$. A spanning tree is called {\bf quasi-local} if for any bounded $\Lambda\subset\C$ there exists a bounded domain $\bar\Lambda(\For,\Lambda)\subset\C$ such that every tree of $\For$ with leaves in $\Lambda$ is contained in $\bar\Lambda$. 

The set of essential spanning forests in $\hat\C$ (with an arbitrary set of vertices $\xi(\For)$) will be denoted by $\Omega$. It is easy to check that $\Omega$ is a closed subset of the Polish space $\Omega^\times$, hence itself is Polish. A simple explicit {\bf metric}, denoted by $d_\Omega$, is given by the restriction from $\Omega^\times$ to $\Omega$ of the sum over $\ell$ of the Hausdorff distance on $\Set^{(\ell)}$ multiplied by the weight $2^{-\ell}$.

For the tori $\Tor_M$, the spaces $\Omega^{(\ell)}_M$, $\Omega^\times_M$, $\Omega_M$ are defined analogously, with the only difference being that any essential spanning forest here is a single tree. The metric $d_{\Omega_M}$ is defined the same way as $d_\Omega$.
\ede

The only way in which two vertices may be disconnected in an essential spanning forest $\For$ in $\hat\C$ is that all the paths between them go through $\infty$; therefore, either $\For$ is a spanning tree, or no component of it is contained in a bounded domain of $\C$. This is the property that the adjective ``essential'' for these spanning forests refers to. (In the setting of discrete infinite graphs, this reduces to saying that all components of the forest are infinite trees.) Also, note that the above definition allows for having more than one path between two vertices, which will in fact happen in the scaling limit of the $\MST$.

\subsection{The quad-crossing topology}\label{ss.quads}

Let us quickly recall the notation and the basic results for the quad-crossing topology of percolation configurations, introduced in \cite{SSblacknoise} and studied further in \cite{GPS2a,DPSL}.

Let $D\subset \hat\C=\C\cup\{\infty\}$ be open, or be equal to the torus $\Tor_M$. A {\bf quad} in the domain $D$ can be considered as a homeomorphism $\Quad$ from $[0,1]^2$ into $D$. The space of all quads in $D$, denoted by $\QUAD_D$, can be equipped with the following metric: $d_\QUAD(\Quad_1,\Quad_2):=\inf_{\phi}\sup_{z\in \p [0,1]^2} |\Quad_1(z)-\Quad_2(\phi(z))|$, where the infimum is over all homeomorphisms $\phi: [0,1]^2 \longrightarrow [0,1]^2$ which preserve the 4 corners of the square. A {\bf crossing} of a quad $\Quad$ is a connected closed subset of $[\Quad]:=\Quad([0,1]^2)$ that intersects both $\p_1\Quad=\Quad(\{0\}\times[0,1])$ and $\p_3\Quad=\Quad(\{1\}\times[0,1])$. We say that $\Quad$ has a {\bf dual crossing} between $\p_1\Quad$ and $\p_3\Quad$ by some closed subset $S\subseteq[\Quad]$ if there is no crossing in $S$ between $\p_2\Quad=\Quad([0,1]\times\{0\})$ and $\p_4\Quad=\Quad([0,1]\times\{1\})$.

From the point of view of crossings, there is a natural partial order on $\QUAD_D$: we write $\Quad_1 \leq \Quad_2$ if any crossing of $\Quad_2$ contains a crossing of $\Quad_1$.  
Furthermore, we write $\Quad_1 < \Quad_2$ if there are open neighborhoods $\mathcal{N}_i$ of $\Quad_i$ (in the uniform metric) such that $ N_1\leq N_2$ holds for any $N_i\in \mathcal{N}_i$.  A subset $S\subset \QUAD_D$ is called {\bf hereditary} if whenever $\Quad\in S$ and $\Quad'\in\QUAD_D$ satisfies $\Quad' < \Quad$, we also have $\Quad'\in S$. The collection of all closed hereditary subsets of $\QUAD_D$ will be denoted by $\HH_D$. Any discrete percolation configuration $\omega_\eta$ of mesh $\eta>0$, considered as a union of the topologically closed percolation-wise open hexagons in the plane, naturally defines an element $S(\omega_\eta)$ of $\HH_D$: the set of all quads for which $\omega_\eta$ contains a crossing. In particular, near-critical percolation at level $\lambda\in\R$, as defined in Definition~\ref{d.NCE}, induces a probability measure on $\HH_D$, which will be denoted by $\P_\eta^\lambda$.

By introducing a natural topology, $\HH_D$ can be made into a compact metric space. Indeed, let 
$$
\boxminus_\Quad:=\{S\in\HH_D:\Quad\in S\} \qquad \text{for any }\Quad\in \QUAD_D\,,
$$ 
and let 
$$
\boxup_U:=\{S\in\HH_D: S\cap U=\emptyset\}\qquad \text{for any open }U\subset \QUAD_D\,.
$$ 
Then, define $\T_D$ to be the minimal topology that contains every $\boxminus_\Quad^c$ and $\boxup_U^c$ as open sets. It is proved in \cite[Theorem 3.10]{SSblacknoise} that for any nonempty open $D$, the topological space $(\HH_D,\T_D)$ is compact, Hausdorff, and metrizable. Furthermore, for any dense $\QUAD_0 \subset \QUAD_D$, the events $\{\boxminus_\Quad : \Quad\in\QUAD_0\}$ generate the Borel $\sigma$-field of $\HH_D$. An arbitrary metric generating the topology $\T_D$ will be denoted by $d_\HH$. Now, since Borel probability measures on a compact metric space are always tight, we have subsequential scaling limits of $\P_\eta^\lambda$ on $\HH_D$, as $\eta=\eta_k\to 0$. Moreover, the following convergence of probabilities holds. For critical percolation, $\lambda=0$, it is Corollary 5.2 of \cite{SSblacknoise}; for general $\lambda$, the exact same proof works, using that the RSW estimates hold in near-critical percolation.

\begin{lemma}\label{l.crossconv}
For any $\lambda\in\R$, any subsequential scaling limit $\P^\lambda_{\eta_k} \to \P^\lambda_\infty$, and any quad $Q\in \QUAD_D$, one has $\P_\infty^\lambda[\p \boxminus_Q] =0$. Therefore, by the weak convergence of $\P^\lambda_{\eta_k}$ to $\P^\lambda_\infty$,
\[
\P_{\eta_k}^\lambda[\boxminus_Q] \to \P_\infty^\lambda[\boxminus_Q]\,. 
\] 
\end{lemma}

For the case of site percolation on $\eta\Tg$, we know much more than just the existence of subsequential limits. As explained in \cite[Subsection 2.3]{GPS2a}, the existence of a unique quad-crossing scaling limit for $\lambda=0$ follows from the loop scaling limit result of \cite{Smirnov,CNfull}. The case of general $\lambda$ is Theorem 1.4 of \cite{DPSL}:

\begin{theorem}[Near-critical scaling limit]\label{t.NCSL}
For any $\lambda\in\R$, there is a unique measure $\P^\lambda_\infty$ for percolation configurations $\omega^\lambda_\infty$ in $(\HH_D,\T_D)$ such that the weak convergence $\omega^\lambda_\eta \xrightarrow{d} \omega^\lambda_\infty$ holds.
\end{theorem}

\begin{figure}[htbp]
\SetLabels
(0.66*0.1)\small{$\p_2 A=\gamma$}\\
(0.85*0.75)\red{\small{$\overline U_\eps$}}\\
\endSetLabels
\centerline{
\AffixLabels{
\includegraphics[width=0.8\textwidth]{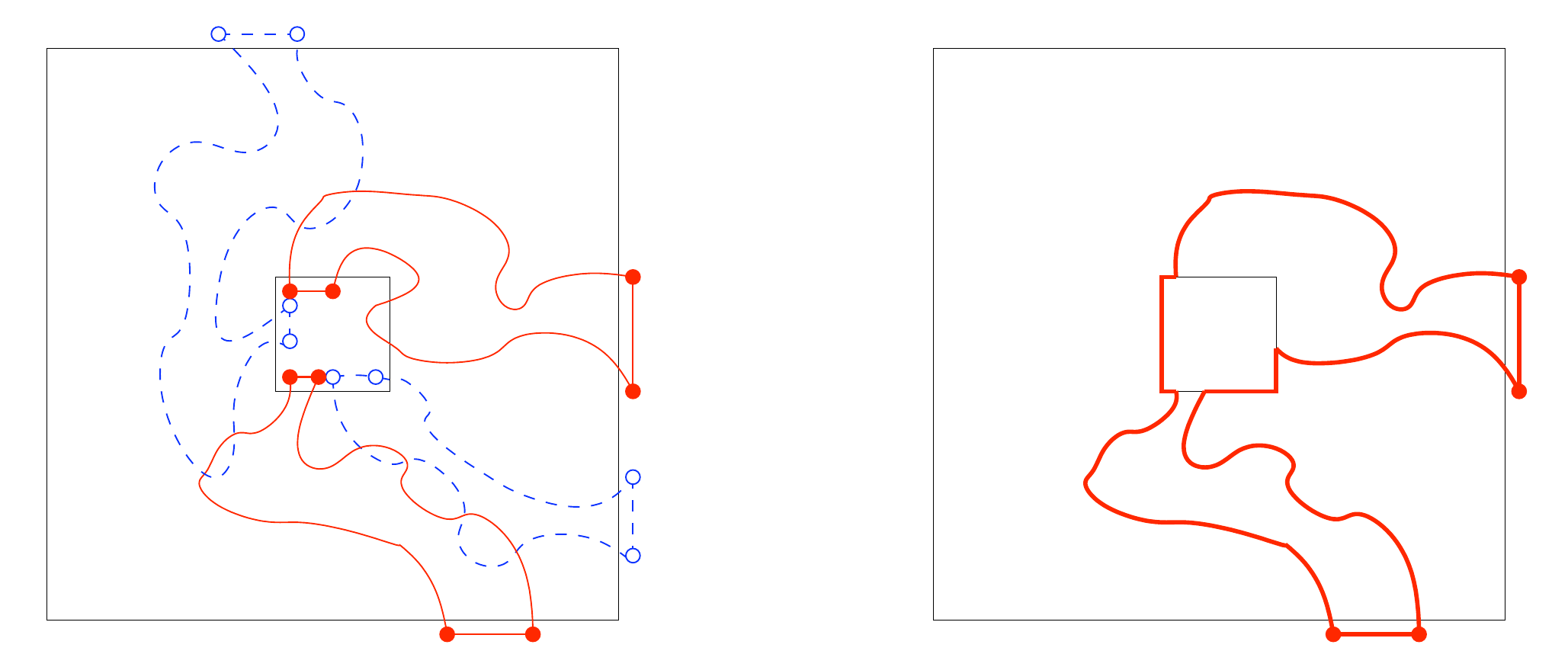}
}}
\caption{Defining the alternating 4-arm event (in Subsection~\ref{ss.quads}) and the color of a pivotal point (in Subsection~\ref{ss.pivotals}) using quad-crossings. Quads with a primal crossing are in solid red, quads with a dual crossing are in dashed blue.} 
\label{f.PivColor}
\end{figure}

We have shown in \cite{GPS2a} that the arm events between the boundary pieces of an annulus are measurable w.r.t.~the quad-crossing topology, and the convergence of probabilities (analogous to Lemma~\ref{l.crossconv}) holds. Namely, for any topological annulus $A\subset D$ with piecewise smooth inner and outer boundary pieces $\p_1 A$ and $\p_2 A$ (and for the case of $D=\Tor_M$, we also require $A$ to be null-homotopic), we define the {\bf alternating 4-arm event} in $A$ as $\A_4=\bigcup_{\delta>0} \A_4^\delta$, where $\A_4^\delta$ is the existence of quads $\Quad_i \subset D$, $i=1,2,3,4$, with the following properties (see the left side of Figure~\ref{f.PivColor}):
\begin{itemize}
\item[(i)] $Q_1$ and $Q_3$ are disjoint and are at distance at least $\delta$ from each other; the same for $Q_2$ and $Q_4$; 
\item[(ii)] for $i\in \{1,3\}$, the sides $\p_1Q_i=Q_i(\{0\}\times [0,1])$  lie inside $\p_1 A$ and the sides $\p_3Q_i=Q_i(\{1\}\times [0,1])$ lie outside  $\p_2 A$; for $i\in \{2,4\}$, the sides $\p_2Q_i=Q_i([0,1] \times \{0\})$ lie inside $\p_1 A$ and the sides $\p_4Q_i=Q_i([0,1]\times \{1\})$ lie outside  $\p_2 A$; all these sides are at distance at least $\delta$ from the annulus $A$ and from the other $Q_j$'s;
\item[(iii)] the four quads are ordered cyclically around $A$ according to their indices;
\item[(iv)]  For $i\in \{1,3\}$, we have $\omega\in \boxminus_{Q_i}$, while for $i\in \{2,4\}$, we have $\omega\in  \boxminus_{Q_i}^c$. In plain words, the quads $Q_1,Q_3$ are crossed, while the quads $Q_2,Q_4$ are dual crossed between the boundary pieces of $A$, with a margin $\delta$ of safety.
\end{itemize}

The definitions of general {\bf (mono-  or polychromatic) $k$-arm events} in $A$ are of course analogous: for arms of the same color we require the corresponding quads to be completely disjoint, and we still require all the boundary pieces lying outside the annulus $A$ to be disjoint.

The following lemma is proved for critical percolation in Lemma 2.9 
of \cite{GPS2a}. For near-critical percolation, the same proofs work, using the stability of multi-arm probabilities (see Lemma 8.4 and Proposition 11.6 of \cite{DPSL}, or \cite{Kesten}), together with the existence of the near-critical scaling limit \cite[Theorem 1.4]{DPSL}.

\begin{lemma}\label{l.arms}
Let $A\subset D$ be a piecewise smooth topological annulus (with finitely many non-smooth boundary points). Then the 1-arm, the alternating 4-arm and any polychromatic 6-arm event in $A$, denoted by $\A_1$, $\A_4$ and $\A_6$, respectively, are measurable w.r.t.~the scaling limit of critical percolation in $D$, and one has 
\[
\lim_{\eta \to 0}\P_\eta^\lambda[\A_i]=\P_\infty^\lambda[\A_i]\,.
\]
Moreover, in any coupling of the measures $\{\P_{\eta}^\lambda\}$ and $\P_\infty^\lambda$ on $(\HH_D,\T_D)$ in which $\omega_{\eta}^\lambda \to\omega^\lambda_\infty$ a.s.~as $\eta\to 0$,  we have 
\begin{align}\label{e.Delta}
\Pb{\{ \omega_\eta^\lambda\in \A_i \} \triangle \{\omega^\lambda_\infty \in \A_i\}} \to 0\qquad (\text{as }\eta\to 0)\,.
\end{align}
\end{lemma}



\subsection{Pivotals and pivotal measures}\label{ss.pivotals}

In \cite{DPSL}, we managed to describe the changes of macroscopic connectivities in a percolation configuration under the stationary or the asymmetric near-critical dynamics using just the pivotal measures of \cite{GPS2a}, without making explicit use of notions like clusters or the set of pivotal sites in continuum percolation. Unfortunately, the situation is slightly more complicated for the models in the present paper, hence we need some foundational work in addition to what was done in \cite[Section 2.4]{GPS2a}.


Let $x$ be a point surrounded (with a positive distance) by a piecewise smooth Jordan curve $\gamma\subset D$, where ``surrounded'' means that $D\setminus\gamma$ has two connected components, with the one containing $x$ being homeomorphic to a disk. For any $\eps>0$, fix a lattice $\eps\Z^2$ in $D$, and let $B_\eps(x)$ be the $\eps$-square $[i,i+1)\eps \times [j,j+1)\eps$ in the lattice that contains $x$. We say that  $x$ is {\bf pivotal for $\gamma$} in $\omega_\infty^\lambda$ if, for any $\eps>0$ such that $B_\eps(x)$ is surrounded by $\gamma$, the alternating 4-arm event occurs in the annulus with boundary pieces $\p B_\eps(x)$ and $\gamma$, as defined in Subsection~\ref{ss.quads}. We let $\Piv^\gamma$ denote the set of pivotal points for $\gamma$ in $D$.
Furthermore, we can identify the {\bf color of a pivotal point} $x\in\Piv^\gamma$ as {\bf open} (black) versus {\bf closed} (white, empty), as follows. We let $\PivO^{\gamma,\eps}$ denote the set of points $x$ for which $\gamma$ surrounds $x$ without intersecting or touching $B_\eps(x)$, and there exist quads $Q_{\eps,i}$, $i=1,2,3,4$, exhibiting the 4-arm event from $\p B_\eps(x)$ to $\gamma$ such that the quad $\overline U_\eps$, given by taking the union of $U_\eps:=Q_{\eps,1}\cup Q_{\eps,3} \cup B_\eps(x)$ and the bounded components of $\C\setminus U_\eps$, is crossed between the boundary pieces $Q_{\eps,1}(\{1\}\times [0,1])$ and $Q_{\eps,3}(\{1\}\times [0,1])$; see the right side of Figure~\ref{f.PivColor} in the previous subsection. Then, we let the set of open pivotals for $\gamma$ be 
$$
\PivO^{\gamma}:=\left\{ x \in D : x\in\PivO^{\gamma,\eps} \textrm{ for all }\eps>0\text{ s.t.~}B_\eps(x)\textrm{ is surrounded by }\gamma\right\}.
$$ 
Clearly, the event $x\in\PivO^{\gamma}$ is measurable w.r.t.~the quad-crossing topology. We will use the notation $x\in\PivO^{\gamma,\eps,\delta}$ for the event that all the crossing events in $\PivO^{\gamma,\eps}$ are satisfied even with a $\delta$ margin of safety. Finally, we set $x\in \PivC^\gamma$ if the analogous dual crossing holds in the quad given by $Q_{\eps,2}\cup Q_{\eps,4} \cup B_\eps(x)$, for each small enough $\eps>0$. 

Note that for a discrete percolation configuration $\omega_\eta^\lambda$ the above definitions do not work: instead of taking all small enough $\eps>0$, we just need to take the annulus between $\gamma$ and the hexagon of the point $x\in D$. And here it is clear what the sets $\PivO^\gamma(\omega_\eta^\lambda)$ and $\PivC^\gamma(\omega_\eta^\lambda)$ are: their disjoint union is the set of pivotal hexagons $\Piv^\gamma(\omega_\eta^\lambda)$, and the color is determined by the color of the hexagon itself. We will also use notation like $x\in\PivO^{\gamma,\eps}(\omega_\eta^\lambda)$: it has the meaning given above, using quad-crossings, and of course it cannot hold unless $\eta$ is small enough, say $\eps>2\eta$, so that $\p B_\eps(x)$ already intersects at least four $\eta$-hexagons.

\bpr[The set of pivotals, with colors]\label{p.PivSetCol}
In any coupling of the measures $\{\P_\eta^\lambda\}$ and $\P_\infty^\lambda$ on $(\HH_D,\T_D)$ in which $\omega_{\eta}^\lambda\xrightarrow{a.s.}\omega_\infty^\lambda$ as $\eta\to 0$,  for any piecewise smooth null-homotopic Jordan curve $\gamma\subset D$ we have the following statements:
\begin{itemize}
\item[{\bf (i)}] $\PivO^\gamma(\omega_\eta^\lambda)$ converges in probability  to $\PivO^\gamma(\omega_\infty^\lambda)$ in the Hausdorff metric of closed sets. Same for $\PivC^\gamma$.
\item[{\bf (ii)}]  Almost surely, $\PivO^\gamma(\omega_\infty^\lambda)\cup\PivC^\gamma(\omega_\infty^\lambda)=\Piv^\gamma(\omega_\infty^\lambda)$, a disjoint union.
\item[{\bf (iii)}] Almost surely, whenever $x\in\Piv^\gamma(\omega_\infty^\lambda)$ for some $\gamma$, the color of $x$ is the same for all such $\gamma$.
\end{itemize}
\epr

Note that (ii) is not a tautology (neither that the two colored sets are disjoint, nor that their union is the set of all the pivotals), since in $\omega_\infty^\lambda$ we did not define the set of closed pivotals as the complement of open pivotals. 

The main difficulty in proving~(i) is that the event $x\in\PivO^\gamma$ is not an open set in the quad-crossing topology $(\HH_D,\T_D)$: perturbing a configuration even by an arbitrary small amount may destroy a pivotal for $\gamma$, making the 4-arm event happen only from a strictly positive distance $\eps>0$ to $\gamma$. In terms of discrete percolation configurations, if there is an open pivotal connecting two halves of a cluster, then making the connection between the two halves a bit thicker is a small change w.r.t.~the quad-crossing topology, but it kills the pivotal. In particular, the harder direction in (i) will be to prove that there are ``enough'' pivotals in $\omega^\lambda_\infty$, since this requires controlling all scales simultaneously. 

\proof  For (i), we need to prove that for any $\eps>0$, if $\eta>0$ is small enough, then with probability at least $1-\eps$, for every $x_{\eta} \in \PivO^\gamma(\omega_{\eta}^\lambda)$ there exists some $x\in\PivO^\gamma(\omega_\infty^\lambda)$ within distance $\eps$ from $x_\eta$, and vice versa, for every $x\in\PivO^\gamma(\omega_\infty^\lambda)$ there exists $x_{\eta} \in \PivO^\gamma(\omega_{\eta}^\lambda)$.

There will be two key ingredients. Firstly, for any small $\alpha,\eps>0$ there exists $\delta,\bar\eta>0$ such that for all $0<\eta<\bar\eta$,
\be\label{e.margin}
\Pb{\PivO^{\gamma,\eps}(\omega^\lambda_\eta) = \PivO^{\gamma,\eps,\delta}(\omega^\lambda_\eta)} > 1-\alpha\,.
\ee
The existence of a $\delta$ that still depends on $x\in \PivO^{\gamma,\eps}(\omega^\lambda_\eta)$, or rather on its lattice square $B_\eps(x)$, is just  a special case of  \cite[Corollary 2.10]{GPS2a}. Then, taking the probability $\alpha$ of the error much smaller than $\eps^2$, we can find a $\delta>0$ that, with large probability, works for all points in $\PivO^{\gamma,\eps}(\omega^\lambda_\eta)$ simultaneously, proving~(\ref{e.margin}).

The point of introducing the $\delta$ margin of safety is that now~(\ref{e.margin}) immediately implies that there exists some monotone function $f=f_{\alpha,\eps}: [0,\infty)\lora [0,\infty)$ that could be described using the dyadic uniformity structures of \cite[Lemma 2.5]{GPS2a} and \cite[Proposition 3.9]{DPSL}) such that
\be\label{e.wobble}
\begin{aligned}
\P\big[ \forall\, x \in \PivO^{\gamma,\eps}(\omega^\lambda_\eta) \text{ and }\forall\, \tilde\omega\in \HH_D\text{ with } d_\HH(\tilde\omega,\omega^\lambda_\eta) < f(\delta),&  \\
\text{ we have } x\in\PivO^{\gamma,\eps,\delta/2}(\tilde\omega)\big] & > 1-\alpha\,,
\end{aligned}
\ee
for some $\delta>0$ and any $0<\eta<\bar\eta$, as given by~(\ref{e.margin}). 

The second key ingredient is that for any small $\alpha,\beta>0$, if $\eps,\hat\eta>0$ are small enough, then
\be\label{e.close}
\Pb{\forall\,x\in\PivO^{\gamma,\eps}(\omega^\lambda_\eta) \ \exists \, \tilde x\in\PivO^{\gamma}(\omega^\lambda_\eta)\text{ with }d(\tilde x,x)<\beta} > 1-\alpha
\ee 
for all $0<\eta<\hat\eta$. Before proving this, let us see how (\ref{e.wobble}) and~(\ref{e.close}) imply item~(i). We start with the first direction.

Fix $\alpha,\beta>0$ small.
Corresponding to them,~(\ref{e.close}) gives some $\eps_0,\hat\eta_0 >0$.
Now, corresponding to $\alpha$ and this $\eps_0$, there are $\delta_0,\bar\eta_0>0$ given by (\ref{e.wobble}). Take $0<\eta_0<\bar\eta_0\wedge\hat\eta_0$ so small that 
\be\label{e.converg}
\Pb{d_\HH(\omega^\lambda_\eta,\omega^\lambda_\infty)<f(\delta_0)/2\text{ holds for all }\eta\leq \eta_0} > 1-\alpha
\ee
in the coupling $\omega_{\eta}^\lambda\xrightarrow{a.s.}\omega_\infty^\lambda$ that we have. If the event of~(\ref{e.converg}) holds, then we have $d_\HH(\omega^\lambda_\eta,\omega^\lambda_{\eta_0})<f(\delta_0)$ for all $\eta\leq\eta_0$, and hence, together with (\ref{e.wobble}) and~(\ref{e.close}), we get
$$
\Pb{\forall\,x\in\PivO^{\gamma,\eps_0}(\omega^\lambda_{\eta_0}) \ \exists \, \tilde x\in\PivO^{\gamma}(\omega^\lambda_\eta)\text{ with }d(\tilde x,x)<\beta} > 1-3\alpha\,.
$$

Similarly, for $k\geq 1$, corresponding to $\alpha/2^k$ and $\beta/2^k$, there are $\eps_k,\hat\eta_k>0$ given by~(\ref{e.close}); we can make sure that $\eps_k < \eps_{k-1}/2$. Then, corresponding to $\alpha/2^k$ and $\eps_k$, there are $\delta_k,\bar\eta_k>0$ given by (\ref{e.wobble}). Take $0<\eta_k<\eta_{k-1}/2 \, \wedge\, \bar\eta_k \, \wedge \, \hat\eta_k$ so small that $d_\HH(\omega^\lambda_\eta,\omega^\lambda_\infty)<f(\delta_k)/2$ is satisfied for all $\eta\leq\eta_k$ with probability at least $1-\alpha/2^k$. Then, for all $\eta<\eta_k$, (\ref{e.wobble}) and~(\ref{e.close}) together give that 
$$
\Pb{\forall\,x\in\PivO^{\gamma,\eps_k}(\omega^\lambda_{\eta_k}) \ \exists \, \tilde x\in\PivO^{\gamma}(\omega^\lambda_\eta)\text{ with }d(\tilde x,x)<\beta/2^k } > 1-3\alpha/2^{k}\,.
$$

Iterating this procedure, we get that there exist sequences $\eta_k\to 0$ and $\eps_k\to 0$ such that with probability at least $1-3\alpha\sum_{k\geq 0} 2^{-k}=1-6\alpha$, for any $x_0 \in\PivO^{\gamma,\eps_0}(\omega^\lambda_{\eta_0})$ there exist 
\be\label{e.seq}
x_k\in\PivO^{\gamma,\eps_k}(\omega^\lambda_{\eta_k})\text{ for }k=0,1,2,\dots,\text{ satisfying }d(x_{k+1},x_k)<\beta/2^k\,.
\ee
These points have a limit $x_k \to \tilde x_0$, which satisfies $d(x_0,\tilde x_0)<2\beta$. Unsurprisingly, we claim that $\tilde x_0\in \PivO^\gamma(\omega^\lambda_\infty)$. Indeed, otherwise there would exist some $\tilde\eps>0$ such that $\tilde x_0 \not\in \PivO^{\gamma,\tilde\eps} (\omega^\lambda_\infty)$, but for some small enough $\eps$, this would clearly contradict the existence of an $\omega^\lambda_{\eta_k}$ satisfying $d_\HH(\omega^\lambda_{\eta_k},\omega^\lambda_\infty)<\eps$ and having an almost-pivotal $x_k \in\PivO^{\gamma,\eps}(\omega^\lambda_{\eta_k})$ at distance $d(x_k,\tilde x_0)<\eps$, which we have from~(\ref{e.seq}). Since we can take $\alpha$ and $\beta$ arbitrarily small, this finishes the proof of the first direction of item~(i).

For the other direction, if $x\in\PivO^\gamma(\omega_\infty^\lambda)$, then, by definition, for all $\eps>0$ with $B_\eps(x)$ surrounded by $\gamma$, there is some $\delta>0$ such that $x\in\PivO^{\gamma,\eps,\delta}(\omega_\infty^\lambda)$. Now, if $\omega_\eta^\lambda$ is close enough to $\omega_\infty^\lambda$ (again quantifiable in the sense of dyadic uniformity structures), then $x\in\PivO^{\gamma,\eps,\delta/2}(\omega_\eta^\lambda)$ also occurs. By~(\ref{e.close}), if $\eps>0$ is small enough, this means with large probability that there is an actual pivotal of $\omega_\eta^\lambda$ close to $x$, as required.

\begin{figure}[htbp]
\SetLabels
(0.1*0.08)\small{$\gamma$}\\
(0.2*0.21)\small{$\beta$}\\
(0.22*0.36)\small{$\eps$}\\
\endSetLabels
\centerline{
\AffixLabels{
\includegraphics[width=\textwidth]{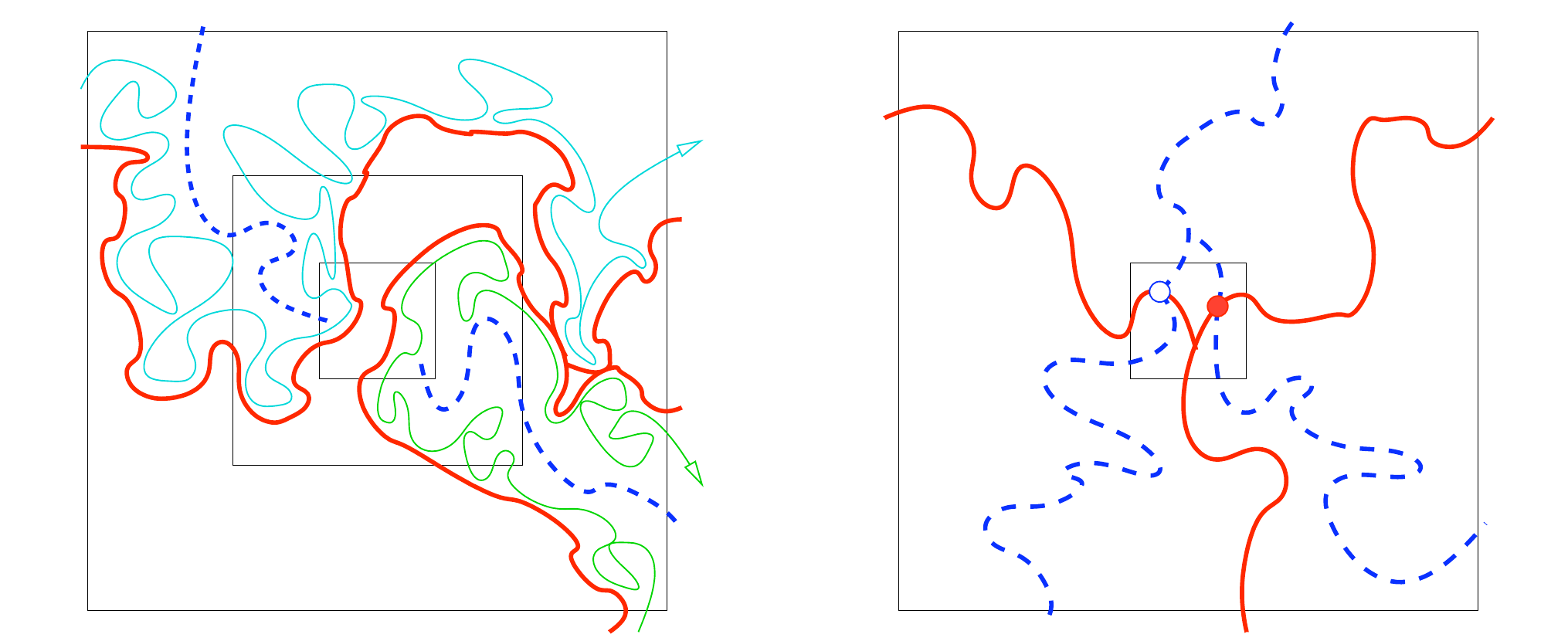}
}}
\caption{Left: an open $\eps$-almost pivotal event without actual pivotals in the $\beta$-square implies a 6-arm event (four open and two closed arms) between radii $\eps$ and $\beta$. Right: having both an open and a closed $\gamma$-pivotal in an $\eps$-square implies a 6-arm event from $\eps$ to $\gamma$.} 
\label{f.PivExists}
\end{figure}

We still owe the proof of~(\ref{e.close}). Assume that $x\in\PivO^{\gamma,\eps}(\omega^\lambda_\eta)$ but there are no open pivotal sites in $B_{\beta}(x)$. This implies that there is a 6-arm event from $\p B_\eps(x)$ to $\p B_{\beta}(x)$: the interfaces between the open and closed arms cannot touch each other within $B_{\beta}(x)$, hence their open sides form two disjoint open paths, creating four open arms besides the two closed ones; see the left side of Figure~\ref{f.PivExists}. 
Since the 6-arm exponent is strictly larger than 2 at any fixed near-critical level $\lambda$ (see \cite[Corollary A.8]{SchSt} for $\lambda=0$, and \cite[Proposition 11.6]{DPSL} or Proposition~\ref{p.stab} in the present paper for general $\lambda$), we can take $\beta:=\eps^\zeta$ with $\zeta>0$ small enough for the polychromatic 6-arm probability satisfy $\alpha_6(\eps,\beta)=o(\eps^2)$, and then the probability that such a 6-arm event occurs anywhere in the domain tends to zero as $\eps\to 0$, and we are done.

In item (ii), the fact that the union of the two colored sets gives all the pivotals follows immediately from the discrete analogue and item (i). To prove the disjointness claim, by part (i) it is enough to prove that the probability of having a closed and an open pivotal for $\gamma$ within distance $\eps$ from each other goes to 0 as $\eps\to 0$. But this event implies the existence of  6 disjoint arms from $\eps$ to $\gamma$ (see the right side of Figure~\ref{f.PivExists}), and hence, as usual, the 6-arm exponent being larger than 2 implies the claim.

For item (iii), if $\gamma_1$ and $\gamma_2$ both surround $x$, with $x\in \PivO^{\gamma_1}(\omega^\lambda_\infty) \cap  \PivC^{\gamma_2}(\omega^\lambda_\infty)$, then we would also have $x\in \PivO^{\gamma}(\omega^\lambda_\infty) \cap  \PivC^{\gamma}(\omega^\lambda_\infty)$, where $\gamma$ is the part of $\gamma_1\cup\gamma_2$ that is visible from $x$. However, this is impossible by item (ii). 
\QED
\medskip

Beyond the set of pivotals, we are also interested in the normalized counting measure on them. In \cite[Subsection~2.6]{DPSL}, for any fixed $\eps>0$, we defined the set of {\bf $\eps$-important points} $\Pimp^\eps(\omega_\eta)$ of any discrete percolation configuration in a bounded domain $D\subset \hat\C$, relative to the $(\eps,3\eps)$-annuli given by a fixed lattice $\eps\Z^2$. Namely, for $x\in D$, we let $B_\eps(x)$ be the lattice square as before, let $\tilde B_\eps(x)$ be the $3\eps$-square centered at $B_\eps(x)$, and let $x\in\Pimp^\eps$ if{f} $x\in\Piv^{\p \tilde B_\eps(x)}$. Then we considered the normalized counting measure  $\mu^\eps(\omega_\eta)$ on this set $\Pimp^\eps$. Of course, the same discrete definition works for near-critical percolation configurations $\omega_\eta^\lambda$. Then, the main result of  \cite{GPS2a} is the following convergence of $\mu^\eps$ for $\lambda=0$, extended to general $\lambda\in\R$ by \cite[Theorem 11.5]{DPSL}:

\begin{theorem}\label{t.PivMeas}
For any $\eps>0$, there exists a measurable map $\mu^\eps: \HH_D \lora \mathfrak{M}_D$, into the space of finite Borel measures on $D$, such that, for $\lambda\in\R$, as $\eta\to 0$,
\begin{equation*}
(\omega_\eta^\lambda, \mu^\eps(\omega_\eta^\lambda)) \overset{d}{\longrightarrow} (\omega_\infty^\lambda, \mu^\eps(\omega_\infty^\lambda))
\end{equation*}
in the quad-crossing topology $(\HH_D,\T_D)$ in the first coordinate and in the L\'evy-Prokhorov distance of measures in the second one. Furthermore, the above Proposition~\ref{p.PivSetCol} implies immediately the convergence 
$$
\big(\PimpO^\eps(\omega_\eta^\lambda),\PimpC^\eps(\omega_\eta^\lambda)\big) \overset{d}{\longrightarrow} 
\big(\PimpO^\eps(\omega_\infty^\lambda),\PimpC^\eps(\omega_\infty^\lambda)\big)
$$ 
in the Hausdorff metric of closed sets.
\end{theorem}

\section{Enhanced networks and cut-off forests built from the near-critical ensemble}\label{s.enhanced}

The pivotal measures of \cite{GPS2a} that we recalled in Theorem~\ref{t.PivMeas} were used in \cite{DPSL} as the intensity measures for the Poisson point processes of pivotal sites that switch as the near-critical parameter $\lambda\in\R$ changes. Here is the exact notation that we will use:

\begin{definition}\label{d.PPP}
Let $\vl=(\lambda,\lambda')\in\R^2$ be any pair of near-critical parameters with $\lambda<\lambda'$, and let $\eps>0$ be fixed. Let $\omega^\lambda$ be a near-critical configuration $\omega_\eta^\lambda$ or $\omega_\infty^\lambda$ in $\Tor_M$. We will denote by $\PPP^\eps_\vl=\PPP_\vl(\mu^\eps(\omega^\lambda))$  the Poisson point process 
$$
\PPP^\eps_\vl = \{ (x_i,t_i), 1\le i \le p \} \subset \Pimp^\eps(\omega^\lambda) \times [\lambda,\lambda']
$$
of intensity measure $\mu^\eps(\omega^\lambda)(dx)\times \1_{[\lambda,\lambda']}(t)\, dt $. The set $\{x_1,\dots,x_p\}$ of 
pivotals will usually be denoted by $X$.
For the case of $\omega^\lambda_\eta$, the process $\PPP^\eps_\vl$ can clearly be constructed measurably from $\omega_\eta^{[\lambda,\lambda']}$, and we will always work in this natural coupling.
\end{definition}

In Section 6 and Subsection 11.2 of \cite{DPSL}, for any quad $\Quad\subset \C$, any $\eps>0$, any discrete or continuum near-critical percolation configuration $\omega^\lambda$ and the associated Poisson point process $\PPP^\eps_\vl(\omega^\lambda)$, we constructed an edge-colored graph $\Net_Q(\omega^\lambda,\PPP^\eps_\vl)$, called an {$\eps$-network}, whose vertex set was the Poisson point set $X=\{x_1,\dots,x_p\}$ of pivotals together with the four boundary arcs of $\Quad$, and whose edge set was given by the primal and dual connections in $\omega^\lambda$ between the vertices. Since in this paper we are primarily interested in spanning trees, not in quad-crossings, it will be useful to change the boundary conditions in the definition slightly (but still using the quad-crossing topology). We will also need to add a bit more structure to these networks: roughly, we will need to know which pivotals in $X$ are connected together by an open cluster of $\omega^\lambda\setminus X$, and will need to know the colors of these pivotals in $\omega^\lambda$. The resulting structures will be called enhanced networks. Just as in \cite{DPSL}, we start with the following simple definition: 

\begin{definition}[A nested family of dyadic coverings]\label{d.nested}
For any $b>0$ in $2^{-\N}=\{2^{-k} : k=0,1,2,\dots\}$, let $G_b$ be a disjoint covering of $\Tor_M$ using the lattice $b$-squares 
$\big\{ [k,k+1)b \times [\ell,\ell+1)b : (k,\ell)\in\Z^2 \big\}$. 
Now, for any $r\in 2^{-\N}$ and any finite subset $X=\{x_1,\ldots, x_p\} \subset \Tor_M$, one can associate uniquely $r$-squares 
$B_{x_1}^r,\ldots,B^r_{x_p}$ in the following manner: for all $1\le i \le p$, there is a unique square $\tilde B_{x_i} \in G_{r/2}$ which contains $x_i$ and we define $B_{x_i}^r$ to be the $r$-square in the grid $r\Z^2 - (r/4,r/4)$ centered around the $r/2$-square $\tilde B_{x_i}$. We will denote by $B^r(X)$ this family of $r$-squares. This family has the following two properties:
\bi
\item[(i)] Each point $x_i$ is at distance at least $r/4$ from $\p B_{x_i}^r$. 
\item[(ii)] For any set $X$, $\{ B^r(X) \}_{r\in 2^{-\N}}$ forms a nested family of squares in the sense that for any $r_1<r_2$ in $2^{-\N}$, and any $x\in X$, we have 
$
B^{r_1}_{x} \subset B^{r_2}_{x}
$.
\ei
\end{definition}

For a finite set of points $X\subset \Tor_M$, let $r^*(X)>0$ denote one-tenth of the smallest distance between any pair $x_i\not=x_j\in X$. With minor changes from the case of a domain with a boundary to the case of a torus, it is proved in \cite[Proposition 5.2]{DPSL} that for $X$ being the pivotals in $\PPP^\eps_\vl$, the random variable $r^*(\PPP^\eps_\vl)$ is almost surely positive (with a small abuse of notation, since $\PPP^\eps_\vl$ is a subset of $\Tor_M \times [\lambda,\lambda']$).

\bde\label{d.rNetwork}
For $0<r<r^*(\PPP^\eps_\vl)$, the {\bf $r$-mesoscopic $\eps$-network} $\Net^{r\meso}_M(\omega^\lambda,\PPP^\eps_\vl)$ associated to a near-critical percolation configuration $\omega^\lambda$ in the torus $\Tor_M$ and the Poisson point process $\PPP^\eps_\vl$ of Definition~\ref{d.PPP} is the graph with vertex set $\{x_1,\dots,x_p\}$ and two types of edges, labelled primal or dual, with a primal edge connecting $x_i$ and $x_j$ if there exists a quad $R$ such that $\p_1 R$ and $\p_3 R$ remain strictly inside $B_{x_i}^r$ and $B_{x_j}^r$, and $R$ remains strictly away from the squares $B_{x_k}^{r}, k \notin \{i,j\}$, and for which $\omega^\lambda \in \boxminus_R$. Dual edges are defined analogously (still w.r.t.~$\omega^\lambda$).

We consider two $r$-mesoscopic networks to be the same if the $r$-squares for the vertices (as embedded in $\Tor_M$) and the labelled graph structures coincide. For $r_1<r_2$, we can compare an $r_1$-mesoscopic network with an $r_2$-mesoscopic network by considering the unique $r_2$-squares containing the $r_1$-squares of the first network.
\ede


We will now take $r\to 0$, get a network $\Net_M(\omega^\lambda,\PPP^\eps_\vl)$, and then compare these networks for $\omega^\lambda_\eta$ and $\omega^\lambda_\infty$. The following results were proved in \cite[Theorem 6.14]{DPSL} and \cite[Subsection 7.4]{DPSL} for $\lambda=0$, extended to general $\lambda$ in \cite[Subsection 11.2]{DPSL}, for networks defined using slightly different boundary conditions than here, but with the same proofs working fine:

\bpr[$r$-stabilization and $\eta$-convergence of networks]\label{p.netconv}\ 
\bi
\item[{\bf (i)}]
There exists a measurable scale $0<r_M=r_M(\omega^\lambda_\infty,\PPP^\eps_\vl)<r^*(\PPP^\eps_\vl(\omega^\lambda_\infty))$ such that for all $r\in (0,r_M)$ we get the same $r$-mesoscopic $\eps$-network $\Net_M^{r\meso}(\omega^\lambda_\infty,\PPP^\eps_\vl)$. This stabilized network will be called the {\bf $\eps$-network} $\Net^{\vl,\eps}_\infty=\Net_M(\omega^\lambda_\infty,\PPP^\eps_\vl)$. For discrete percolation configurations, the definition of $\Net^{\vl,\eps}_\eta=\Net_M(\omega^\lambda_\eta,\PPP^\eps_\vl)$ is the obvious one.
\item[{\bf (ii)}]
For any $\alpha>0$ there is a scale $r_\alpha=r_\alpha(M,\vl,\eps)$ such that in any coupling with $\omega_\eta^\lambda \xrightarrow{a.s.}\omega_\infty^\lambda$ in $\Tor_M$, for all sufficiently small $\eta>0$ there is a coupling of $\PPP^\eps_\vl(\omega_\eta^\lambda)$ and $\PPP^\eps_\vl(\omega_\infty^\lambda)$ such that with probability at least $1-\alpha$ the following holds: $r_\alpha$ is less than both $r_M(\omega^\lambda_\infty,\PPP^\eps_\vl) < r^*(\PPP^\eps_\vl(\omega^\lambda_\infty))$ and $r^*(\PPP^\eps_\vl(\omega^\lambda_\eta))$, and for all $r<r_\alpha$ we have 
$$
\Net_M^{r\meso}(\omega^\lambda_\eta,\PPP^\eps_\vl(\omega^\lambda_\eta))=
\Net_M^{r\meso}(\omega^\lambda_\infty,\PPP^\eps_\vl(\omega^\lambda_\infty))\,;
$$
in this sense, $\Net^{\vl,\eps}_\eta=\Net_M(\omega_\eta^\lambda,\PPP^\eps_\vl)$ coincides with $\Net^{\vl,\eps}_\infty=\Net_M(\omega_\infty^\lambda,\PPP^\eps_\vl)$. (Only in this sense, not exactly, since the vertex sets $\PPP^\eps_\vl(\omega^\lambda_\infty)$ and $\PPP^\eps_\vl(\omega^\lambda_\eta)$ are only close to each other, but do not coincide.)
\ei
\epr

\begin{figure}[htbp]
\centerline{
\AffixLabels{
\includegraphics[width=0.8 \textwidth]{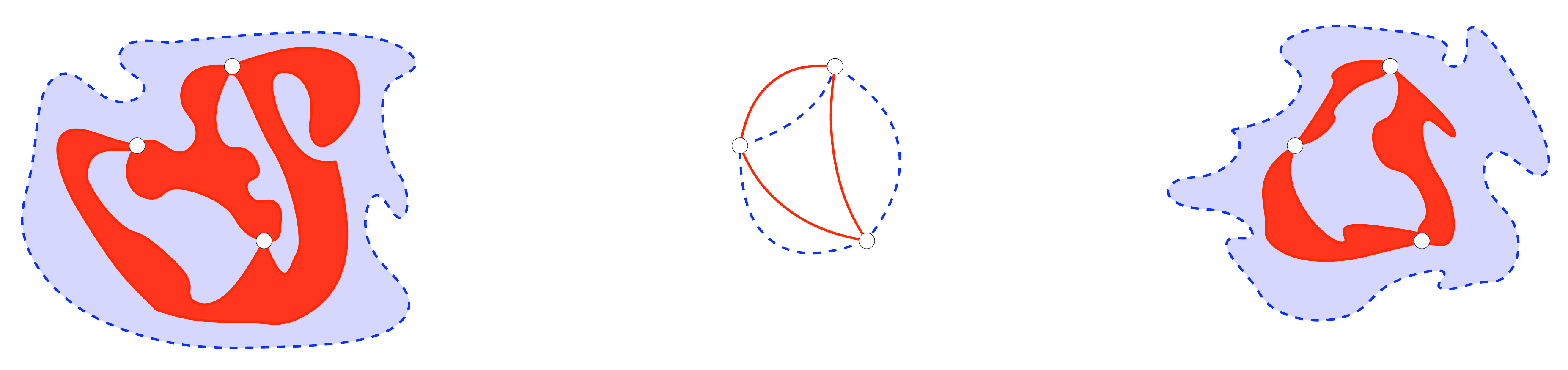}
}}
\caption{The same graph structure in a network (the middle picture) may correspond to very different cluster structures (on the two sides).}
\label{f.NeedRouters}
\end{figure}

Note that a network in itself may completely fail to describe the structure of clusters: see Figure~\ref{f.NeedRouters}. 
This is a bit of a problem for the purposes of the present paper, hence we are going to add some extra structure to our networks that will be measurable w.r.t.~the quad-crossing topology (in particular, it makes sense for $\omega^\lambda_\infty$), while it describes how the pivotals of $\PPP^\eps_\vl$ are connected to each other in $\omega^\lambda$.

\bde[Mesoscopic sub-routers]\label{d.mesosub}
Fix $0<r<\rho<\infty$. Utilizing the notation introduced in Definition~\ref{d.nested}, let $B^r(\Tor_M)$ be the finite covering of $\Tor_M$ by overlapping $r$-squares. Given a subset $Y$ of the set $X=\{x_1,\dots,x_p\}$ of the pivotals in $\PPP^\eps_\vl$, with $|Y|\geq 2$, an {\bf $(r,\rho)$-mesoscopic sub-router for $Y$} is an $r$-square $B\in B^r(\Tor_M)$ with the following properties:
\bi
\item it is at distance at least $2\rho$ from each $x_i\in X$;
\item  there is an open circuit (i.e., no dual arm) in the square annulus with inner face $B$ and outer radius $\rho$; the largest $s$-square with some $s\in 2^{-\N}$ that is concentric with $B$, contains it, and is surrounded by the open circuit will be denoted by $\widehat B$;
\item for each $x_i\in Y$, there exists a quad $R$ with $\p_1 R$ contained in $\widehat B$, $\p_3 R$ contained in $B^r_{x_i}$, remaining strictly away from all the squares $B^{r}_{x_k}$ with $x_k \in X\setminus\{x_i\}$, and for which $\omega^\lambda \in \boxminus_R$.
\ei
\ede

Let $\RR_Y(B)$ denote the event that an $r$-square $B$ is an $(r,\rho)$-mesoscopic sub-router for some $Y\subseteq X$. This is measurable w.r.t.~$\omega^\lambda$, and using Lemmas~\ref{l.crossconv} and~\ref{l.arms}, in the coupling of Proposition~\ref{p.netconv}~(ii), the set of $r$-squares $B$ for which $\RR_Y(B)$ holds in $\omega^\lambda_\eta$ is the same with probability tending to 1 (as $\eta\to0$) as in $\omega^\lambda_\infty$. Furthermore, by choosing $(r,\rho)$ appropriately, this set will turn out to be non-empty with high probability, for all possible $Y$. For this, a key proposition, interesting in its own right, is the following:

\bpr[The volume of clusters]\label{p.clusterdim}
For any $\lambda\in\R$, $M>\rho>0$ and $\zeta>0$ fixed, for percolation $\omega^\lambda_\eta$ in $\Tor_M$, with probability tending to 1 as $\eta\to 0$,  all clusters of diameter at least $\rho$ have at least $(\rho/\eta)^{91/48-\zeta}$ sites. (Note that $91/48$ equals 2 minus the one-arm exponent 5/48 \cite{LSW:1arm}.)

Similarly, with probability tending to 1 as $r\to 0$, uniformly in the mesh $\eta$, all these clusters have a ``large $r$-volume'' in the following sense: the number of $r$-squares in $B^r(\Tor_M)$ that intersect the cluster is at least $(\rho/r)^{91/48-\zeta}$.
\epr

After the first version of this paper was posted, Rob van den Berg pointed out that this proposition follows from~(3.15) of \cite{Jarai}. However, since the proof there is hard to read, we decided to keep our proof for the sake of completeness. Earlier, similar but weaker results were proved in \cite[Lemma 3.20]{KestenRWIIC} and \cite[Theorem 3.3]{BCKS01}. Finally, \cite[Lemma 2.7]{vdBC} gives a bit more elegant version of our argument, but proving a little less; in particular, it is not proved there that all the radial crossings of a $(\rho/3,\rho)$-annulus are everywhere well-separated from each other (see our proof below).

\proof The proof will rely only on multi-arm exponents, hence, in view of Proposition~\ref{p.stab}, the reader may just think of $\lambda=0$. We will do the case of the standard volume (number of sites in the $\eta$-mesh); the proof works the same way for the case of the $r$-volume. 

Take the lattice $(\rho/3) \Z^2$, and centered around each $\rho/3$-square, consider the square of side-length $\rho$ and the annulus between these two square boundaries. It is easy to check that any cluster of diameter at least $\rho$ produces a radial crossing of such a $(\rho/3,\rho)$-annulus. The number of such annuli is $\asymp (M/\rho)^{2}$.

Whether a given $(\rho/3,\rho)$-annulus $A_\rho$ is radially crossed can be decided using the radial exploration process started at any point along the boundary at radius $\rho/3$, with open hexagons on the right side, closed hexagons on the left, stopped when reaching the boundary at radius $\rho$. (See around Figure 2.6 of \cite{GPS2a} or \cite[Section 4.3]{WWperc} for the definition of this exploration process.) If the annulus is crossed, there are two cases: either (a) there is also an open circuit, or (b) there is also at least one radial dual crossing.

\begin{figure}[htbp]
\SetLabels
(0.15*0.36)$\rho/3$\\
(0.04*0.07)$\rho$\\
(0.27*0.85)\textcolor{MyMagenta}{$\tilde A_i \subset \tilde B_i$}\\
(0.26*0.65)\textcolor{Red}{$\Gamma$}\\
(0.29*0.55)\textcolor{MyMagenta}{$\tilde A_{\Gamma,\delta}$}\\
(0.57*0.85)\textcolor{Red}{$\Gamma$}\\
(0.85*0.13)\textcolor{Red}{$\Gamma$}\\
\endSetLabels
\centerline{
\AffixLabels{
\includegraphics[width=1.1 \textwidth]{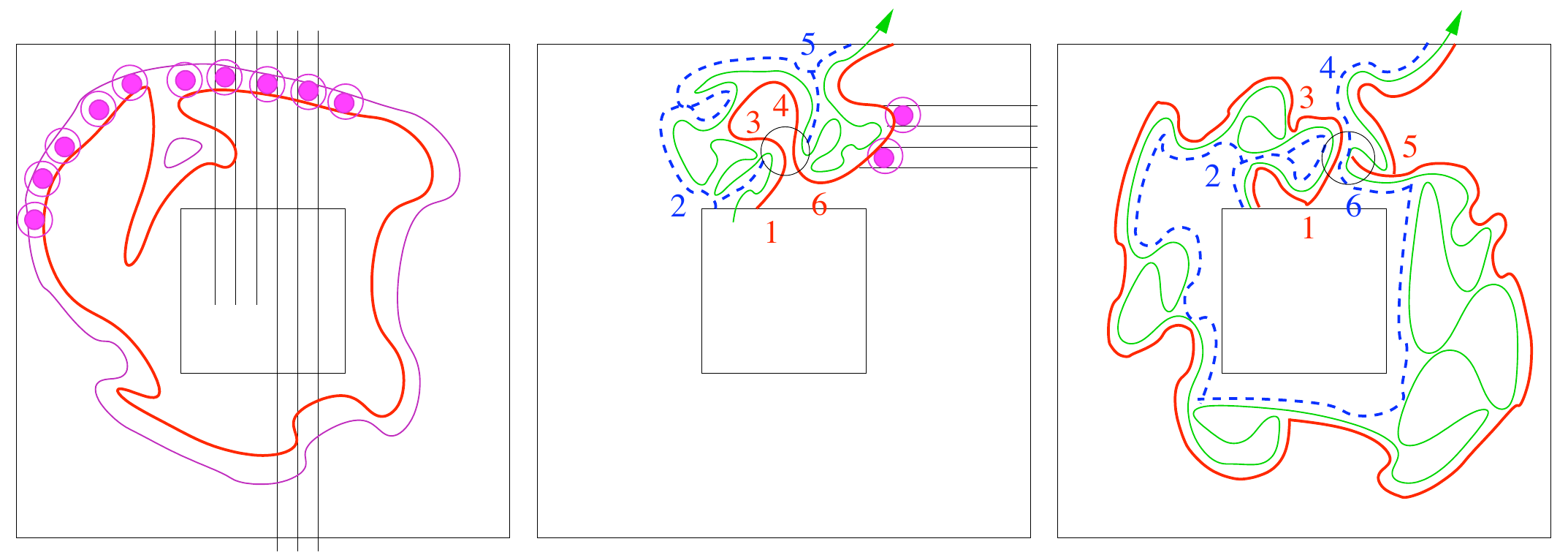}
}}
\caption{If the annulus $A_\rho$ has an open circuit or is crossed radially, then the radial exploration process gives an open path $\Gamma$ that has macroscopically wide unexplored space on one side, collecting large enough volume connected to $\Gamma$ with high probability.}
\label{f.ClusterVol}
\end{figure}
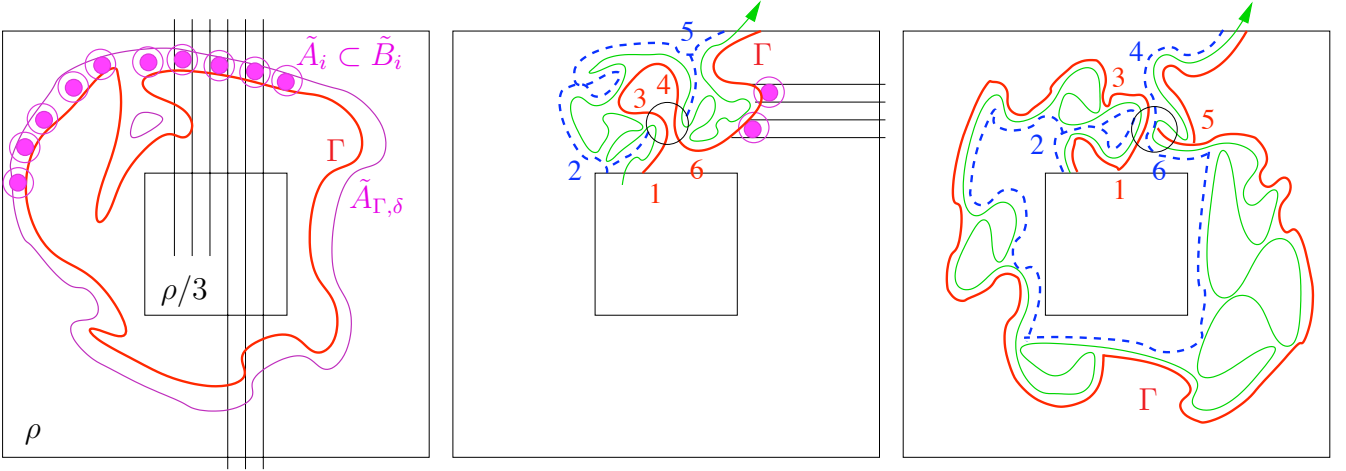

{\bf (a)} Condition on having an open circuit; this is slightly more general than the first of the two above cases, since we do not condition on having also a radial crossing. Condition on the smallest open circuit, $\Gamma$. The radial exploration process finds it from inside, hence the configuration in the annulus between $\Gamma$ and $\p_2 A_\rho$, denoted by $A_\Gamma$, is undisturbed percolation. Moreover, by the half-plane 3-arm exponent being 2, the probability that the distance between $\Gamma$ and $\p_2 A_\rho$ is smaller than $\delta\rho$ is $O(\delta)$. Let this distance be the random variable $\delta_\Gamma\rho$, take any $0<\delta<\delta_\Gamma$, and take the set of points of $A_\Gamma$ whose distance from $\Gamma$ is less than $\delta\rho$. It is clear that this set, denoted by $\tilde A_{\Gamma,\delta}$, contains a collection of $K=K(\delta)\geq c/\delta$ disjoint balls of diameter $\delta\rho$, denoted by $\tilde A_i$, $i=1,\dots, K$, such that all their pairwise distances are at least $\delta\rho$; for instance, take a family of vertical parallel lines with mesh $\delta\rho$, and in every other slab, take the uppermost ball of diameter $\delta\rho$ that touches $\Gamma$. See the first picture in Figure~\ref{f.ClusterVol}. We will still fine-tune the value of $\delta$ later.

If a site in some $\tilde A_i$ has an open arm to distance at least $c \delta\rho$, then with a uniformly positive probability it is connected to $\Gamma$, within the $\delta\rho/2$-neighborhood of $\tilde A_i$ that will be denoted by $\tilde B_i$. Vice versa, most sites in $\tilde A_i$ need to have an arm of length $c\delta\rho$ in order to be connected to $\Gamma$. Thus, letting $X_i(\delta)$ be the number of sites in $\tilde A_i$ that are connected to $\Gamma$ within $\tilde B_i$, and using quasi-multiplicativity of the one-arm probability $\alpha_1(\cdot,\cdot)$ (see Remark~\ref{r.arms}), we have 
$$
\E^\lambda_\eta[ X_i(\delta) ]\asymp (\delta\rho/\eta)^2 \alpha_1(\eta,\delta\rho)=(\delta\rho/\eta)^{91/48+o(1)}.
$$
It is a standard argument using quasi-multiplicativity and a summation over dyadic scales that the second moment of $X_i$ is comparable to the square of the first moment (see, e.g., \cite[Lemma 3.1]{GPS1} for the second moment of the number of pivotals). Thus, by the Paley-Zygmund second moment inequality (a simple consequence of Cauchy-Schwarz; see, e.g., \cite[Section 5.5]{LPbook}), there exists a uniform constant $c=c_\lambda>0$ such that 
$$
\P^\lambda_\eta\big[X_i(\delta) > c\,\E^\lambda_\eta X_i(\delta)\big]>c\,.
$$ 
Using the independence of the variables $X_i$ (conditionally on $\Gamma$) that follows from the disjointness of the neighborhoods $\tilde B_i$, we get that 
\be\label{e.VolumeTail}
\begin{aligned}
\P^\lambda_\eta\Big[\text{cluster of }\Gamma\text{ has volume}\leq c\,\E^\lambda_\eta[ X_i(\delta) ]
 \,\Big|\, \Gamma \Big] & \\
& \hskip -5 cm \leq \P^\lambda_\eta\Big[ X_i < c\,\E^\lambda_\eta[ X_i(\delta) ]
\text{ for all }i=1,\dots,K(\delta) \,\Big|\, \Gamma \Big] \\
& \hskip -5 cm \leq (1-c)^{K(\delta)} = \exp(-c'/\delta)\,.
\end{aligned}
\ee

Now we want to choose $\delta$ such that the bound $c\,\E^\lambda_\eta[ X_i(\delta) ] = (\delta\rho/\eta)^{91/48+o(1)}$ on the cluster size becomes at least $(\rho/\eta)^{91/48-\zeta}$. This means $\delta=(\rho/\eta)^{-48 \zeta/91+o(1)}$, but this choice is allowed only if this value is less than $\delta_\Gamma$. As mentioned above, this fails with probability $(\rho/\eta)^{-48 \zeta/91+o(1)}$, which, for $\eta$ small enough, is much smaller than $(\rho/M)^2$. Therefore, with probability tending to 1 as $\eta\to 0$, in all the at most $O((M/\rho)^{2})$ annuli where case (a) occurs, $\delta_\Gamma$ is large enough and the event of~(\ref{e.VolumeTail}) fails to hold, hence the cluster of $\Gamma$ has volume at least $(\rho/\eta)^{91/48-\zeta}$.

{\bf (b)} Condition on the second case, and let $\Gamma$ be the clockwisemost radial open crossing that the exploration process has found. We claim that, similarly to case (a), there is a random variable $\delta_\Gamma$, uniformly positive in $\eta$, such that no hexagons have been explored in the clockwise $\delta_\Gamma\rho$-neighborhood of $\Gamma$. Indeed, this was already used in \cite[Lemma 2.9]{GPS2a} in the proof of the quad-measurability of the 1-arm event, and the reason is simply that this maximal distance $\delta_\Gamma$ can be less than some $\delta>0$ only if the radial exploration path comes to distance $\delta\rho$ to itself without touching, which would imply a full plane 6-arm event from distance $\delta\rho$ to  distance of order $\rho$ (or a half-plane 3-arm event, if it happens close to one of the boundary components of $A_\rho$). See the second and third pictures in Figure~\ref{f.ClusterVol}. Now, we can repeat the rest of the proof of case (a) within this unexplored space of width $\delta_\Gamma\rho$, and we are almost done: we have just proved that, with very high probability as $\eta\to 0$, the cluster found by the radial exploration process started at some arbitrary (say, uniform random) point at radius $\rho/3$ has large volume. However, we want this for {\it all} clusters that cross $A_\rho$, while the above procedure finds larger clusters with larger probability. 

\begin{figure}[htbp]
\SetLabels
(0.41*0.15)$\rho/3$\\
(0.92*0.15)$\rho$\\
\endSetLabels
\centerline{
\AffixLabels{
\includegraphics[width=0.5 \textwidth]{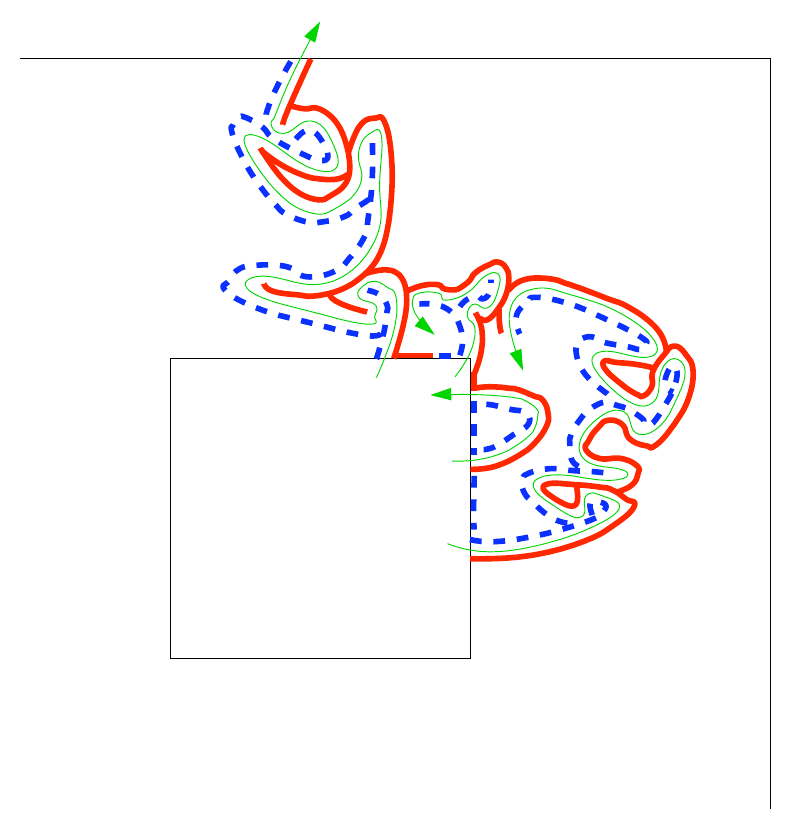}
}}
\caption{Consecutive radial exploration processes.}
\label{f.Bumps}
\end{figure}

To this end, once we have found one crossing cluster, we start a new radial exploration from radius $\rho/3$, at the first point on $\p B_{\rho/3}$ to the right of the last boundary touching point of the first exploration path that has an open site on the right and a closed site on the left side. We stop the process either when it reaches an open site explored by the previous exploration path and hence turns inside, towards $\p B_{\rho/3}$, or when it reaches $\p B_{\rho}$ (which we may call a ``success''). Then we take the next point on $\p B_{\rho/3}$ that has an open site on the right and a closed site on the left side, and so on, until the entire boundary $\p B_{\rho/3}$ has been explored and hence all radially crossing clusters have been found. Now, before each success, the right boundary of what has been built by the sequence of unsuccessful explorations is an open arm from $\p B_{\rho/3}$ to $\p B_{\rho}$, and from each point of this open arm, there is also a closed arm to $\p B_{\rho/3}$. Therefore, if the next successful exploration path comes $\delta\rho$-close to this right boundary, then it creates a full plane 6-arm or a half-plane 3-arm event (the third picture of Figure~\ref{f.ClusterVol} applies locally), which do not happen anywhere in $\A_\rho$ if $\delta$ is small enough. Therefore, all these right boundaries have the open unexplored space to their right that is required for our argument to work. Since each radially crossing cluster has, as a subset, such a right boundary (not necessarily the right boundary of the entire cluster), the proof of Proposition~\ref{p.clusterdim} is complete.
\qed
\medskip

We can now prove that the $(r,\rho)$-mesoscopic sub-routers of Definition~\ref{d.mesosub} exist:

\begin{lemma}\label{l.subrouters}
With probability tending to 1 as $\eta\to 0$ and then $r\to 0$, for any $0<\rho<r^*(X)$, for all $Y\subseteq X$ with $|Y|\geq 2$ whose points are connected together by a single open cluster of $\omega_\eta^\lambda\setminus X$ (more precisely, there is a cluster  of $\omega_\eta^\lambda\setminus X$ that neighbours each hexagon in $Y$), the set of $(r,\rho)$-mesoscopic sub-routers for $Y$ is non-empty.
\end{lemma}

\begin{figure}[htbp]
\SetLabels
(0.07*0.05)$B$\\
(0.*0.8)$\rho \wedge \eps$\\
(0.16*0.7)$r$\\
(0.84*0.7)$r$\\
\endSetLabels
\centerline{
\AffixLabels{
\includegraphics[width=0.5 \textwidth]{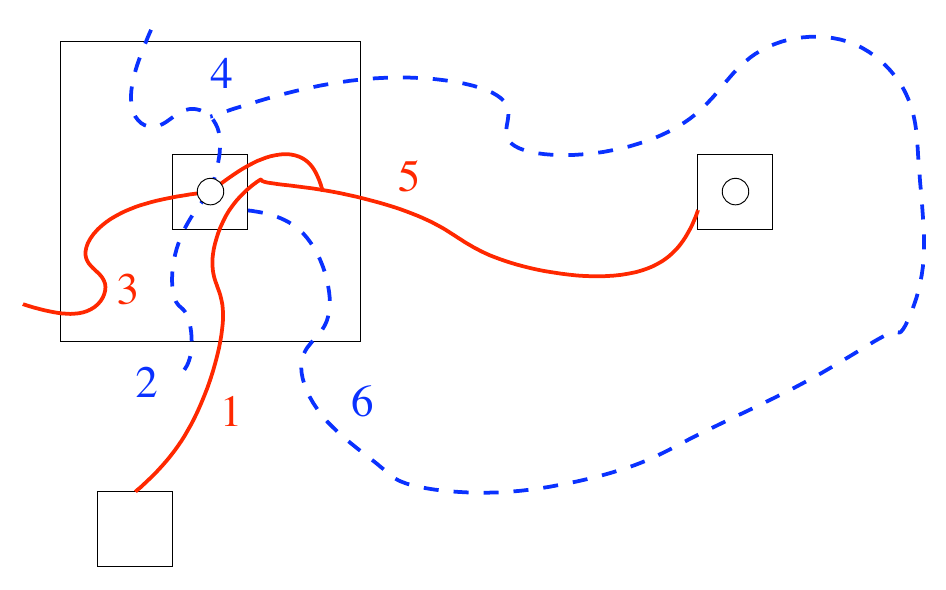}
}}
\caption{Connections from possible sub-routers $B$ can avoid other $r$-squares $B^r_{x_k}$ unless a 6-arm event happens.}
\label{f.subrouter}
\end{figure}

\proof Assume that in a configuration $\omega_\eta^\lambda$, some $Y\subseteq X$ satisfies the above conditions. Let $\rho$ be less than $r^*(X)$, take $r\ll \rho$, and consider any $r$-square $B$ that intersects the cluster and whose distance from $B_r(X)$ is at least $\rho$. By the definition of $r^*(X)$ and by $|Y|\geq 2$, such a $B$ certainly exists. We are going to examine when such a $B$ could be a mesoscopic sub-router. Simply aiming at $\widehat B=B$, the required quad connecting $B$ with an $x_i\in Y$ can fail to exist only if all the connections from $B$ to $B^r_{x_i}$ are $r$-close to some $x_k \in X\setminus\{x_i\}$; however, this would imply a 6-arm event from radius $r$ to $\rho\wedge \eps$ (see Figure~\ref{f.subrouter}), which does not occur anywhere in $\Tor_M$ if $r$ is small enough.

 We still need to show that, among the $r$-squares $B$ as above, there is at least one that also has the open circuit in the $(r,\rho)$-annulus around it.

 If $\rho < r^*(X)$, then any cluster $\calC$ connecting the points of $Y$ has a connected subset $\calC'$ of diameter at least $\rho$ that has a distance at least $\rho$ from all points of $X$. (We used here the definition of $r^*(X)$ and that $|Y|\geq 2$.) For the maximal such $\calC'$, the proof of Proposition~\ref{p.clusterdim} clearly applies, and for $r\ll\rho$, the number of $r$-squares in $B^r(\Tor_M)$ intersected by $\calC'$ is at least $(\rho/r)^{91/48-\zeta}$ with probability tending to $1$ as $r\to 0$. On the other hand, any of these $r$-squares fails to be an $(r,\rho)$-mesoscopic sub-router only if there is no open circuit in the $(r,\rho)$-annulus around $B$. In such a case, we have both a primal and a dual arm in the $(r,\rho)$-annulus, which event has probability $(r/\rho)^{1/4+o(1)}$, uniformly in $\eta>0$, by the 2-arm exponent \cite{SW}. Thus the number of such $r$-squares is $(\rho/r)^{7/4+o(1)}$ in expectation, and by Markov's  inequality, it is unlikely to be much larger, for any of the possible subsets $Y\subseteq X$ (whose number is independent of $r$). Since $(\rho/r)^{7/4+o(1)}$ is negligible compared to the $r$-volume $(\rho/r)^{91/48-\zeta}$ if $\zeta>0$ is small enough, with probability going to 1 as $r\to 0$, we do have $(r,\rho)$ sub-routers in every cluster spanned by some $Y\subseteq X$.
\qed
\medskip

If $B_1$, $B_2$ are $(r,\rho)$ sub-routers for $Y_1,Y_2\subseteq X$, respectively, we will call them connected if there exists a quad $R$ with $\p_1 R$ contained in $\widehat{B_1}$, $\p_3 R$ contained in $\widehat{B_2}$, remaining strictly away from all the squares $B^{r}_{X}$, and for which $\omega^\lambda \in \boxminus_R$. As before, in the coupling of Proposition~\ref{p.netconv}~(ii), for $\rho<r_M$, the relation of being connected converges in probability as $\omega_\eta\to \omega_\infty$, which also implies that it is an equivalence relation even in $\omega_\infty$. If $B_i$ is an $(r,\rho)$ sub-router for $Y_i\subseteq X$, $i=1,2$, and $B_1$ and $B_2$ are connected, then both $B_i$'s are $(r,\rho)$ sub-routers for $Y_1\cup Y_2$, since we can glue the path between $B_1$ and $B_2$, the circuit around $B_2$, and the path from $B_2$ to any of the $r$-squares $B \in B^r(Y_2)$ to get a path from $B_1$ to $B$. Therefore, for each equivalence class of $(r,\rho)$ sub-routers there exists a maximal subset $Y\subseteq X$ for which all elements of the equivalence class are sub-routers. Such a maximal subset $Y$ will sometimes be called a {\bf cluster of pivotals}, and a corresponding equivalence class is said to be {\bf spanned by $Y$}. For instance, in Figure~\ref{f.NeedRouters}, the left configuration has two clusters, spanned by the same three pivotals, while the right configuration has three clusters, each with a maximal $Y$ of two elements. 

In each equivalence class of sub-routers, we want to single out one of them. In order to do this in a way that is typically continuous w.r.t.~$\omega^\lambda$, we need to restrict ourselves to the case when every open cluster of $\omega^\lambda$ in $\Tor_M$ has diameter less than $M/3$; this will turn out to be typically the case when $\lambda$ is very negative. (For continuum percolation configurations, the diameter is the $\limsup_{r\to 0}$ of distances between $r$-boxes that are connected in the usual sense that there is a crossed quad with its opposite sides contained in the $r$-boxes. It is clear from Lemma~\ref{l.crossconv} that, in any coupling with $\omega_\eta^\lambda \xrightarrow{a.s.}\omega_\infty^\lambda$, the event that this diameter is at most $M/3$ converges almost surely.) Then, for $r<M/10$, the set of $(r,\rho)$ sub-routers for any $Y$ has an isometric embedding into $[-M/2,M/2)^2\subset \C$. The leftmost sub-router of the lowermost ones in such an embedding will be the same in any of these embeddings; moreover, its location in $\Tor_M$ can change only a little if we move each sub-router a little, to at most distance $M/10$. The set of these ``leftmost of lowermost'' sub-routers will be the {\bf $(r,\rho)$-mesoscopic routers of $X$}, or, after fixing $\rho=r_M/2$ (from Proposition~\ref{p.netconv}), the set of {\bf $r$-mesoscopic routers}.  
Note that by restricting ourselves to subsets $|Y|\geq 2$, clusters containing only one pivotal from $X$ will not have routers. 

Although we will not really need them, for the sake of symmetry in our presentation, analogously to the above routers that used primal (open) connections, we also define {\bf dual clusters of pivotals} and {\bf dual $r$-mesoscopic routers}.

We can now define the enhanced networks we promised.

\begin{figure}[htbp]
\centerline{
\AffixLabels{
\includegraphics[width=0.9 \textwidth]{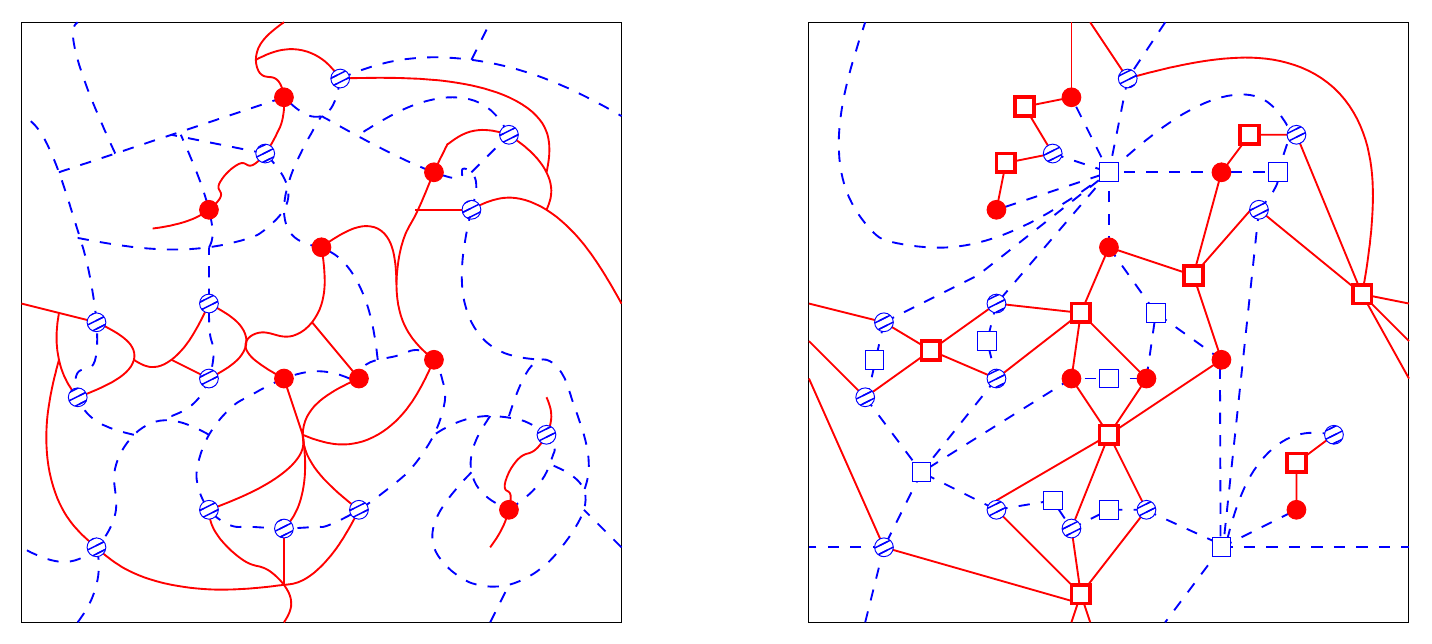}
}}
\caption{A schematic picture of a percolation configuration $\omega^\lambda_\eta$ with the pivotals of $\PPP^\eps_\vl$ on a torus, and the corresponding enhanced network. Pivotals are represented by circles, routers are represented by squares. Primal connections are shown using red solid lines, dual connections are shown using blue dashed lines.}
\label{f.EnNet}
\end{figure}

\bde\label{d.ErNetwork}
Assuming that the diameter of every cluster in $\omega^\lambda$ is at most $M/3$, the {\bf $r$-mesoscopic enhanced $\eps$-network} $\EnNet_M^{r\meso}(\omega^\lambda,\PPP^\eps_\vl)$ is the following vertex- and edge-labeled bipartite graph. One part of the vertex set is the set $X$ of the pivotals of $\PPP^\eps_\vl$, the other part is the $r$-mesoscopic routers of $X$ (both the primal and dual ones). The vertices in $\PPP^\eps_\vl$ are colored open or closed, according to the definitions before Proposition~\ref{p.PivSetCol}; the routers are colored in the obvious way. The edge set consists of the connections between the routers and the elements of their maximal $Y\subset X$, labelled primal or dual according to the color of the router. The edges are drawn on the torus so that they are homotopic (with fixed endpoints) to the connections they represent; clearly, one can also achieve that they do not intersect each other. See Figure~\ref{f.EnNet}.

If the assumption on the diameters is not satisfied, we take  $\EnNet_M^{r\meso}(\omega^\lambda,\PPP^\eps_\vl)$ to be empty.
\ede

Note that the networks of Definition~\ref{d.rNetwork} are measurable functions of these enhanced networks in a very simple way: there exists a primal (or dual) router with edges to $x_i,x_j \in X$ in  $\EnNet_M^{r\meso}(\omega^\lambda,\PPP^\eps_\vl)$ if and only if there is a primal (dual, resp.) edge between $x_i$ and $x_j$ in $\Net_M^{r\meso}(\omega^\lambda,\PPP^\eps_\vl)$. Moreover, the same proof as for Proposition~\ref{p.netconv}, together with Theorem~\ref{t.PivMeas}, implies the following:

\bpr[$r$-stabilization and $\eta$-convergence of enhanced networks]\label{p.ENetconv}\ 
\bi
\item[{\bf (i)}]
There is a measurable scale $\tilde r_M=\tilde r_M(\omega^\lambda_\infty,\PPP^\eps_\vl)\in (0,r^*)$ such that for all $r\in (0,\tilde r_M)$ we get the same $r$-mesoscopic enhanced $\eps$-network $\EnNet_M^{r\meso}(\omega^\lambda_\infty,\PPP^\eps_\vl)$ in the sense that the networks are the same, plus the colors in $X$ and the collections of primal and dual clusters of pivotals are also the same. (The corresponding routers do not exactly stabilize, since for a smaller $r$ new $(r,r_M/2)$ sub-routers can appear; but they cannot disappear, and hence each router does converge to a point in $\Tor_M$ as $r\to 0$.) This stabilized network will be called the {\bf enhanced $\eps$-network} $\EnNet^{\vl,\eps}_\infty=\EnNet_M(\omega^\lambda_\infty,\PPP^\eps_\vl)$. For discrete percolation configurations, the definition of $\EnNet^{\vl,\eps}_\eta=\EnNet_M(\omega^\lambda_\eta,\PPP^\eps_\vl)$ is the obvious one.
\item[{\bf (ii)}]
In any coupling with $\omega_\eta^\lambda \xrightarrow{a.s.}\omega_\infty^\lambda$ in $\Tor_M$, there is a coupling of $\PPP^\eps_\vl(\omega_\eta^\lambda)$ and $\PPP^\eps_\vl(\omega_\infty^\lambda)$ such that with probability tending to 1 as $\eta\to 0$, we have that $\EnNet^{\vl,\eps}_\eta=\EnNet_M(\omega_\eta^\lambda,\PPP^\eps_\vl)$ is the same as $\EnNet^{\vl,\eps}_\infty=\EnNet_M(\omega_\infty^\lambda,\PPP^\eps_\vl)$ in the sense that the vertex sets for $\eta$ and $\infty$ (consisting of the pivotals in $\PPP^\eps_\vl$ and the routers) are arbitrarily close to each other, and the  labelled graph structures coincide.
\ei
\epr

\begin{remark}
These enhanced networks are very useful planar (more precisely, toroidal) representations of the discrete and continuous percolation configurations, which was not a priori obvious how to achieve, since the quad-crossing space allows for non-planar configurations and hence is not ideal to express planarity.
\end{remark}

Using the enhanced networks, we are now going to define a spanning forest $\MSF^{\vl,\eps}$ with vertices being the primal routers in $\EnNet^{\vl,\eps}$. We will show in Section~\ref{s.epsapprox} that, for $\lambda<0$ very negative, $\lambda'>0$ very large, and $\eps>0$ small, this forest has a unique giant tree component, which will be the cut-off tree $\MST^{\vl,\eps}$ that approximates well the macroscopic structure of $\MST$ in $\Tor_M$. This forest $\MSF^{\vl,\eps}$ may have edges that intersect each other besides their endpoints; nevertheless, it consists of trees immersed into the torus $\Tor_M$ in the sense of Subsection~\ref{ss.ESF}, and the fact it turns out to approximate $\MST_\eta$ implies in particular that these possible intersections vanish in the limit.
 

\begin{figure}[htbp]
\centerline{
\AffixLabels{
\includegraphics[width=0.9 \textwidth]{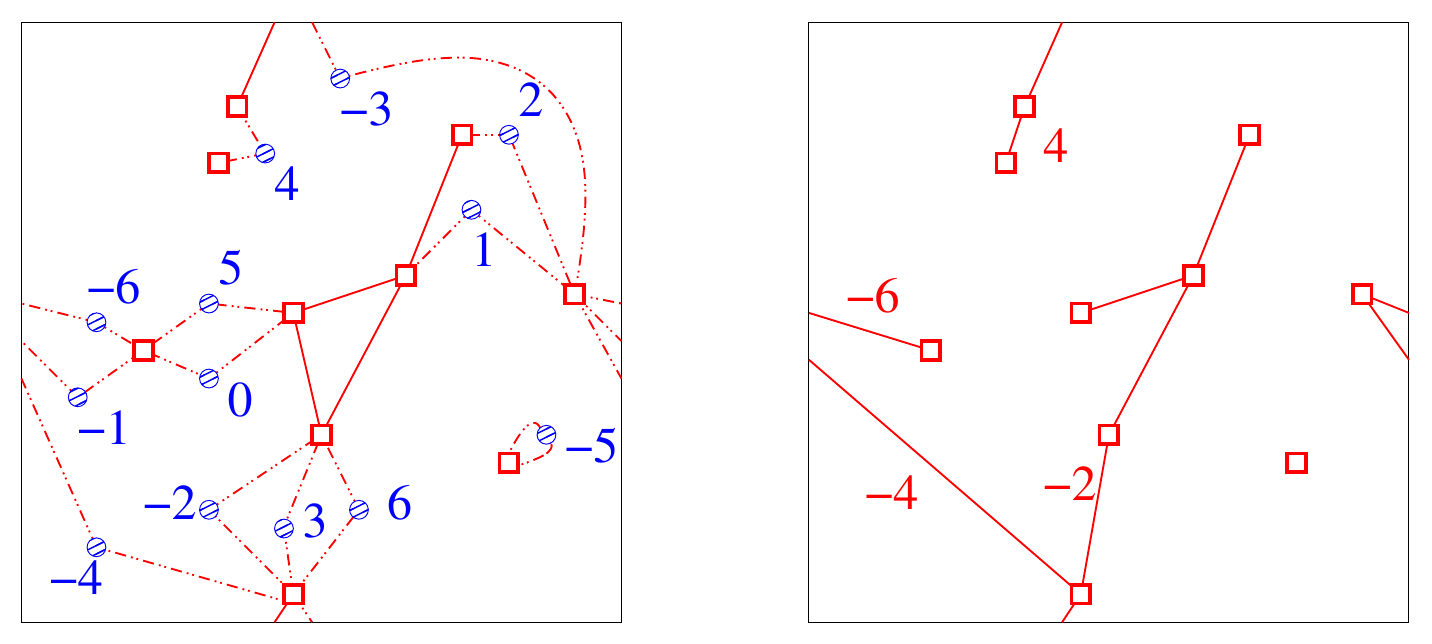}
}}
\caption{Building the cut-off forest $\MSF^{\vl,\eps}$ from the enhanced network of Figure~\ref{f.EnNet}. On level $\lambda$ there was a cycle that  had to be broken. The numbers from $-6$ to $6$ on the closed pivotals of $\omega^\lambda_\eta$ represent their levels $t_i\in(\lambda,\lambda')$ in $\PPP^\eps_\vl$ at which they become open. The resulting spanning forest has two components.}
\label{f.EnNetMST}
\end{figure}

\bde[Constructing the cut-off spanning forest  $\MSF^{\vl,\eps}$ on $\Tor_M$]\label{d.MSF}\
\bnum
\item The vertices are the primal routers in $\EnNet^{\vl,\eps}$. Connect two routers by an edge if they are both connected in $\EnNet^{\vl,\eps}$ to the same  open pivotal of $\omega^\lambda$. The resulting graph usually has several components (e.g., six of them on the left-hand picture of Figure~\ref{f.EnNetMST}), which more-or-less represent the $\lambda$-clusters in $\omega^\lambda$ (this will be made more precise in the next section). 
\item In each component of this graph, choose a spanning tree in an arbitrary deterministic way, and label each edge of this tree by $\lambda$.
\item For each pivotal $x_i$ of $\PPP^\eps_\vl$ that is closed in $\omega^\lambda$, add an edge between the corresponding routers, and label it by its $t_i\in (\lambda,\lambda')$ value. Note that these edges may be loops, as the one labelled by $-5$ on the left-hand picture of Figure~\ref{f.EnNetMST}, for instance.
\item As in the so-called reversed Kruskal algorithm, from each cycle delete the edge with the largest label, and get a minimal spanning tree in each component of the above graph.
\item Draw all the edges of the thus constructed forest as straight line segments, respecting the torus topology (i.e., choosing the line segment on the torus that is homotopic to the concatenation of the embedded edges of $\EnNet_M^{r\meso}(\omega^\lambda,\PPP^\eps_\vl)$ that gave rise to this edge of the forest). See the right-hand picture of Figure~\ref{f.EnNetMST}. Note that the edges may intersect each other besides their endpoints (even if this does not happen on this particular picture).
\enum
\ede

\section{Approximation of $\MST_\eta$ by the cut-off trees $\MST^{\vl,\eps}_\eta$}\label{s.epsapprox}

\subsection{Preparatory lemmas and the definition of $\MST^{\vl,\eps}_\eta$}\label{ss.prep}

Our first lemma is a RSW-type result that is interesting even in the critical case. Nevertheless, the simplest proof we have found uses our dynamical and near-critical stability results from \cite[Section 8]{DPSL}.

\begin{lemma}[Local Ring Lemma]\label{l.AnnBarrier}
There exists $\delta>0$ such that for any $\lambda\leq 0$ and any radius $R < |\lambda|^{-4/3}$, for all small enough mesh $0<\eta<\eta_0(\lambda,R)$, one has 
\[
\P_\eta[\A_{R,\lambda,\delta}]>1-\frac 1 {100}\,,
\]
where $\A_{R,\lambda,\delta}$ stands for the event that there exist $\lambda$-clusters for the restriction of $\omega_\eta^\lambda$ to the annulus $\A_{R, 2R}$:
$\calC_1,\ldots,\calC_N, \calC_{N+1}=\calC_1$ which satisfy the following conditions:
\bnum
\item for each $i\in [1,N]$, $\diam(\calC_i)\geq \delta R$; in particular, the clusters of the percolation configuration non-restricted to $A_{R,2R}$ that contain these $\calC_i$'s also have diameter $\geq \delta R$; 
\item for each $i\in [1,N]$, there exists at least one closed site $y_i$ neighboring both $\calC_i$ and $\calC_{i+1}$;
\item the circuit $\{ \calC_1,\ldots, \calC_N\}$ disconnects the annulus $A_{R,2R}$ in the sense that the two boundaries of the annulus are not connected in the graph $A_{R,2R}\setminus \bigcup_{i=1}^N (\calC_i \cup \{y_i\})$.
\enum
Moreover, we can choose the clusters $\calC_i$ and the points $y_i$ such that all the $y_i$'s are elements of the Poisson point set $\PPP^\eps_\vl$, with $\eps=\delta R$ and $\lambda'$ large enough (depending on $R$).
\end{lemma}

\proof
Consider the near-critical coupling $(\omega^t_\eta)_{t \in \R}$. For $\lambda'\ge \lambda$ large enough (on the order of $R^{-3/4}$), there is a probability at least $995/1000$ that $\omega_\eta^{\lambda'}$ has an open circuit even in the smaller annulus $A_{5R/4, 7R/4}$; this follows from known results on the correlation length, e.g., \cite[Theorem 10.7]{DPSL}. Now sample $\omega^\lambda_\eta$,  
consider some small $\eps>0$ to be fixed in a second, and let $\tilde\omega_\eta^{\lambda'}$ be the configuration where we open only those vertices in the coupling while getting from $\lambda$ to $\lambda'$ that are given in $\PPP^\eps_\vl$.
The assumption $R<|\lambda|^{-4/3}$, below the correlation length given by  \cite[Theorem 10.7]{DPSL}, implies that we can choose $\eta>0$ small enough compared to $R$ so that $\omega^\lambda_\eta$ has 4-arm probabilities inside the domain $A_{R,2R}$ that are comparable to the critical ones. Therefore, the critical case computations of \cite[Section 8]{DPSL} apply uniformly in $\lambda\leq 0$ and $R<|\lambda|^{-4/3}$, and by a straightforward modification of \cite[Proposition 8.6]{DPSL} from quad-crossings to annulus circuits, for $\eps=\delta R>0$ with $\delta>0$ small enough (uniformly in $\lambda$ and $R$), the probability that $\omega_\eta^{\lambda'}$ has an open circuit in $A_{5R/4, 7R/4}$ but $\tilde\omega^{\lambda'}_\eta$ does not have one in $A_{R,2R}$ is less than $5/1000$. Altogether, the probability that $\tilde\omega^{\lambda'}_\eta$ has an open circuit in $A_{R,2R}$ is at least $99/100$. But such a circuit must be composed of $\lambda$-clusters and $\eps$-important points that have become open, which implies that all these $\lambda$-clusters must have diameter at least $\eps$, and the lemma is proved.
\qed

\begin{lemma}[Global Ring Lemma]\label{l.UniformBarrier}
Fix $\delta>0$ as in Lemma~\ref{l.AnnBarrier}. For any $\lambda\leq 0$ and $\alpha>0$, there is a radius $r = r(\lambda,\alpha) < \tfrac{\delta}{2} |\lambda|^{-4/3}$ such that, for any small enough $\eta$, with probability at least $1-\alpha$, one can find around {\em all} points $x\in \Tor_M$ an annulus $A_{R,2R}$ surrounding $x$ with $\bar r= r/ \delta \leq R \leq |\lambda|^{-4/3}$ that satisfies the event $\A_{R,\lambda, \delta}$. (The choice $\bar r= r/\delta $ is made so that the clusters we find are at least of diameter $r$.)
\end{lemma}

\proof  Consider the covering of $\Tor_M$ by the squares given by $\bar r \Z^2$, and around each such $\bar r$-square, consider the dyadic annuli up to scale  $|\lambda|^{-4/3}$. By Lemma~\ref{l.AnnBarrier}, the probability  that there is an $\bar r$-square for which all the dyadic annuli fail to have the required ring of clusters is at most 
$$
O(1) (M/\bar r)^2 (1/100)^{\log_2 \frac {|\lambda|^{-4/3}}{\bar r}} 
 = O(M^2) \, |\lambda|^{4/3 \log_2 100} \bar r^{-2+\log_2 100}\,,
$$
which can be made arbitrarily small as $\bar r \to 0$.
\qed

Part (ii) of the next lemma again has a RSW feeling to it, and is again proved using \cite[Section 8]{DPSL}.

\begin{lemma}[Subcritical lakes joining the supercritical ocean]\label{l.ocean} Consider percolation $\omega^\lambda_\eta$ on $\Tor_M$ with $\lambda<-1$,  and fix an arbitrarily small $\alpha>0$.
\bi
\item[{\bf (i)}] For any $s>0$, if $\lambda<-1$ is small enough, then for all $\eta>0$ small enough, with probability at least $1-\alpha$, all clusters in $\omega^\lambda_\eta$ have diameter less than $s$. 
\item[{\bf (ii)}] Using part (i), take $\lambda<-1$ small enough so that with probability more than $1-\alpha/3$, all clusters in $\omega^\lambda_\eta$ have diameter less than $M/3$. Then, for any $r>0$ there is a $\lambda'_0(\lambda,r,\alpha)>0$ and an $\eps_0(\lambda,r,\alpha)>0$ such that for all $\lambda'\geq \lambda'_0$ and $\eps\leq \eps_0$, for all small enough $\eta>0$, with probability at least $1-\alpha$, all the clusters in $\omega^\lambda_\eta$ of diameter at least $r$ are connected via primal paths in the enhanced network $\EnNet^{\vl,\eps}_\eta= \EnNet_M(\omega_\eta^\lambda,\PPP^\eps_\vl)$ with $\vl=(\lambda,\lambda')$, defined in Proposition~\ref{p.ENetconv}, in the sense that each such cluster contains a primal router and these routers are all connected by primal edges (through closed or open pivotals, as in Definition~\ref{d.MSF}) in the enhanced network.
\ei
\end{lemma}

\proof It is proved in \cite[Theorem 10.7]{DPSL} that, for any fixed $s>0$, as $\lambda'\to\infty$, the probability of having an open circuit in a given annulus $A_{s/3,s}$ in $\omega^{\lambda'}_\eta$ converges to 1. Consider a tiling of $\Tor_M$ by $s/3$-squares, and the annuli of side-length $s$ centered around them. By a union bound, the probability of having open circuits in all of them converges to 1. When all these circuits are present, their union is a single component, and any subset of $\Tor_M$ of diameter at least $s$ intersects this cluster.

Running the above argument for dual circuits in $\omega^{\lambda}_\eta$ with $\lambda\to-\infty$ gives that, with probability tending to 1, the diameter of the largest open cluster must be less than $s$, which proves item~(i). 

For item~(ii), we will use the first paragraph with $s=r/10$. Note that any cluster $\calC$ of $\omega^\lambda_\eta$ with diameter at least $r$ will radially cross two such $(s/3,s)$-annuli at distance at least $r/2$ from each other, $A_1$ and $A_2$. Moreover, assuming that $\calC$ has diameter at most $M/3$ (which is satisfied for all clusters $\calC$ with probability more than $1-\alpha/3$), we can choose $A_1$ and $A_2$ with the additional property that, for each of them, not all of the eight neighboring inner squares are intersected by $\calC$. On the other hand, the first paragraph says that if $\lambda'_0$ is large enough, then with probability more than $1-\alpha/3$, all the $(s/3,s)$-annuli will have open circuits in $\omega^{\lambda'_0}_\eta$. Now we use \cite[Proposition 8.1]{DPSL}, which implies that in the configuration $\tilde\omega^{\lambda'_0}_\eta$ that we get by starting from the configuration $\omega^\lambda_\eta$ and opening only the pivotal points of $\PPP^\eps_{\vl_0}$, with $\vl_0=(\lambda,\lambda'_0)$, all these open circuits in the $(s/3,s)$-annuli will already be there with probability more than $1-2\alpha/3$, provided that $0<\eps \leq \eps_0(\lambda'_0,r,\alpha)$ is small enough. This can happen in our two above annuli $A_1,A_2$ only if $\calC$ neighbours points of $\PPP^\eps_{\vl_0}$ in both annuli that are closed in $\omega^\lambda_\eta$ but open in $\tilde\omega^{\lambda'_0}_\eta$. These $\eps$-pivotal points appear in the enhanced network  $\EnNet^{\vl_0,\eps}_\eta=\EnNet_M(\omega_\eta^\lambda,\PPP^\eps_{\vl_0})$, and $\calC$ has a primal router connecting these pivotals with high probability (as $\eta\to 0$). 

Furthermore, the open circuits of $\tilde\omega^{\lambda'_0}_\eta$ in the $(s/3,s)$-annuli must be composed of clusters of $\omega^\lambda_\eta$ with diameter at least $\min\{s/3,\eps\}$, joined by $\eps$-pivotals. The expected number of disjoint such clusters is bounded from above by some $K=K(M,\lambda,r,\eps)<\infty$, uniformly in $\eta$ (as follows from \cite{AizIIC}), hence in all of them we have primal routers with probability tending to 1, as $\eta\to 0$. That is, altogether, with probability at least $1-\alpha$, all the $\lambda$-clusters of diameter at least $r$ have diameter at most $M/3$, and they are all connected in the enhanced network $\EnNet^{\vl_0,\eps_0}_\eta$, for $\lambda'_0$ chosen large enough, then $\eps_0(\lambda'_0,r,\alpha)$ and then $\eta$ chosen small enough. 

Finally, it is clear by construction that, if $\lambda' \geq \lambda'_0$ and $\eps \leq \eps_0$, then $\EnNet^{\vl_0,\eps_0}_\eta$ is a graph minor of $\EnNet^{\vl,\eps}_\eta$, if $\eta$ is small enough. Therefore, the conclusion of the previous paragraph holds also for $\EnNet^{\vl,\eps}_\eta$, and we are done.
\qed

Using the above lemmas, we can now see why there is typically a unique giant component in the cut-off forests $\MSF^{\vl,\eps}_\eta$ and $\MSF^{\vl,\eps}_\infty$ of Definition~\ref{d.MSF}:

\bl[Defining the cut-off trees $\MST^{\vl,\eps}_\eta$ and $\MST^{\vl,\eps}_\infty$]\label{l.giant}
For any $M>0$, any small $s>0$ and $\alpha>0$, if $\lambda < -1$ is very negative, $\eps>0$ is small, and $\lambda' > 1$ is large enough, then with probability at least $1-\alpha$ for any small enough mesh $\eta>0$, there is a {\bf unique giant component} in the cut-off forest $\MSF^{\vl,\eps}_\eta$ (and hence, by Proposition~\ref{p.ENetconv}, in $\MSF^{\vl,\eps}_\infty$), with the properties that it comes to distance at most $s$ from any point of $\Tor_M$, while all other components of $\MSF^{\vl,\eps}_\eta$ have diameter at most $s$. This giant tree component will be our {\bf approximating cut-off tree}, denoted by $\MST^{\vl,\eps}_\eta$ and $\MST^{\vl,\eps}_\infty$; whenever the above large probability event fails to occur, we set $\MST^{\vl,\eps}$ to be a single point in $\Tor_M$, and call this tree degenerate.
\el

\proof
Take $\lambda<-1$ such that $2 |\lambda|^{-4/3}<s$ holds and the diameter bound of Lemma~\ref{l.ocean}~(i) applies. By Lemma~\ref{l.UniformBarrier}, with probability at least $1-\alpha/2$, every point of $\Tor_M$ has in its $s$-neighborhood a ring of $\lambda$-clusters of diameter at least $r(\lambda,\alpha)$ each (possibly a single cluster, but still of diameter in $\Tor_M$ less than $s$). Now, if we take $\eps>0$ small and $\lambda'>0$ large, then Lemma~\ref{l.ocean}~(ii) says that with probability at least $1-\alpha/2$ all $\lambda$-clusters of diameter at least $r(\lambda,\alpha)$ get connected in the enhanced network $\EnNet^{\vl,\eps}_\eta$. Therefore, with probability altogether at least $1-\alpha$, there is a component of the graph of Definition~\ref{d.MSF} that has distance at most $s$ from any point of $\Tor_M$, while all other components have diameter at most $r(\lambda,\alpha)\leq s$. The spanning trees of these components inherit these properties, hence we are done.
\qed

So, we have finally defined the trees $\MST^{\vl,\eps}_\eta$ and $\MST^{\vl,\eps}_\infty$. Naturally, they are also immersed spanning trees in the sense of Definition~\ref{d.ESF}. Now, using Proposition~\ref{p.ENetconv}, we can easily show that they are close to each other in the space of essential spanning forests, introduced in that definition:

\bc\label{cor.MSTABNW}
For any $M>0$ and $s,\alpha>0$, if $\lambda<-1$ is very negative, $\eps>0$ is small, and $\lambda'>1$ is large enough, then, in the coupling of Proposition~\ref{p.ENetconv}~(ii) between $(\omega^\lambda_\eta,\PPP^\eps_\vl)$ and $(\omega^\lambda_\infty,\PPP^\eps_\vl)$, we have
$$
\Pb{ d_{\Omega_M}(\MST^{\vl,\eps}_\eta,\MST^{\vl,\eps}_\infty) < s} > 1-\alpha\,,
$$
for all $\eta>0$ small enough.
\ec

\proof The parameters $\lambda,\lambda',\eps$ can be set so, by Lemma~\ref{l.giant}, that both $\MST^{\vl,\eps}_\eta$ and $\MST^{\vl,\eps}_\infty$ are non-degenerate with probability at least $1-\alpha/2$, for any $\eta>0$ small enough. Now, by Proposition~\ref{p.ENetconv}~(ii), we can take $\eta>0$ so small that, with probability at least $1-\alpha/2$, the enhanced networks $\EnNet^{\vl,\eps}_\eta$ and $\EnNet^{\vl,\eps}_\infty$ agree as graphs and the Hausdorff distance between their vertex sets is less than $s$. On the event that both trees exist, the networks agree, and the vertex sets are closer than $s$ to each other, which occurs with probability at least $1-\alpha$, the uniform distance between the corresponding $\ell$-trees is always less than $s$, and hence the sum with the weights $2^{-\ell}$ is also less than $s$, and we are done. 
\qed

\subsection{Approximation as $\eps\to 0$ and $(\lambda,\lambda')\to (-\infty,\infty)$.}\label{ss.eps}

After these preparations, we can turn to approximating $\MST_\eta$ on the lattice $\eta\Tg \cap \Tor_M$ by the cut-off trees $\MST^{\vl,\eps}_\eta$. First of all, note that if two vertices of $\eta\Tg \cap \Tor_M$ are in the same $\lambda$-cluster, then the path in $\MST_\eta$ that connects them remains in that $\lambda$-cluster. This also means that for any two $\lambda$-clusters of $\eta\Tg\cap \Tor_M$, there is a unique path in $\MST_\eta$ that connects them.

The next lemma is a key step in approximating $\MST_\eta$ by cut-off trees:

\begin{lemma}[Paths through macroscopic clusters]\label{l.macro}
Fix any $\alpha>0$. If $\lambda<-1$ is small enough, $\rho>0$ is arbitrary, $\eps > 0$ is small and $\lambda'>0$ is large enough, then with probability at least $1-\alpha$, for any two clusters of diameter at least $\rho$ in $\omega^\lambda_\eta$, there is a unique path in $\MSF^{\vl,\eps}_\eta$ that connects the two clusters, and the unique path in $\MST_\eta$ doing the same goes through the same closed pivotals of $\PPP^\eps_\vl$, and hence the distance of these two paths in the uniform metric is at most the maximal diameter of all $\lambda$-clusters.
\end{lemma}

\begin{figure}[htbp]
\SetLabels
(0.2*0.86)$\rho$\\
(0.33*0.72)$\eps_1$\\
(0.21*0.39)\blue{$x$}\\
(0.65*0.15)\textcolor{Maroon}{$\calC$}\\
\endSetLabels
\centerline{
\AffixLabels{
\includegraphics[width=0.5 \textwidth]{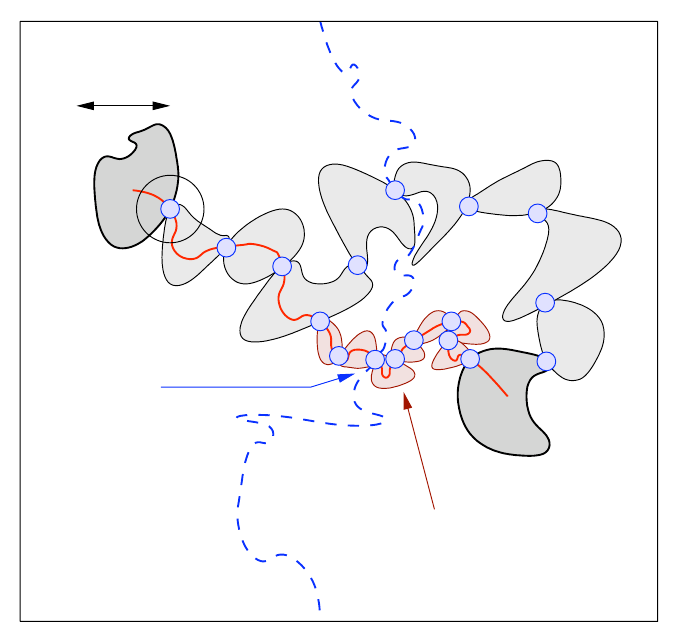}
}}
\caption{The path in $\MST_\eta$ connecting two large $\lambda$-clusters does not go through very small $\lambda$-clusters, basically because of the near-critical stability of 4-arm probabilities.}
\label{f.largeclusters}
\end{figure}
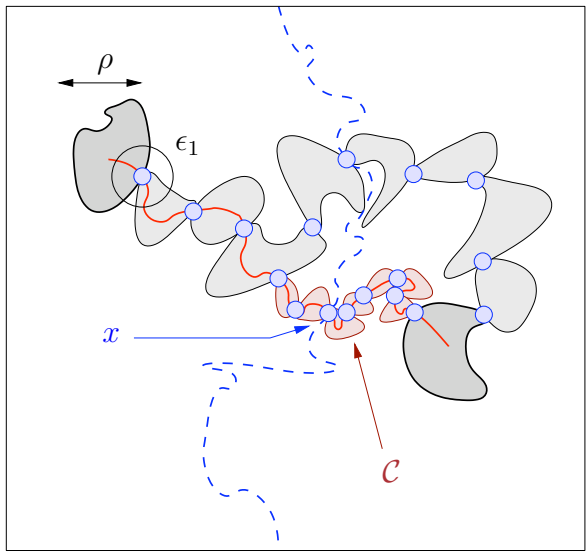

\proof
Choose $\lambda<-1$ then $\eps_1>0$ and $\lambda'_1>1$ so that the event of Lemma~\ref{l.ocean}~(ii), with $r:=\rho$ there, occurs with probability at least $1-\alpha/3$. Condition on this event, set $\vl_1=(\lambda,\lambda'_1)$, and consider the path in $\MSF^{\vl_1,\eps_1}_\eta$ that connects two of the clusters. There is a corresponding path in $\eta\Tg$, going through the same finitely many $\eps_1$-important points of $\omega^\lambda_\eta$ and some $\lambda$-clusters, using labels at most  $\lambda'_1$. Therefore, the true path in $\MST_\eta$ also uses labels at most $\lambda'_1$. Assume now that this path goes through some $\lambda$-cluster $\calC$ of diameter at most $r\ll\rho$. 
This path must go through a vertex $x$ of $\eta\Tg$, neighboring $\calC$, with the following properties (see Figure~\ref{f.largeclusters}):
\bi
\item it is closed in $\omega^\lambda_\eta$ but open in $\omega^{\lambda'_1}_\eta$;
\item it has two $\lambda$-closed arms emanating from it, which together separate the two clusters of diameter at least $\rho$ that we started with;
\item on the side of these two closed arms that contains $\calC$, there is a $\lambda$-open arm from $x$ only to distance at most $r$. (Note here that $x$ might be neighbouring $\lambda$-clusters of diameter more than $r$, but we can choose the closed arms so that they separate those clusters from $\calC$.)
\ei

If $x$ had the alternating 4-arm event to a distance more than $r$ in $\omega^\lambda_\eta$, that could happen only if the two open arms out of these four were on the side of the two long closed arms that does not contain $\calC$. This would altogether yield a 5-arm event to distance $r$, of type closed-open-closed-open-closed, where the extra closed arm in the middle is needed to get an {\it alternating} 4-arm event. Moreover, since the labels along the path in $\MST_\eta$ are all at most  $\lambda'_1$, we would get a $(\lambda,\lambda_1')$-near-critical six-arm event from $x$ to distance $r$, as defined in Proposition~\ref{p.stab} (the sixth arm is of length at least $\rho/2$, because of being $\lambda_1'$-connected to the $\lambda$-cluster of diameter at least $\rho$). By that proposition and by the 6-arm exponent being larger than two (see \cite[Corollary A.8]{SchSt}) this happens with very small probability if $\eta$ is small enough. So, we can basically assume that $x$ is not $r$-pivotal in $\omega^\lambda_\eta$. On the other hand, if we now change all the labels above $\lambda$ along the path in $\MST_\eta$ to $\lambda$, then, in the new configuration, $x$ will have the alternating 4-arm event to distance at least $\eps_1$. Since the labels we have changed are all in $[\lambda,\lambda'_1]$, we can apply another version of near-critical stability, Lemma 8.5 of \cite{DPSL}, saying that the probability that the importance of a given vertex $x$ can be changed from $r$ to $\eps_1$ by such label changes is at most  
$$
O_{\lambda,\lambda_1'}(1)\,r^2\,\alpha_4(\eta,r)\,\alpha_4(\eps_1,1)^{-1}\,.
$$
Additionally, the probability that the status of this vertex $x$ is different in $\omega^{\lambda'_1}_\eta$ than in $\omega^\lambda_\eta$ is $O_{\lambda,\lambda_1'}(1)\, \eta^2/\alpha_4(\eta,1)$, independently of everything. The number of possible vertices $x\in\Tor_M$ is $O_M(\eta^{-2})$. Altogether, the expected number of vertices $x$ satisfying this scenario is 
$$
O_{M,\lambda,\lambda_1'}(1)\,r^2\,\alpha_4(r,1)^{-1}\,\alpha_4(\eps_1,1)^{-1}\,.
$$
This is arbitrarily small if $r$ is small, hence the probability that this scenario happens is also small. Summarizing, there exists $r>0$ depending on $M,\alpha,\lambda$ and $\rho$, such that for all small enough $\eta>0$, with probability at least $1-\alpha/3$, the path in $\MST_\eta$ connecting any two $\lambda$-clusters of diameter at least $\rho$ does not go through $\lambda$-clusters smaller than $r$. 

Now choosing $\eps>0$ small and $\lambda'>1$ large, again by Lemma~\ref{l.ocean}~(ii), the enhanced network $\EnNet^{\vl,\eps}_\eta$ will connect all the $\lambda$-clusters of diameter at least $r$ with probability at least $1-\alpha/3$. Altogether, with probability at least $1-\alpha$, for any two $\lambda$-clusters of diameter at least $\rho$, the unique paths in $\MST_\eta$ and $\MSF^{\vl,\eps}_\eta$ both go through the same $\lambda$-clusters, connected by $\lambda$-closed pivotals of importance at least $r$. The last half sentence of the lemma follows immediately from the way Definition~\ref{d.MSF} is done.
\qed

We can now easily prove the main result of this section:

\begin{proposition}\label{p.epsapprox}
For any $M>0$ and $s,\alpha>0$, if $\lambda<-1$ is very negative, $\eps>0$ is small, and $\lambda'>1$ is large enough, then we have
$$
\Pb{ d_{\Omega_M}(\MST_\eta,\MST^{\vl,\eps}_\eta) < s} > 1-\alpha\,,
$$
for all $\eta>0$ small enough.
\end{proposition}

\proof 
As in the proof of Lemma~\ref{l.giant}, take $\lambda<-1$ such that with probability at least $1-\alpha/2$, all $\lambda$-clusters in $\Tor_M$ have diameter less than $s$, and every point of $\Tor_M$ has in its $s/2$-neighborhood a ring of $\lambda$-clusters of diameter at least $r$ each, for some $0<r<s/2$, uniformly in $\eta$, as provided by Lemma~\ref{l.UniformBarrier}. Now, if we take $\eps>0$ small and $\lambda'>0$ large, then, with probability at least $1-\alpha/2$, all $\lambda$-clusters of diameter at least $r$ are connected in $\MST^{\vl,\eps}_\eta$, and, by Lemma~\ref{l.macro}, the paths connecting them are at a uniform distance at most $s$ from the corresponding paths of $\MST_\eta$. We will assume that both events of probability at least $1-\alpha/2$ hold. 

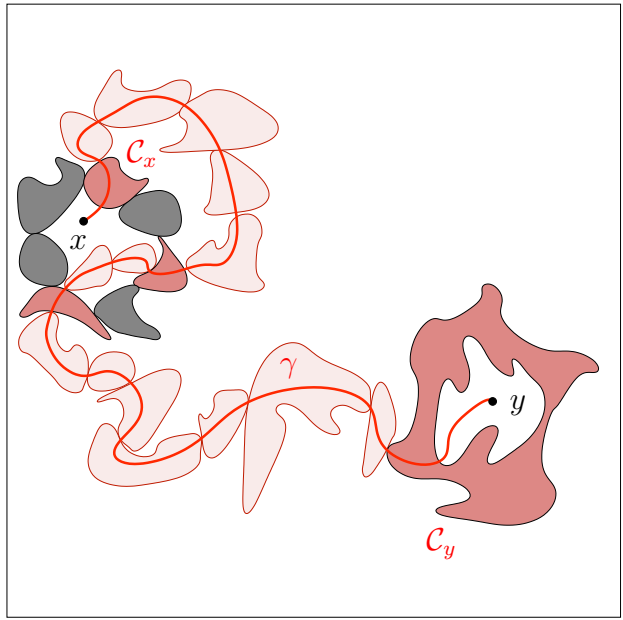
\begin{figure}[htbp]
\SetLabels
(0.13*0.6)$x$\\
(0.82*0.34)$y$\\
(0.46*0.4)\red{$\gamma$}\\
(0.23*0.74)\red{$\calC_x$}\\
(0.7*0.11)\red{$\calC_y$}\\
\endSetLabels
\centerline{
\AffixLabels{
\includegraphics[width=0.5 \textwidth]{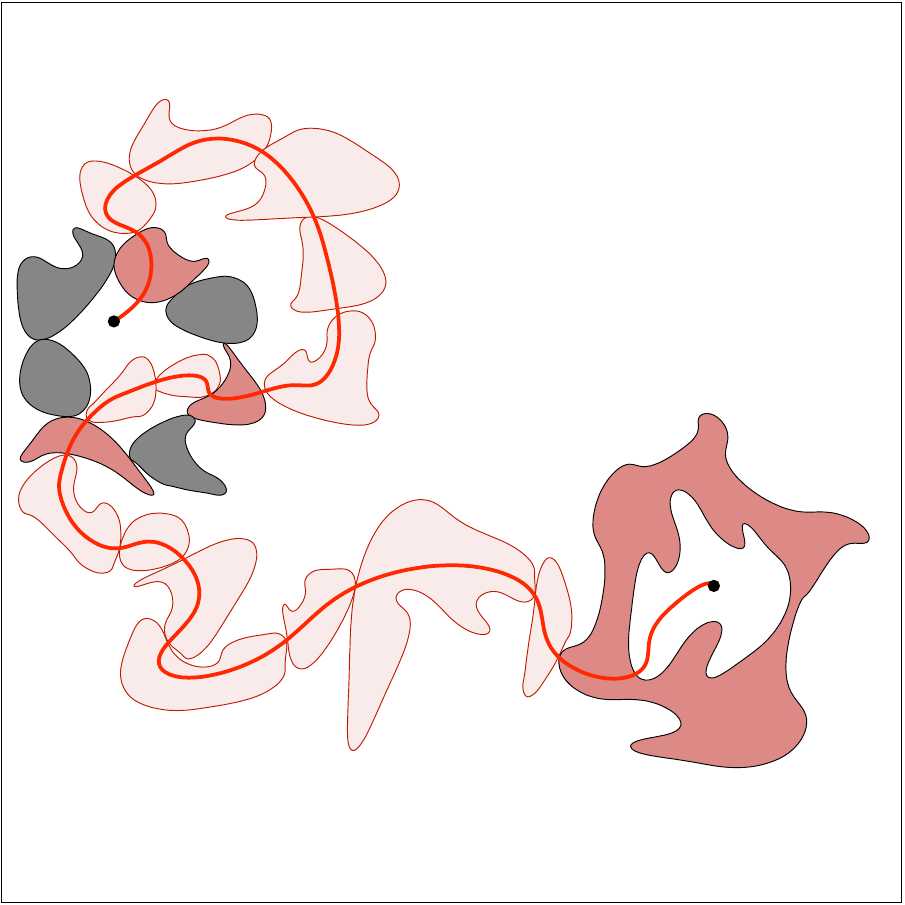}
}}
\caption{Paths in $\MST_\eta$ can be approximated by paths through macroscopic $\lambda$-clusters.}
\label{f.epsapprox}
\end{figure}

Consider any path $\gamma$ of $\MST_\eta$ connecting some $x,y\in\Tor_M$.  Both $x$ and $y$ have the above-mentioned ring of macroscopic $\lambda$-clusters around them, and $\gamma$ must intersect at least one member of each ring. See Figure~\ref{f.epsapprox}. But then, the part of $\gamma$ connecting the intersected members closest to $x$ and $y$, denoted by $\calC_x$ and $\calC_y$, respectively, by the previous paragraph, is uniformly $s$-close to a path in $\MST^{\vl,\eps}_\eta$, denoted by  $\gamma^{\vl,\eps}$. And this  $\gamma^{\vl,\eps}$ is of course $s$-close to the entire $\gamma$, since the parts of $\gamma$ going from $x$ to $\calC_x$ and from $y$ to $\calC_y$ are contained in the $s$-neighborhoods of $\calC_x$ and $\calC_y$.

In the other direction, consider any path $\gamma^{\vl,\eps}$ in $\MST^{\vl,\eps}_\eta$, connecting two routers. The clusters of pivotals corresponding to these routers have diameter at most $s$, but could be rather small. Nevertheless, fixing one point in each cluster, there is a ring of macroscopic $\lambda$-clusters around each, which certainly contains a cluster of pivotals that  $\gamma^{\vl,\eps}$ goes through. The rest of the proof is just as above.

Now that we have good approximations for paths in the two trees connecting any two vertices, the extension to trees with $\ell >2$ leaves is straightforward.
\qed

\section{Proof of the main result}\label{s.main}

\subsection{Putting the pieces together for $\MST$ on tori $\Tor_M$}\label{ss.pieces}
 
In this subsection, we prove convergence in any fixed torus $\Tor_M$.

\bth[Limit of $\MST_\eta$ and $\MST_\infty^{\vl,\eps}$ in $\Tor_M$]\label{t.torus}
In the metric space $\Omega_M$ of spanning trees in the torus $\Tor_M$, as defined in Definition~\ref{d.ESF}, the spanning tree $\MST_\eta$ on the lattice $\eta\Tg \cap \Tor_M$ converges in law to a translation invariant $\MST_\infty$, which is also the distributional limit of the cut-off trees $\MST_\infty^{\vl,\eps}$, as $\vl\to(-\infty,\infty)$ and $\eps\to 0$.
\eth

\proof 
Using the results of the previous section, the proof is classical; e.g., the exact same strategy was used in \cite[Section 9]{DPSL}. By Proposition~\ref{p.epsapprox}, for any $k\in\N$ there exists $\vl_k=(\lambda_k,\lambda'_k)$ and $\eps_k>0$, such that, for all $0<\eta<\eta_k$ sufficiently small,
\be\label{e.etaeta}
\Pb{ d_{\Omega_M}(\MST_\eta,\MST^{\vl_k,\eps_k}_\eta) < 2^{-k}} > 1-2^{-k}\,.
\ee
Now, by Corollary~\ref{cor.MSTABNW}, there is a coupling between  $(\omega^\lambda_\eta,\PPP^\eps_\vl)$ and $(\omega^\lambda_\infty,\PPP^\eps_\vl)$, and by the same token, between $\omega^{[\lambda,\lambda']}_\eta$ and $(\omega^\lambda_\infty,\PPP^\eps_\vl)$, such that, for all $0<\eta<\eta'_k$ sufficiently small,
\be\label{e.etainfty}
\Pb{ d_{\Omega_M}(\MST^{\vl_k,\eps_k}_\eta,\MST^{\vl_k,\eps_k}_\infty) < 2^{-k} } > 1-2^{-k}\,.
\ee
Combining~(\ref{e.etaeta}) and~(\ref{e.etainfty}) using the triangle inequality, in the same coupling, 
$$
\Pb{ d_{\Omega_M}(\MST_\eta,\MST^{\vl_k,\eps_k}_\infty) < 2^{-k+1} } > 1-2^{-k+1}\,.
$$
We can now couple all the trees $\MST^{\vl_k,\eps_k}_\infty$ to $\MST_\eta$ one-by-one, and given $\MST_\eta$, conditionally independently to each other, such that, again using the triangle inequality,
\be\label{e.kl}
\Pb{ d_{\Omega_M}(\MST^{\vl_k,\eps_k}_\infty,\MST^{\vl_\ell,\eps_\ell}_\infty) < 2^{-k+2}\text{ for all }\ell \geq k\text{ simultaneously} } > 1-2^{-k+2}\,.
\ee
Using the Borel-Cantelli lemma for the events appearing on the left hand side of (\ref{e.kl}) shows that the sequence $\MST^{\vl_k,\eps_k}_\infty$ is almost surely a Cauchy sequence in $\Omega_M$. The space is complete, hence there is an almost sure limit $\MST_\infty$. Of course, this limit may a priori depend on the sequences $\{\vl_k\}$, $\{\eps_k\}$ and on the coupling. However, using the triangle inequality again, going through $\MST^{\vl_k,\eps_k}_\infty$, we have that for any $\delta>0$, if $\eta>0$ small enough, then
$$
\Pb{ d_{\Omega_M}(\MST_\eta,\MST_\infty) < \delta } > 1-\delta\,.
$$
Therefore, in this coupling, $\MST_\eta$ converges in probability, and hence in law, to $\MST_\infty$, in the metric space $\Omega_M$. Since $\MST_\eta$ and the metric $d_{\Omega_M}$ are translation invariant, the limit $\MST_\infty$ is also invariant.

To prove the convergence of $\MST_\infty^{\vl,\eps}$, note that the bounds~(\ref{e.etaeta}) and~(\ref{e.etainfty}) hold not just for $\vl_k$ and $\eps_k$, but for all $\eps<\eps_k$ and  $\lambda<\lambda_k$ and $\lambda' > \lambda'_k$, thus we have that $\MST^{\vl,\eps}_\infty$ is close in distribution in the $d_{\Omega_M}$-metric to $\MST_\eta$ and hence to $\MST_\infty$.
\qed

\subsection{Extension to the full plane; invariance under translations, scalings and rotations}\label{ss.full}

We are now ready to prove the main result of this paper.

\proofof{Theorem~\ref{t.main}}
We will use the notation $\MST_\eta^M$ and $\MST_\infty^M$ for $\MST_\eta$ and its scaling limit on the torus $\Tor_M$. We will also use the approximations $\MST^{\vl,\eps,M}$.

It was proved in \cite[equation~(8.1)]{ABNW} that $\MST_\eta$ is {\bf uniformly quasi-local} in the sense that for any $\delta>0$ and compact $\Lambda\subset \C$ there exists a $\bar\Lambda_\delta\subset \C$ such that for any small enough $\eta>0$, with probability at least $1-\delta$, all trees with leaves in $\Lambda$ are contained in $\bar\Lambda_\delta$. Since this event is measurable w.r.t.~the percolation ensemble inside $\bar\Lambda_\delta$, by taking $\delta>0$ small and $M>0$ so large that  $\bar\Lambda_\delta \subset [-M,M]^2$, we get that  the law of $\MST_\eta$ restricted to $\Lambda$ is close in total variation distance to the law of $\MST_\eta^M$ restricted to $\Lambda$. By the $(\vl,\eps)$-approximation result Proposition~\ref{p.epsapprox}, the same holds for $\MST^{\vl,\eps,M}_\eta$, and by the uniformity in $\eta>0$, also for $\MST^{\vl,\eps,M}_\infty$. In the proof of Theorem~\ref{t.torus}, we have constructed $\MST^M_\infty$ as a limit of $\MST^{\vl,\eps,M}_\infty$, thus we also have that the law of $\MST^M_\infty$ restricted to $\Lambda$ converges as $M\to\infty$, in the metric $d_{\Omega_M}$ that is based on the flat Euclidean metric on $\Tor_M$.

Now we take $\Lambda=[-L,L]^2$, with $L\to\infty$. As pointed out at the beginning of Subsection~\ref{ss.ESF}, the metric defined in~(\ref{e.metric}) for $\hat\C$ is equivalent to the Euclidean metric in bounded domains, while the distance between any two points in $\hat\C\setminus [-L,L]^2$  is at most $O(1/L)$. Thus the uniform distance between any two trees embedded in $\hat\C\setminus[-L,L]^2$ is at most $O(1/L)$, and if two essential spanning forests are $\delta$-close in the metric $d_{\Omega_M}$ restricted to $[-L,L]^2$, then their distance in $d_\Omega$ is $O_L(\delta)+O(1/L)$. Therefore, the convergence in $d_{\Omega_M}$ for any given $\Lambda\subset\C$, established in the previous paragraph, implies convergence in $d_\Omega$.

Translation invariance of the limit measure $\MST_\infty$ follows from a standard trick: for any compact $\Lambda\subset \C$, quasi-locality implies that the limit of $\MST_\infty^{[-M,M]^2}$ restricted to $\Lambda$, as $M\to\infty$, is the same as the limit of $\MST^{[-M+x,M+x]^2}_\infty$ restricted to the same $\Lambda$, for any $x\in\R$, and hence $\MST_\infty$ restricted to $\Lambda$ has the same distribution as restricted to $\Lambda-x$.

To prove scale-invariance, consider the scaling $f_\alpha(z):=\alpha z$. The conformal covariance of the pivotal measures, proved in \cite[Theorem 6.1]{GPS2a}, says that
\be\label{e.fmu}
(f_\alpha)_*(\mu^\eps(\omega^\lambda_\infty)) ~=~ \alpha^{-3/4} \mu^{\alpha\eps}(f_\alpha(\omega^\lambda_\infty))\,.
\ee
Also, by the conformal covariance of $\omega^\lambda_\infty$, proved in \cite[Theorem 10.3]{DPSL}, we have 
\be\label{e.fomega}
f_\alpha(\omega^\lambda_\infty) ~\overset{d}{=}~ \omega^{\alpha^{-3/4}\lambda}_\infty\,.
\ee
Scaling the spatial intensity measure of a Poisson point process by $\alpha^{-3/4}$ as in~(\ref{e.fmu}) is the same as scaling the time duration by the same factor, in the sense that there is a natural coupling in which the spatial coordinates of the arrivals are the same, and there is a simple scaling between the time coordinates. Thus, combining~(\ref{e.fmu}) and~(\ref{e.fomega}), and denoting the notion of ``same'' in the previous sentence by $\approx$, we have  
\be\label{e.fPPP}
f_\alpha\big(\PPP^\eps_\vl(\omega^\lambda_\infty)\big) ~\overset{d}{\approx}~ \PPP^{\alpha\eps}_{\alpha^{-3/4}\vl}(\omega^{\alpha^{-3/4}\lambda}_\infty)\,.
\ee
Since our constructions of $\MSF^{\vl,\eps}_\infty$ and $\MST^{\vl,\eps}_\infty$ in Definition~\ref{d.MSF} and Lemma~\ref{l.giant} are equivariant under spatial and time scalings, the identities~(\ref{e.fomega}) and~(\ref{e.fPPP}) imply that
$$
f_\alpha(\MST^{\vl,\eps,M}_\infty) ~\overset{d}{=}~ \MST^{\alpha^{-3/4}\vl ,\, \alpha\eps ,\, \alpha M}_\infty\,.
$$
Since we obtained $\MST_\infty$ as a limit of $\MST^{\vl,\eps,M}_\infty$ with $\vl \to(-\infty,\infty)$, $\eps\to 0$, $M\to\infty$, the last identity gives that $f_\alpha(\MST_\infty) \overset{d}{=} \MST_\infty$.

Next, let $f_\theta:\C\lora\C$ be the rotation by angle $\theta$. Now \cite[Theorem 6.1]{GPS2a} and \cite[Theorem 10.3]{DPSL} give for full plane configurations that 
\be\label{e.rotPPP}
f_\theta\big( \omega^\lambda_\infty, \PPP^\eps_\vl(\omega^\lambda_\infty) \big) ~\overset{d}{=}~
\big(\omega^\lambda_\infty, \PPP^{\eps, \theta}_\vl(\omega^\lambda_\infty) \big)\,,
\ee
where $\PPP^{\eps, \theta}_\vl$ is constructed using a rotated grid to define $\eps$-importance. (As pointed out in \cite[Remark 6.3]{GPS2a}, this rotational equivariance of the $\eps$-importance measure and hence the Poisson point process is not a tautology, since the normalization factor in the definition of the measure is not changed with the rotation.) Now, if we want to consider $\MST^{\vl,\eps}$ on the torus $\Tor_M$, the rotated $\eps$- and $r$-grids cannot be exactly defined; nevertheless, we can consider the squares in the grid fully contained in $f_\theta([-M,M]^2)$, and make some arbitrary definition close to the boundary --- due to quasi-locality, this will not matter. Hence, from~(\ref{e.rotPPP}) we get that for large $M>0$, the distribution of $f_\theta(\MST^{\vl,\eps,M}_\infty)$ restricted to some fixed domain $\Lambda$, which is close to $f_\theta(\MST_\infty)$ restricted to $\Lambda$, is close to $\MST^{\vl,\eps,M,\theta}_\infty$ restricted to $\Lambda$. On the other hand,  Corollary~\ref{cor.MSTABNW} and Proposition~\ref{p.epsapprox} work fine with the rotated grids, giving that $\MST^{\vl,\eps,M,\theta}_\infty$ is close to $\MST_\eta^{\vl,\eps,M,\theta}$, and the latter is close to $\MST_\eta$. Finally, since $\MST_\eta$ is close to $\MST_\infty$, after taking all the limits we get that $f_\theta(\MST_\infty)$ agrees with $\MST_\infty$ in distribution.
\qed


\section{Geometry of the limit tree $\MST_\infty$}\label{s.MSTgeom}

\subsection{Degree types and pinching}\label{ss.degtype}

\begin{figure}[htbp]
\centerline{
\includegraphics[width=0.9 \textwidth]{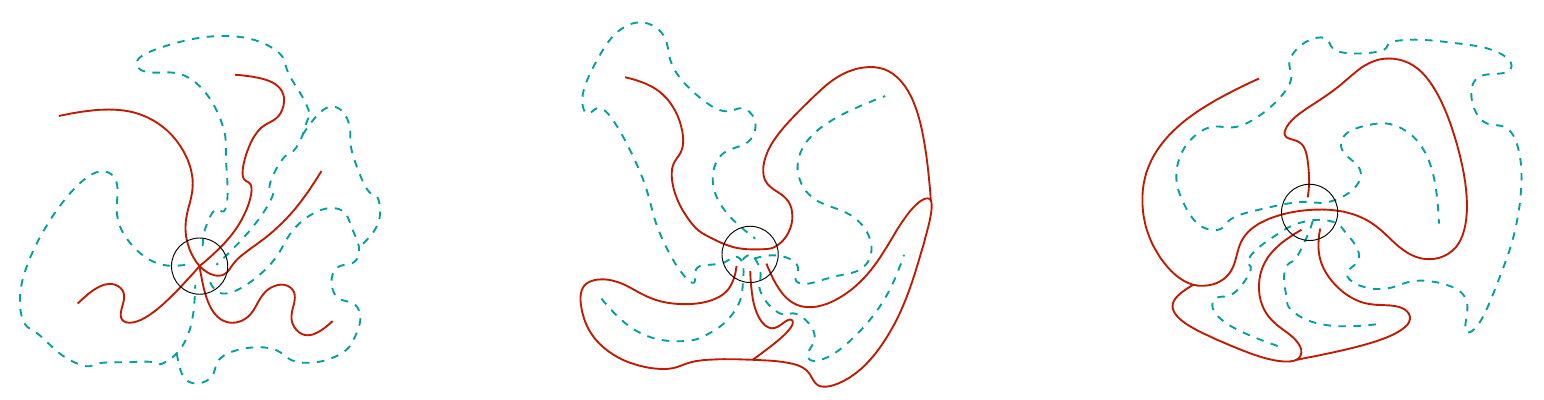}
}
\caption{Degree type $(5)$ and two examples of $(2,1,1,1)$ in a spanning tree of the plane, giving degree types $(1,1,1,1,1)$, $(4,1)$ and $(3,2)$ in a dual spanning tree.}
\label{f.degtype}
\end{figure}

The {\bf degree} of a point $x\in\hat\C$ in an immersed tree $f:\tau\lora\hat\C$, where we assume that $f$ is locally injective from each edge of $\tau$, is 
\be\label{e.deg}
\deg_{f}(x):=\sum_{f(v)=x}\deg_\tau(v)\,,
\ee
where the sum is over all points $v$ of $\tau$, meaning a vertex in $V(\tau)$ or a point on an edge in $E(\tau)$, and $\deg_\tau(v)$ is  the combinatorial degree in the first case, while equals 2 in the second case. (Note that for any immersed tree $f:\tau\lora\hat\C$ there exists a minimal $\tau^*\prec\tau$ on which $f$ can be defined, and here $f$ is equivalent to a locally injective immersion $f^*:\tau^*\lora\hat\C$. Hence the local injectivity assumption is not a real restriction.) For an
essential spanning forest $\For$,
\be\label{e.degF}
\deg_\For(x):=\sup_{\ell\geq 1} \sup_{f\in\For^{(\ell)}} \deg_f(x)\,.
\ee
The {\bf degree type} of a point $x$ in an immersed tree $f:\tau\lora\hat\C$ is the vector of summands in~(\ref{e.deg}), ordered in decreasing order, and the degree type in an essential spanning forest $\For$ is the supremum as in~(\ref{e.degF}), now w.r.t.~a natural partial order on the vectors of degree types: after padding vectors with zeros at the end, use the lexicographic ordering. The supremum in this partial order exists because of condition~(iii) of Definition~\ref{d.ESF}. (This supremum vector may start with a few $\infty$ entries if $\deg_\For(x)=\infty$, but this is fine.) See Figure~\ref{f.degtype} for a few examples (ignoring at this point the dual trees on the pictures).

For instance, saying that $x\in\C$ is a {\bf pinching point} for $\For$ if $\For^{(2)}$ includes a path which passes through $x$ twice without terminating there can be expressed as saying that $x$ has degree type at least $(2,2)$. If one of the two branches terminates at $x$, the other does not, i.e., degree type at least $(2,1)$, then we talk about a {\bf figure of 6}, while degree type at least $(1,1)$ is called a point of non-uniqueness, or a {\bf loop} at $x$. Points of degree type at least $(2)$ constitute the {\bf trunk} of $\For$: the union of curves in $\For^{(2)}$ excluding the endpoints. A {\bf branching point} is a point with degree type at least $(3)$.

\bl[Dual spanning tree]\label{l.dual}
There is a spanning tree $\MST_\infty^\dagger$ of $\C$ coupled with $\MST_\infty$ that is dual in the sense that none of its paths cross any of the paths of $\MST_\infty$, and whose distribution is again that of $\MST_\infty$.
\el

Note that we are not claiming that $\MST_\infty^\dagger$ is measurable w.r.t.~$\MST_\infty$, nor that there is a unique such spanning tree. These claims should be possible to prove, but we will not need them. For all subsequential scaling limits of the Uniform Spanning Tree on $\Z^2$, they were proved in \cite{Schramm} via first establishing that the trunk is a topological tree that is everywhere dense in $\C$, and then defining the dual tree in the complement of the trunk. 

\begin{figure}[htbp]
\centerline{
\includegraphics[width=0.35 \textwidth]{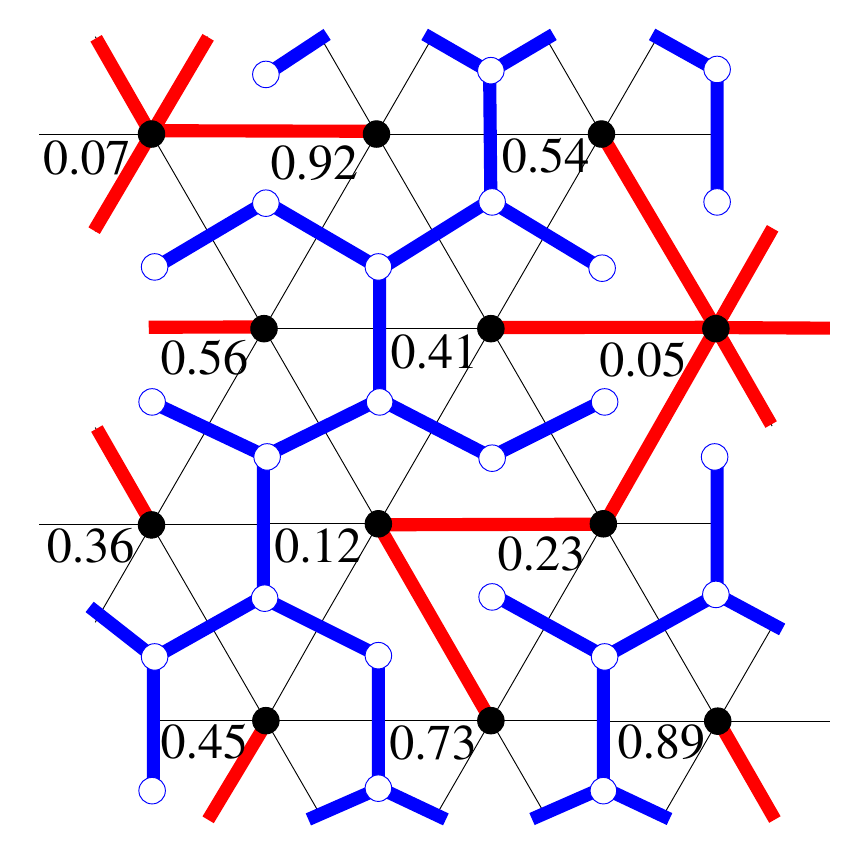}
\hskip 0.1 cm
\includegraphics[width=0.35 \textwidth]{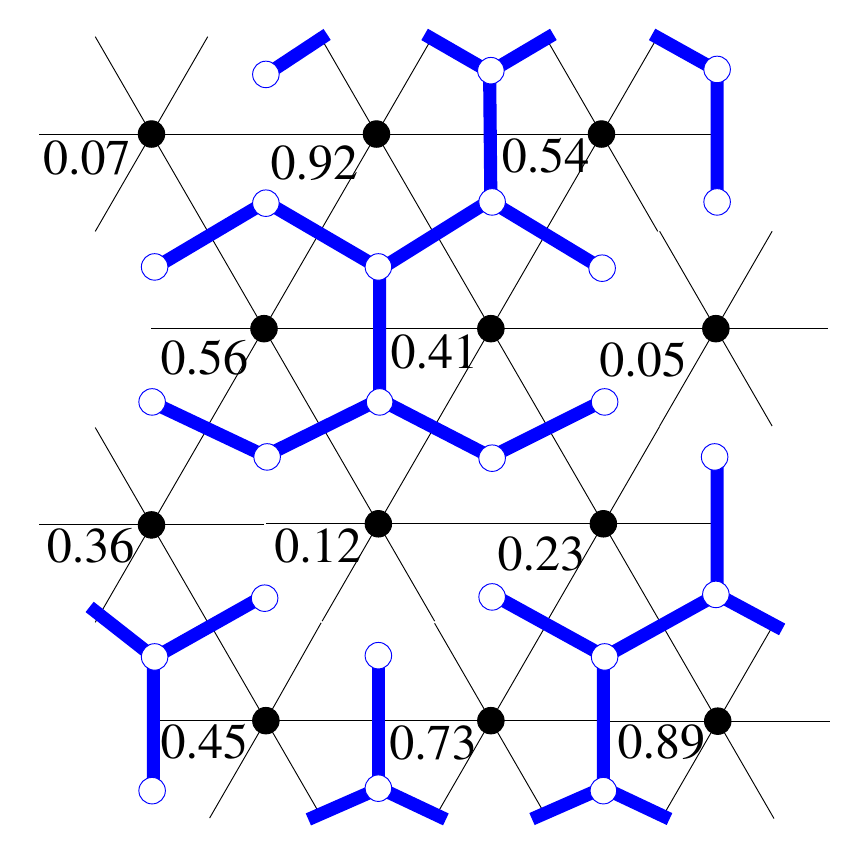}
\hskip 0.1 cm
\includegraphics[width=0.35 \textwidth]{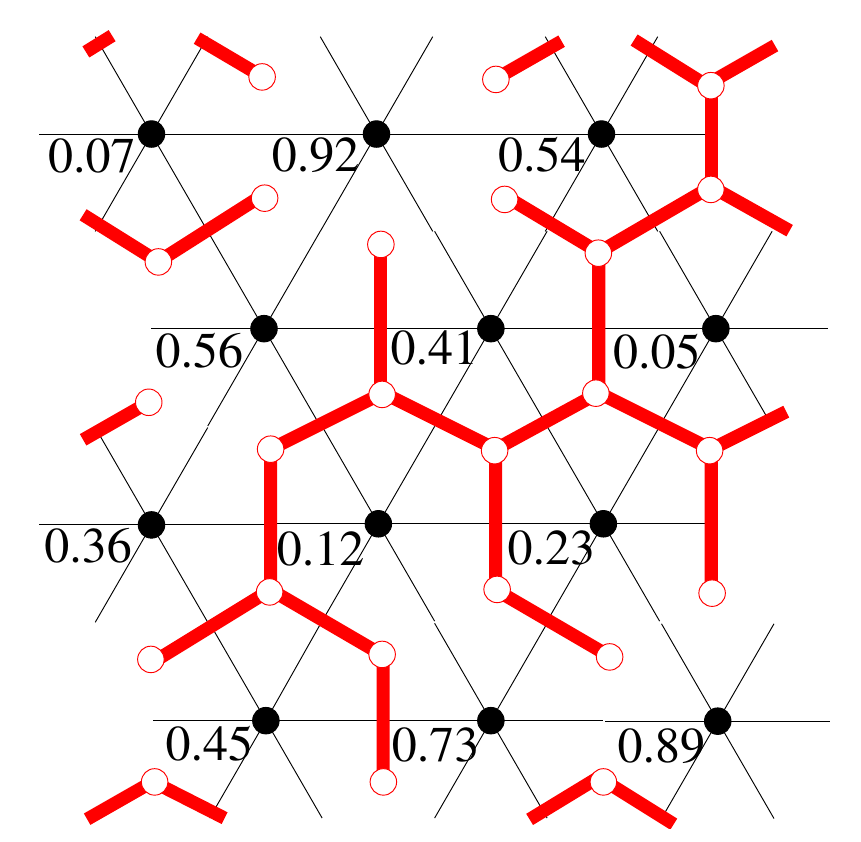}
}
\caption{(a) $\MST^M_\eta$ and its dual almost-tree ${\MST_\eta^M}^\dagger$. (b) The tree ${\mathsf{MaxST}_\eta^M}^*$ associated with the same vertex labels, which is just ${\MST_\eta^M}^\dagger$ minus two edges. (c) The tree ${\MST_\eta^M}^*$.}
\label{f.HexaTree}
\end{figure}

\proofof{Lemma~\ref{l.dual}} 
The planar dual of the triangular lattice $\Tg$ is the hexagonal lattice $\Td$, and, as usual, $\MST^M_\eta$ on $\eta\Tg \cap \Tor_M$ has a dual graph on $\eta\Td\cap\Tor_M$, denoted by ${\MST_\eta^M}^\dagger$. Because of the torus geometry, this dual has some cycles, but it is easy to check that for any null-homotopic cycle in $\eta\Td\cap\Tor_M$, the edge whose dual in $\eta\Tg \cap \Tor_M$ has the minimal weight must be present in $\MST_\eta^M$, hence must be missing from ${\MST_\eta^M}^\dagger$, and thus we are almost talking about the {\it Maximal} Spanning Tree on $\eta\Tg^*\cap\Tor_M$, denoted by ${\mathsf{MaxST}_\eta^M}^*$, which is defined from the Unif$[0,1]$ vertex labels $\{V(x)\}$ on $\Tg$ analogously to  $\MST^M_\eta$: each edge $e^* \in E(\eta\Td\cap\Tor_M)$ is the dual edge to some $e=(x,y) \in  E(\eta\Tg\cap\Tor_M)$, then we let
$$
U(e^*):=\big(V(x)\vee V(y), V(x)\wedge V(y)\big)\,,
$$ 
and get ${\MST_\eta^M}^*$ (or ${\mathsf{MaxST}_\eta^M}^*$) by removing the maximal (or minimal) edge from each cycle of $E(\eta\Td\cap\Tor_M)$ w.r.t.~the lexicographic ordering. See Figure~\ref{f.HexaTree}.

To make the connection between the macroscopic geometry of ${\MST_\eta^M}^\dagger$ and that of ${\mathsf{MaxST}_\eta^M}^*$ stronger, note that taking a pair of null-homotopic domains $\Lambda \subset \bar\Lambda$, the probability that all the paths of ${\MST_\eta^M}^\dagger$ connecting vertices in $\Lambda$ stay inside $\bar\Lambda$ is the same as in  ${\mathsf{MaxST}_\eta^M}^*$, and conditioning both measures on this event, the distribution of these paths agree. Of course,  ${\mathsf{MaxST}_\eta^M}^*$ has the same distribution as ${\MST^M_\eta}^*$ on $\eta\Tg^*\cap\Tor_M$, hence if we can understand the geometry of ${\MST^M_\eta}^*$, and prove, for instance, uniform quasi-locality, then we can use that, just as in Subsection~\ref{ss.full}, to understand the geometry of ${\MST_\eta^M}^\dagger$.

Now, we claim that, in the spirit of the remark after Figure~\ref{f.vertexTree}, the macroscopic structure of ${\MST^M_\eta}^*$ on $\eta\Tg^*\cap\Tor_M$ can be described using the near-critical ensemble on $\eta\Tg$. Indeed, if $x^*$ and $y^*$ are vertices in $\eta\Tg^*\cap\Tor_M$ such that the triangular faces they represent have vertices $x$ and $y$ in the same percolation $p$-cluster of $\eta\Tg\cap\Tor_M$, then it is immediate to see that the path in ${\MST^M_\eta}^*$ between $x^*$ and $y^*$ cannot cross any cycle of $\eta\Tg\cap\Tor_M$ whose vertices all have labels above $p$. In this sense, the path does not leave the $p$-cluster. It follows that two distinct $p$-clusters cannot be connected by two different paths of ${\MST^M_\eta}^*$ (we would otherwise get a cycle in ${\MST^M_\eta}^*$), and hence our entire paper applies to this version of the Minimal Spanning Tree. Thus we get that ${\MST^M_\eta}^*$ and ${\MST_\eta^M}^\dagger$ have the same unique scaling limit as $\eta\to 0$ then $M\to\infty$, denoted by $\MST_\infty^\dagger$, with the same distribution as $\MST_\infty$.

The fact that the paths of $\MST_\infty^\dagger$ do not cross the paths of $\MST_\infty$ is clear from obtaining them as scaling limits of discrete dual graphs.
\qed

It was proved in \cite{ABNW} that any subsequential limit of $\MST_\eta$ in $\hat\C$ is a spanning tree of $\hat\C$, and hence, using Lemma~\ref{l.dual}, it has one end (a single route to infinity). Furthermore, regularity properties of $\MST$ paths proved in that paper implied that the degrees in $\MST_\infty$ are almost surely bounded from above by some absolute deterministic constant $k_0\in\N$, and that the set of points with loops has Hausdorff dimension strictly between 1 and 2. It was also shown, using a Burton-Keane-type argument with trifurcation points and the amenability of the graph $\Z^2$ (see \cite{BurtonKeane} or \cite[Section 7.3]{LPbook}) that the set of branching points is at most countable. It was conjectured in \cite{ABNW} that there are no branching points of degree 4 or larger, and that there are no pinching points. We are now able to establish the latter conjecture, and get close to the former:

\bth[Degree types in $\MST$]\label{t.types} Almost surely in $\MST_\infty$ on $\C$:
\bi
\item[{\bf (i)}] there are no points of degree type at least $(2,2)$; in other words, for any two points $x,y\in\C$, none of the paths connecting the two vertices has a {pinching point}; 
\item[{\bf (ii)}]  there are no points of degree at least 5 (with any degree type);
\item[{\bf (iii)}]  the set of points of degree 4 (with any degree type) is at most countable.
\ei
These hold not only for the scaling limit on $\eta\Tg$ but also for any subsequential limit on $\eta\Z^2$.
\eth

\proof
{\bf (i)} We want to show that, in any given unit square in $\C$ and for any $0<\rho< 1$, the probability of $\MST_\infty$ having a point of degree type $(2,2)$ with the four strands not being connected within the radius $\rho$ ball of the pinching point is zero. (Note that if there was no positive radius in which the four strands are not connected, then this would in fact be a point of degree type at least 4, not (2,2). On the other hand, the strands must be connected somewhere in $\MST_\infty$, hence we can just assume that they are parts of one path connecting two vertices, which explains the second part of the statement.) For this, it is enough to show that for any $M>1$, the probability in $\MST_\eta=\MST^M_\eta$ that there is an $r$-square $B \in B^r([0,1]^2)$ (as in Definition~\ref{d.nested}, with $r<\rho$) such that there is a path $\gamma$ with two disjoint subpaths, $\gamma_1$ and $\gamma_2$, that both enter $B$ but are not connected to each other in $\MST_\eta$ within the $\rho$-neighbourhood of $B$ tends to 0 as $r\to 0$, uniformly in $\eta>0$.

Fix $\alpha>0$ arbitrarily small. As in the proof of Proposition~\ref{p.epsapprox}, we can take $\lambda<-1$, $\eps>0$, and $\lambda'>0$  such that with probability at least $1-\alpha/2$, all $\lambda$-clusters in $\Tor_M$ have diameter less than $\rho/10$, every point of $\Tor_M$ has in its $\rho/20$-neighborhood a ring of $\lambda$-clusters of diameter at least $\delta$ each, for some $0<\delta<\rho/20$ (uniformly in $\eta$), and all $\lambda$-clusters of diameter at least $\delta$ are connected in $\MST^{\vl,\eps}_\eta$, with these paths going through the same closed pivotals of $\PPP_\vl^\eps$ as the corresponding paths of $\MST_\eta$. We will assume that this event of probability at least $1-\alpha/2$ holds, and also that the above $r$-square $B$ exists, with some $r\ll\rho$ to be determined later.

\begin{figure}[htbp]
\SetLabels
(0.19*0.24)$x_1$\\
(0.8*0.17)$x_2$\\
(0.33*0.49)\red{$\gamma_1$}\\
(0.68*0.49)\red{$\gamma_2$}\\
(0.46*0.57)$r$\\
(0.2*0.65)$\rho$\\
(0.57*0.93)$<\rho/10$\\
(0.4*0.09)$>\delta$\\
(0.39*0.05)$<\rho/10$\\
(0.53*0.16)$<\rho/20$\\
\endSetLabels
\centerline{
\AffixLabels{
\includegraphics[width=0.6 \textwidth]{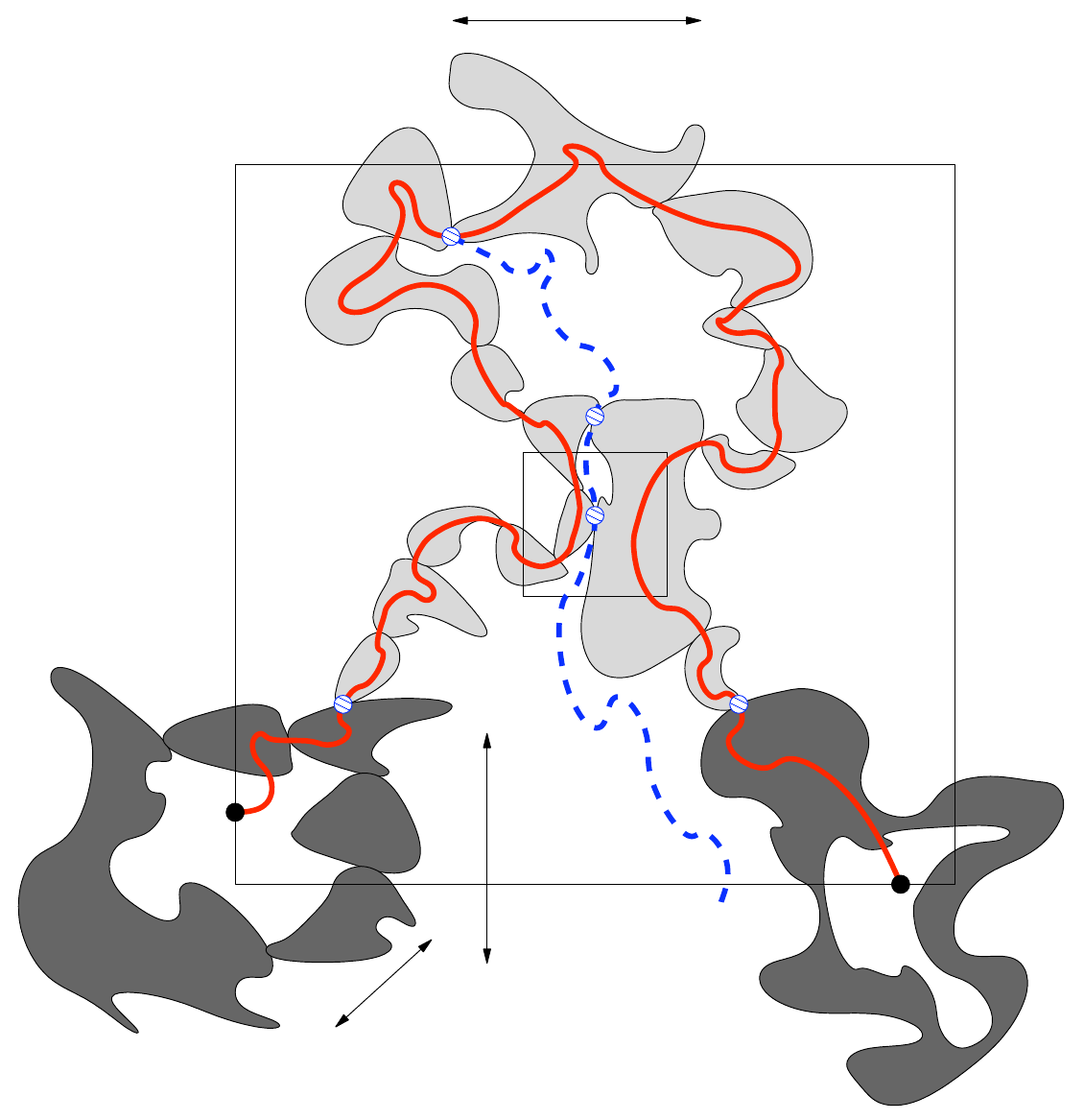}
}
}
\caption{Pinching would imply a near-critical 6-arm event.}
\label{f.pinching}
\end{figure}
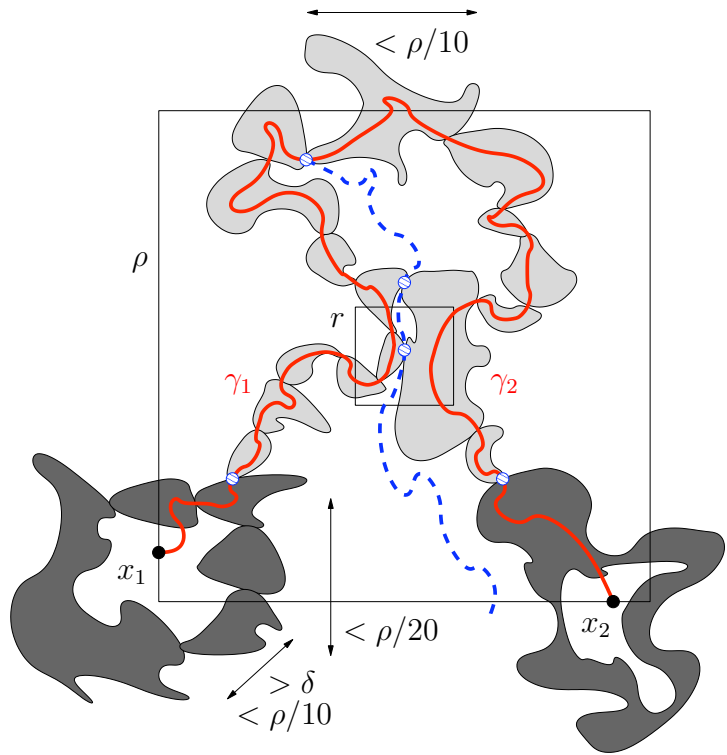

Any path in $\MST_\eta$ that connects two points in the same $\lambda$-cluster must stay in that cluster. Thus, the paths $\gamma_1$ and $\gamma_2$ that are connected in $\MST_\eta$ but not inside the $\rho$-neighbourhood of $B$ (denoted by $B_\rho$), must go through disjoint $\lambda$-clusters inside $B_\rho$. These $\lambda$-clusters all have diameter at most $\rho/10$, connected by $\lambda$-closed pivotals. Close to each end of each $\gamma_i$, there is such a $\lambda$-closed pivotal, at distance at least $\rho-\rho/5$ from $B$. Thus there must exist two $\lambda$-closed paths, separating the $\lambda$-clusters of $\gamma_1$ from those of $\gamma_2$, going through $B$, of radius at least $4\rho/5$. See Figure~\ref{f.pinching}.

We would like to bound now the labels from above on the $\MST_\eta$ paths. To this end,
let $x_i$ be the point where $\gamma_i$ leaves the $\rho$-neighborhood of $B$, at the end of $\gamma_i$ that is opposite from $\gamma_{3-i}$ along $\gamma$, for $i=1,2$. Around each $x_i$, there is a ring of macroscopic $\lambda$-clusters, the $\MST_\eta$ path from $x_1$ to $x_2$ must intersect at least one $\lambda$-cluster from each ring, and the part of the path connecting the two rings must go through $\lambda$-clusters connected by pivotals with labels at most $\lambda'$. Thus, besides the two $\lambda$-closed arms between radii $r$ and $4\rho/5$ we also have four $\lambda'$-open arms between the same radii. By the near-critical stability of 6-arm probabilities, Proposition~\ref{p.stab}, the probability of this happening anywhere in $\Tor_M$ is smaller than $\alpha/2$ if $r/\rho$ is chosen small enough. Therefore, the probability of the existence of $B$ is less than $\alpha$ if $r>0$ is chosen small enough, uniformly in the mesh $\eta>0$, and we are done.


{\bf (ii)} It is proved in \cite{BeffNol} that the critical monochromatic 5-arm exponent is strictly larger than the polychromatic one, which is 2 (see \cite[Corollary A.8]{SchSt}). Therefore, near-critical stability for the monochromatic 5-arm exponent (again, Proposition~\ref{p.stab}) tells us that no near-critical monochromatic 5-arm event between radii $r$ and $\rho$ happens anywhere in $[0,1]^2$ if $r/\rho$ is small enough. Based on this, as before, we will exclude the existence of an $r$-square $B\in B^r([0,1]^2)$ with degree 5 to distance at least $\rho$.

\begin{figure}[htbp]
\centerline{
\includegraphics[width=1 \textwidth]{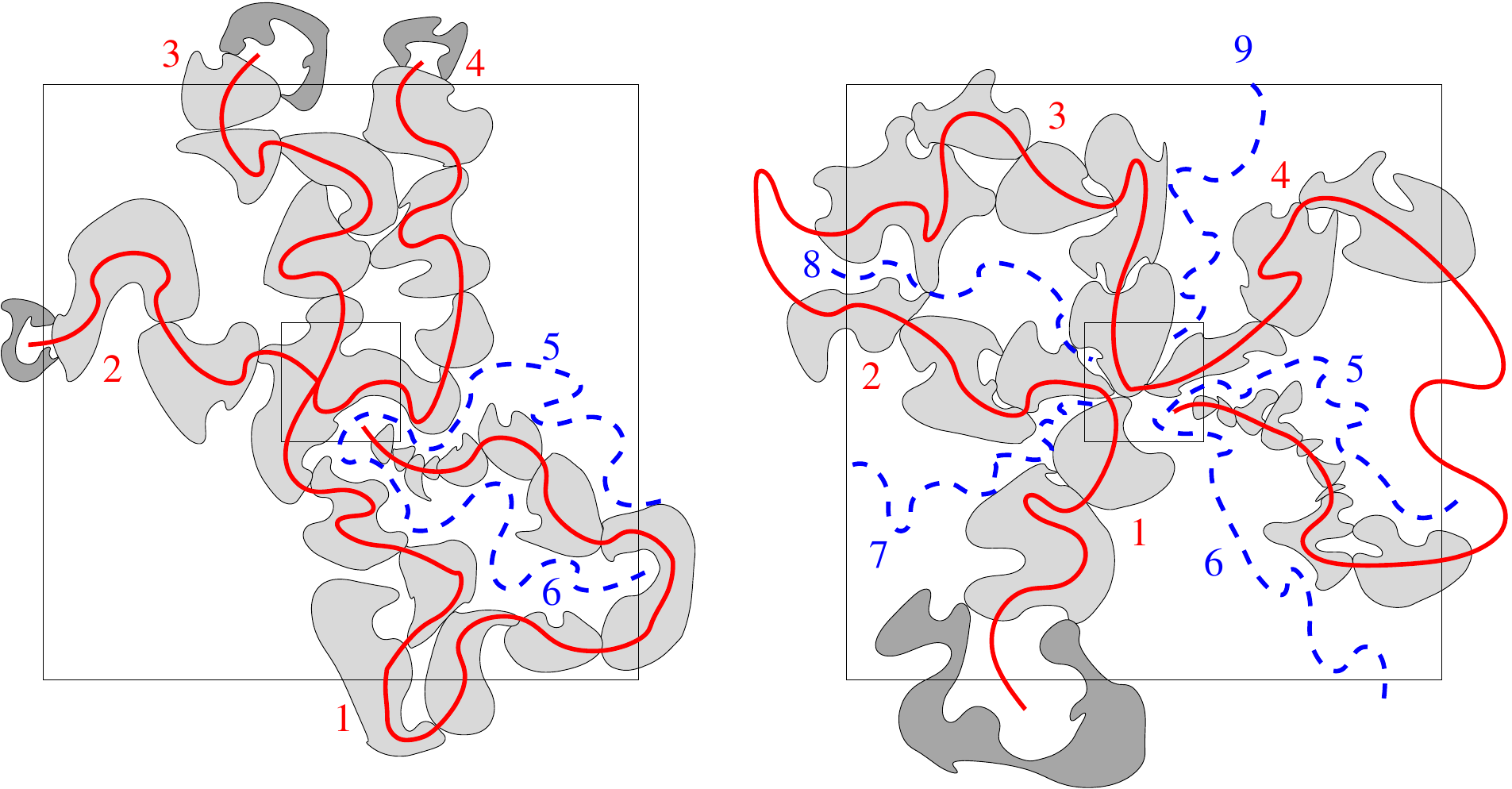}
}
\caption{Degree 5 would imply a near-critical monochromatic 5-arm or a polychromatic 6-arm event.}
\label{f.deg5}
\end{figure}

We look at the $\lambda$-clusters traversed by the five branches, for some small $\lambda<-1$. As in part~(i), the branches contributed by components at least 2 in the vector of the degree type traverse macroscopic $\lambda$-clusters, and hence the labels of their $\lambda$-closed pivotals are all at most some uniform $\lambda'$. That is, a degree component $k\geq 2$ implies $k$ $\lambda'$-open arms from $r$ to $\rho$ (e.g., the open arms labelled 1 to 4 on the left side of Figure~\ref{f.deg5}) and if there are more than one such components, we also have $\lambda$-closed arms separating them (e.g., the closed arms labelled 6 to 9 on the right side of Figure~\ref{f.deg5}). On the other hand, if we have $\ell \geq 1$ branches contributed by components of size 1 in the vector of the degree type, they are necessarily separated from the other branches by $\lambda$-closed paths. Thus, they
\begin{itemize}
\item either contribute $\ell+1$ closed arms to the existing 4 open arms,
\item or if we have at least two degree components with at least 2 branches, then they raise the number of closed arms by $\ell$,
\item or if we do not have any degree components with at least 2 branches, then $\ell\geq 5$, the degree type is all 1's, and we again get $\ell$ closed arms.
\end{itemize}
 Altogether, we either have at least 5 $\lambda'$-open arms, or at least 5 $\lambda$-closed arms, or a $(\lambda,\lambda')$ near-critical polychromatic event with at least 6 arms. None of these happens if $r/\rho$ is small enough, and we are done.

{\bf (iii)} Degree 4 points can have five different degree types: $(4)$, $(3,1)$, $(2,2)$, $(2,1,1)$, $(1,1,1,1)$. The countability of the first two types follows from the countability of branching points proved in \cite{ABNW}. Points of the third type do not exist, by part~(i) above. 
At a point of the fourth type, the dual $\MST_\infty$ tree defined in Lemma~\ref{l.dual} would either have a branching of degree 3, for which we already know countability, or a degree type $(2,2)$, which does not exist by part~(i).  (See Figure~\ref{f.degtype} for examples of dual degree types.) Finally, if a point has degree type $(1,1,1,1)$, then the dual tree has a branching point of degree 4 there, so we have countability again.

Since the well-known 5- and 6-arm bounds and Proposition~\ref{p.stab} hold also for $\Z^2$, all the above arguments work fine for subsequential limits of $\MST_\eta$ on $\eta\Z^2$, as well.
\qed

It is tempting to try and argue that a figure of 6 should imply 5 arms with labels bounded suitably by $\lambda$ and $\lambda'$, and hence by the near-critical stability of the 5-arm exponent (which is 2), the set of points with degree type $(2,1)$ should be at most countable, but we did not manage to make this argument work.

\subsection{A dimension bound for the trunk}\label{ss.trunk}

Our present techniques reveal very little about the dimension of different subsets of interest in $\MST_\infty$. It was proved in \cite{ABNW} that all the curves connecting any two points almost surely have Hausdorff dimension at least some unspecified deterministic $d_\mathrm{min}>1$ and at most another constant $d_\mathrm{max}<2$. Note that, having a countable number of branching points, the trunk is a countable union of such curves, hence we can equivalently talk about the dimension of the trunk. We will now slightly improve the upper bound to $d_\mathrm{max}=2-\alpha_2' < 7/4$, where $\alpha_2'$ is the monochromatic two-arm (or backbone) exponent of critical percolation, shown to be strictly larger than the polychromatic two-arm exponent $\alpha_2=1/4$ in \cite{BeffNol}. According to simulations, the true value of the Hausdorff dimension is close to $1.22^-$ \cite{WiWi, SwM}, while $2-\alpha_2'$ is close to $79/48 = 1.646^-$ \cite{BeffNolSimu}.

\bth\label{t.pathdim}
The Hausdorff dimension of the trunk of $\MST_\infty$ is almost surely at most $2-\alpha_2'<7/4$, where $\alpha_2'$ is the monochromatic two-arm exponent of critical percolation.
\eth

\proof
For any $\rho>0$, let $\Trunk^\rho_\eta$ (resp.~$\Trunk^\rho_\infty$) be the set of points in $[0,1]^2$ that have a path of $\MST_\eta$ (resp.~$\MST_\infty$) passing through them, going to distance at least $\rho$ in both directions. Since the trunk of $\MST_\infty$ is a countable union of sets of the form $\Trunk^\rho_\infty$, it is enough to prove the dimension bound on each $\Trunk^\rho_\infty$. Consider our usual grid $B^r([0,1]^2)$ of $r$-squares, with $r\ll\rho$. The subset of those $r$-squares that are intersected by $\Trunk^\rho_\eta$ (resp.~$\Trunk^\rho_\infty$) will be denoted by $\Trunk^{\rho,r}_\eta$ (resp.~$\Trunk^{\rho,r}_\infty$), and it is clear that in any coupling where $\MST_\eta$ converges to $\MST_\infty$ almost surely, for small enough $\eta>0$ we have $|\Trunk^{\rho,r}_\eta|/9 \leq |\Trunk^{\rho,r}_\infty| \leq 9|\Trunk^{\rho,r}_\eta|$ with probability close to 1, where the factors of 9 accommodate the possibility of the points of $\Trunk^\rho$ moving across the boundaries of $r$-squares. Therefore, it suffices to prove that for any $\beta>0$ there is a sequence $r_k\to 0$ such that
\be\label{e.trunkdr}
\Pb{|\Trunk^{\rho,r_k}_\eta| > r_k^{-2+\alpha_2'-\beta}} < 2^{-k}
\ee
for all small enough $\eta=\eta_k>0$, because then Borel-Cantelli gives that the Hausdorff (even the lower Minkowski) dimension of $\Trunk^\rho_\infty$ is almost surely at most $2-\alpha_2'+\beta$.

To prove (\ref{e.trunkdr}), take $\vl_k$ and $\eps_k$ such that with probability at least $1-3^{-k}$ all $\lambda_k$-clusters have diameter at most $\rho/10$, all points have a ring of $\lambda_k$-clusters of diameter at least $\delta>0$ in their $\rho/20$-neighborhood, and all $\lambda_k$-clusters of diameter at least $\delta$ are connected in $\MST^{\vl_k,\eps_k}_\eta$. Condition on this event, denoted by $\G_k$. Then, just as in the proof of Theorem~\ref{t.types}, every element of $\Trunk^{\rho,r}_\eta$ has a $\vl_k$-near-critical monochromatic 2-arm event from radius $r$ to $\delta/2$. From near-critical stability, we know that, for any $B\in B^r([0,1]^2)$, denoting this 2-arm event by $\A_2'(B,r,\delta/2,\vl_k)$, we have 
$$
\Pb{ \A_2'(B,r,\delta/2,\vl_k) }< C_{\delta,k} \, r^{\alpha_2'}\,.
$$
Since $1/(1-3^{-k})<2$, the previous line gives
$$
\Pb{ \A_2'(B,r,\delta/2,\vl_k) \md \G_k } < 2 C_{\delta,k} \, r^{\alpha_2'}\,,
$$
and, summing up over $B$,
$$
\Eb{ |\Trunk^{\rho,r}_\eta| \md \G_k } < 2C_{\delta,k} \, r^{-2+\alpha'_2}\,.
$$
Then, by Markov's inequality, for any $\beta>0$,
$$
\Pb{ |\Trunk^{\rho,r}_\eta| > 2C_{\delta,k} \, r^{-2+\alpha_2'-\beta/2} \md \G_k } < r^{\beta/2}\,.
$$
By taking $r_k>0$ so small that $2C_{\delta,k} < r_k^{-\beta/2}$ and $r_k^{\beta/2} < 3^{-k}$, we get that
$$
\Pb{|\Trunk^{\rho,r_k}_\eta| > r_k^{-2+\alpha_2'-\beta} \md \G_k } < 3^{-k}\,.
$$
Since we have 
\begin{align*}
\Pb{|\Trunk^{\rho,r_k}_\eta| > r_k^{-2+\alpha_2'-\beta} } &< \Pb{|\Trunk^{\rho,r_k}_\eta| > r_k^{-2+\alpha_2'-\beta} \md \G_k } +\Pb{\G_k^c}\\
& < 3^{-k} + 3^{-k}\,,
\end{align*}
we have verified~(\ref{e.trunkdr}) and completed the proof.
\qed

\section{Invasion percolation}\label{s.IP}

The Invasion Tree in a finite graph is simply the $\MST$ itself, hence it cannot provide us with a good finite approximation to $\IP$ in the infinite plane. Instead, 
we will consider the following finite versions:
\bi
\item $\IP_\eta^{M,\p}$ will be the tree built by the invasion process started from the origin, stopped at the first time that it reaches $\p [-M,M]^2$.
\item For a fixed vertex $x\in V(\eta\Tg)$ and $M$ large enough so that $x\in [-M,M]^2$, we will denote by $\IP_\eta^{M,x}$ the invasion process in the torus $\Tor_M$, started from the origin, stopped at the first time when it reaches $x$.
\ei
When $M\to\infty$, the weak limits of the above measures are $\IP_\eta=\IP_\eta(0)$ and $\IP_\eta(0) \cup \IP_\eta(x)$, respectively. Of course, the latter coincides with $\IP_\eta(0)$ with positive probability, and $\IP_\eta(0) \,\triangle\, \IP_\eta(x)$ is almost surely finite. These results are classical \cite{ChChN1, AleMol, Alexander,LPS}.

Given the enhanced $\eps$-networks $\EnNet^{\vl,\eps}_\eta$ and $\EnNet^{\vl,\eps}_\infty$ defined in Proposition~\ref{p.ENetconv}, the cut-off versions of the above invasion trees, both in the discrete case and in the continuum, can be defined quite similarly to $\MSF^{\vl,\eps}$ (done in Definition~\ref{d.MSF}) and $\MST^{\vl,\eps}$ (in Lemma~\ref{l.giant}):

\bde[The cut-off invasion trees $\IP^{\vl,\eps,s}$ in $\Tor_M$, with target set $\p$ or $x$]\label{d.epsIP}\ 
\bnum
\item Consider the edge-labelled graph defined in steps 1-3 of Definition~\ref{d.MSF} on the set of the primal routers of $\EnNet^{\vl,\eps}$ as vertices.
\item Take its giant component, which exists with large probability by Lemma~\ref{l.giant}. On the bad event that this giant cannot be defined, our cut-off invasion trees will be just degenerate one-point trees. 
\item Take the router closest to the origin $0$; in case of a tie, decide in some arbitrary but fixed manner. This will be called the origin router. Furthermore, consider all routers that are at most distance $s>0$ from the target set $\p [-M,M]^2$ or $x$. By Lemma~\ref{l.giant}, for any $s>0$, if $\lambda$ is very negative, $\lambda'$ is very positive, and $\eps$ is small, then with high probability the set of these target  routers is not empty. When it is empty, the invasion tree will consist of just the origin router.
\item Take the invasion tree process in the above graph, started from the origin router, stopped when reaching any of the target routers. There may be steps in the invasion process when more than one edge with label $\lambda$ lead out of the invaded set; in such a case, all these edges get invaded simultaneously.
\enum
\ede

It was proved already in \cite{ABNW} that the set of points with degree larger than 1 (i.e., points in the trunk or having a loop) in any subsequential scaling limit of $\MST_\eta$ is of zero measure. Therefore, almost surely there is a unique path of $\MST_\infty$ that goes to the origin, and hence we did not lose any information in the above definition by taking the router closest to the origin instead of considering all routers that are $s$-close to it.

Given this definition, we immediately have the following analogues of Corollary~\ref{cor.MSTABNW} and Proposition~\ref{p.epsapprox}. Note the double meaning of the parameter $s$: if we want to reach precision $s>0$ in $d_{\Omega_M}$, it is enough to get $s$-close to the target sets.

\bl\label{l.IPABNW}
For any $M>0$, target set $\p$ or $x\in\Tor_M$, and any $s,\alpha>0$, if $\lambda<-1$ is very negative, $\eps>0$ is small, and $\lambda'>1$ is large enough, then, in the coupling of Proposition~\ref{p.ENetconv}~(ii) between $(\omega^\lambda_\eta,\PPP^\eps_\vl)$ and $(\omega^\lambda_\infty,\PPP^\eps_\vl)$, for all $\eta>0$ small enough,
$$
\Pb{ d_{\Omega_M}(\IP^{\vl,\eps,s}_\eta,\IP^{\vl,\eps,s}_\infty) < s} > 1-\alpha\,.
$$
\el

\bl\label{l.epsapprox}
For any $M>0$, target set $\p$ or $x\in\Tor_M$, and $s,\alpha>0$, if $\lambda<-1$ is very negative, $\eps>0$ is small, and $\lambda'>1$ is large enough, then, for all $\eta>0$ small enough,
$$
\Pb{ d_{\Omega_M}(\IP_\eta,\IP^{\vl,\eps,s}_\eta) < s} > 1-\alpha\,,
$$
where, of course, $\IP_\eta$ is only a shorthand now for $\IP^{M,\p}_\eta$ or $\IP^{M,x}_\eta$.
\el

Using these lemmas, the proof of the following theorem follows exactly the proofs of Theorem~\ref{t.torus} and Theorem~\ref{t.main}.

\bth\label{t.IP} 
For any $M>0$, the invasion trees $\IP_\eta^{M,\p}$ and $\IP_\eta^{M,x}$ started at the origin of $\eta\Tg\cap\Tor_M$ converge in  distribution as $\eta\to 0$, in the metric $d_{\Omega_M}$ of Definition~\ref{d.ESF}, to the unique scaling limits $\IP_\infty^{M,\p}$ and $\IP_\infty^{M,x}$, respectively.

The invasion tree $\IP_\eta$ started at the origin of $\eta\Tg$ converges in distribution to a unique scaling limit $\IP_\infty$ that is invariant under scalings and rotations. 

As $M\to\infty$, the weak limit of $\IP_\infty^{M,\p}$ is $\IP_\infty$ and the weak limit of $\IP_\infty^{M,x}$ is $\IP_\infty(0)\cup\IP_\infty(x)$.
\eth


\section{Questions and conjectures}\label{s.conj}

We start with a very natural and interesting open problem:

\bcj\ 
\bi
\item[{\bf (i)}]  Show that $\MST_\infty$ is not conformally invariant. In particular, show that it is different from the scaling limit of the Uniform Spanning Tree, described in \cite{LSWUST}. 
\item[{\bf (ii)}] Show that $\IP_\infty$ is not conformally invariant.
\ei
\ecj

This is of course supported by simulation results \cite{RGB}. Moreover, it was explained in \cite{MST} why our description of these scaling limits using the near-critical ensemble gives serious support to this conjecture, and why it is nevertheless not at all an easy issue. The case of $\IP_\infty$ might be simpler, using the results of \cite{Damron1}.

Probably the simplest open problem in this section is the following one, left open by Lemma~\ref{l.dual}:

\bcj
Show that there is a unique dual tree $\MST_\infty^\dagger$, measurable w.r.t.~$\MST_\infty$.
\ecj

The following questions are left open by Theorem~\ref{t.types}:

\bq[Topology of $\MST_\infty$]\ 
\bi
\item[{\bf (i)}] Are there non-simple paths giving figures of 6, i.e., points with degree type $(2,1)$? 
\item[{\bf (ii)}] Show that almost surely there are no points of degree 4.
\ei
\eq
 
Finally, sharpening the bound of Theorem~\ref{t.pathdim} would probably require new techniques:

\bq Find the Hausdorff dimension of the paths of $\MST_\infty$.
\eq

\bibliographystyle{alpha}
\addcontentsline{toc}{section}{References}

\def\arXiv#1#2{\href{http://front.math.ucdavis.edu/#1}{{\tt arXiv:#1 [#2]}}} 
\def\arXivo#1{\href{http://front.math.ucdavis.edu/#1}{{\tt [arXiv:#1]}}} 



\ \\

\noindent{\bf Christophe Garban}\\
Institut Camille Jordan, Universit\'e Lyon 1\\
\url{http://math.univ-lyon1.fr/~garban}\\
\\
{\bf G\'abor Pete}\\
Alfr\'ed R\'enyi Institute of Mathematics, Hungarian Academy of Sciences, Budapest,\\
and Institute of Mathematics, Budapest University of Technology and Economics\\
\url{http://www.math.bme.hu/~gabor}\\
\\
{\bf Oded Schramm} (December 10, 1961 -- September 1, 2008)\\
Microsoft Research\\ 
\url{http://research.microsoft.com/en-us/um/people/schramm/}

\end{document}